\documentclass[11pt]{article}
\usepackage{arxiv}

\usepackage{graphicx}
\usepackage{float} 
\usepackage{mathrsfs}
\usepackage{amsmath,amssymb,amsthm}
\usepackage{algorithm}
\usepackage{algpseudocode}
\usepackage{caption}
\usepackage{subcaption}
\usepackage{mathtools}
\usepackage{multicol} 
\usepackage{booktabs}
\usepackage{xcolor}       
\usepackage[T1]{fontenc}    
\usepackage[colorlinks=true,linkcolor=blue,citecolor=blue]{hyperref}
\usepackage{cleveref} 
\usepackage{url} 

\newtheorem{remark}{Remark}

\def\bv{\boldsymbol{v}}
\def\by{\boldsymbol{y}}
\def\bz{\boldsymbol{z}}
\def\bw{\boldsymbol{w}}
\def\bphi{\boldsymbol{\phi}}
\def\bmu{\boldsymbol{\mu}}
\def\bsigma{\boldsymbol{\sigma}}
\def\bepsilon{\boldsymbol{\epsilon}}
\DeclareMathOperator*{\argmin}{argmin} 
\DeclareMathOperator*{\argmax}{argmax} 
\def\bcM{\boldsymbol{\mathcal{M}}}
\def\bcV{\boldsymbol{\mathcal{V}}}
\def\bcY{\boldsymbol{\mathcal{Y}}}
\def\bcW{\boldsymbol{\mathcal{W}}}
\def\b{\textbf}

\begin{document}

\title{InVAErt networks for amortized inference and identifiability analysis of lumped parameter hemodynamic models}
\author{
Guoxiang Grayson Tong$^{1}$, Carlos A. Sing Long$^{2}$ and Daniele E. Schiavazzi$^{1*}$ \\
$^{1}$Department of Applied and Computational Mathematics and Statistics, University of Notre Dame, \\
Notre Dame, 46556, IN, United States\\
$^{2}$Institute of Mathematical and Computational Engineering, Pontificia Universidad Católica de Chile, \\
Santiago, Chile}
\date{}

\maketitle



\begin{abstract}
Estimation of cardiovascular model parameters from electronic health records (EHR) poses a significant challenge primarily due to lack of identifiability.
Structural non-identifiability arises when a manifold in the space of parameters is mapped to a common output, while practical non-identifiability can result due to limited data, model misspecification, or noise corruption.
To address the resulting ill-posed inverse problem, optimization-based or Bayesian inference approaches typically use regularization, thereby limiting the possibility of discovering multiple solutions.
In this study, we use inVAErt networks, a neural network-based, data-driven framework for enhanced digital twin analysis of stiff dynamical systems.
We demonstrate the flexibility and effectiveness of inVAErt networks in the context of physiological inversion of a six-compartment lumped parameter hemodynamic model from synthetic data to real data with missing components.
\end{abstract}

\section{Introduction}\label{sec:introduction}

Lumped parameter (a.k.a. zero-dimensional or 0D) hemodynamic models provide computationally inexpensive representations of the human cardiovascular system and are widely used either as low-fidelity surrogates or to provide closed-form boundary conditions in 3D cardiovascular simulations~\cite{shi2011review,  brown2024modular, kim2010patient, tran2017automated}.
Such representations model the evolution of bulk hemodynamic quantities of interest (QoI), simulating the time evolution of organ-level pressure, volume and flow~\cite{shi2011review, pfaller2022automated}, disregarding any spatial dependence. They can be represented as circuit models arranged by compartments, and formulated as systems of ODEs whose solution is computed using numerical solvers.

Despite their simplicity, working with 0D models can be challenging. First, these models are often described by stiff ODE systems~\cite{davis1991teaching, marquis2018practical}, requiring the careful selection of the numerical algorithm and time step size. Second, due to their reduced computational cost, these systems tend to be overparametrized, which typically induces structural non-identifiability~\cite{wieland2021structural,raue2009structural} in inverse problems aiming to assess cardiovascular function from the observed physiological response~\cite{schiavazzi2017patient, harrod2021predictive, pironet2016structural, pironet2019practical, marquis2018practical}. Additionally, practical non-identifiability may also manifest as a result of poor data quality, model misspecification and aleatoric uncertainty. 
Lack of identifiability can induce solution non-uniqueness, excessive uncertainty, and convergence failure in numerical algorithms~\cite{wieland2021structural,raue2013joining}.

Parameter estimation, sensitivity, and identifiability analysis of lumped parameter hemodynamic models are extensively studied in the literature. An incomplete list of contributions include the analysis of structural identifiability in a six-compartment 0D cardiovascular model in~\cite{pironet2016structural}, which shows that global identifiability can be achieved by adjusting the clinical targets and modifying the state equations accordingly. The method of profile likelihood~\cite{raue2009structural} is applied to investigate practical non-identifiability in a three-compartment hemodynamic model~\cite{pironet2019practical}, where likelihood-based confidence intervals indicate whether a parameter is practically identifiable. In~\cite{tran2017automated, harrod2021predictive, schiavazzi2017patient,marquis2018practical}, Bayesian inference and an eigenanalysis of the Fisher information matrix are used jointly for model calibration, uncertainty quantification and local identifiability analysis in the context of lumped parameter hemodynamics. Finally, a recent neural network-based approach for simulation-based inference (SBI~\cite{cranmer2020frontier}) has been explored in~\cite{wehenkel2023simulation}, where Neural Posterior Estimation (NPE) is applied to a one-dimensional whole-body cardiovascular model.

In this paper, we perform \emph{model synthesis} for a six-compartment lumped parameter hemodynamic model~\cite{davis1991teaching} using the recently proposed framework of inVAErt networks~\cite{tong2024invaert}.
InVAErt networks are designed to learn \emph{enhanced} digital twin representations, which extend simple emulation with the ability to solve amortized inverse problems and to assess indentifiability.
To handle ill-posed inverse problems with multiple solutions, a variational encoder learns an input-dependent latent space, which is used to restore bijectivity in the input-to-output map. Decoding the combination of latent space realizations and observations allows one to obtain entire manifolds of solutions. This is in sharp contrast with classical approaches based on regularization (see, e.g., ~\cite{vogel2002computational}), that enforce uniqueness through augmented penalized objectives, preventing the possibility of characterizing solution multiplicity. 

Bayesian methods can potentially find multiple solutions through posterior analysis. However, they may encounter challenges due to their heavy reliance on prior distributions, the efficiency of MCMC-type samplers in high-dimensional problems, and poor mixing in non-identifiable models~\cite{raue2013joining,wieland2021structural}.

While previous literature on inVAErt networks related to numerical examples from linear and non linear maps, dynamical systems, and spatio-temporal PDEs, this paper demonstrates extensions to applications involving noisy and missing data. 
Noise is injected during training to enhance generalization as well as the ability of inVAErt networks to deal with \emph{practical} identifiability. In addition, we also leverage flow-based density estimation to perform inversion tasks with missing data, a typical feature of electronic health record (EHR) datasets.

Here we briefly summarize our main contributions:
\begin{enumerate}
    \item We perform data-driven model synthesis and generate an \emph{enhanced digital twin} for a lumped parameter hemodynamic system (CVSim-6) using inVAErt networks.
    \item We analyze in detail the stiffness of the system, showing that inVAErt networks can be used for model synthesis of stiff ODE systems.
    \item We show the ability of an inVAErt network to characterize high-dimensional non-identifiable manifolds for complex systems of ODEs, and suggest a visualization technique that offers an improved assessment of the system identifiability.
    \item We demonstrate the ability of inVAErt networks to handle practical non-identifiability and to work with real data from an EHR dataset with missing attributes. 
\end{enumerate}

The content of this paper is organized as follows: In Section~\ref{sec: statement}, we briefly state the problem of interest and introduce the mathematical notations.
The lumped-parameter hemodynamic model, the CVSim-6 system, is introduced in Section~\ref{sec:cvsim6}, with its associated ODE stiffness analysis covered in Section~\ref{sec: cvsim6-stiffness}.
Section~\ref{sec:inVAErt} discusses the architecture of our proposed neural networks, optimization tasks, and training details.
Two inference tasks are discussed in Section~\ref{sec: cvsim6-str} and Section~\ref{sec: cvsim6-ehr}.
The first task focuses on using noiseless synthetic data to explore the structural identifiability of the CVSim-6 system, including the study of non-identifiable manifolds (Section~\ref{sec: cvsim6-str-manifold}), and missing data analysis (Section~\ref{sec: cvsim6-str-missing}).
The second task deals with real-world clinical measurements from an EHR dataset, and show how the addition of training noise can effectively improve performance in inference tasks (Section~\ref{sec: cvsim6-ehr}).
For the interested reader, our dataset and code can be accessed at \url{https://github.com/desResLab/InVAErt4Cardio}.

\subsection{Background and notation}\label{sec: statement}

The input-to-output (forward) map for an hypothetical cardiovascular model is denoted as $f: \bcV \mapsto \bcY$, mapping a vector of parameters $\bv \in \bcV \subset \mathbb{R}^{\dim(\bv)}$ to a set of clinical targets $\by \in \bcY \subset \mathbb{R}^{\dim(\by)}$, with $\bcV$ and $\bcY$ representing the abstract input and output spaces, respectively.
In this context, the operator $f$ involves both solving a system of ODEs and post-processing quantities of clinical interest from the time-dependent solution, e.g., extracting minimum, maximum (diastolic and systolic) or time-average values over the cardiac cycle, computing acceleration times, etc.
The input $\bv$ provides a zero-dimensional, simplified characterization of the physiology through the hydrodynamic analogy~\cite{milivsic2004analysis}, relating major viscous losses to Ohm resistance, vascular compliance to electric capacitance and blood flow inertia to inductance. Despite the simplicity of this approximation, varying $\bv$ can span a broad spectrum of physiological conditions from healthy to disease~\cite{abdi2015lumped,harrod2021predictive}.

Determining lumped parameters point estimates from observations (clinical data) is possible by solving the inverse problem $\bv = h(\by)$, but, in practice a number of challenges are associated with this task. 
First, the inverse or pseudo-inverse operator $h: \bcY \mapsto \bcV$ is often ill-posed due to non-uniqueness in structurally non-identifiable models. Specifically, for a given output $\by$, there may exist $\bcM_{\by} \subset \bcV$ such that $f(\bv) = \by$ for all $\bv \in \bcM_{\by}$. We call $\bcM_{\by}$ the \emph{non-identifiable manifold}, embedded in $\bcV$ and associated with the observation $\by$, where the complexity in determining this manifold depends on the nonlinearity of the forward map $f$, the properties of its gradient, and the dimensionality mismatch between $\bv$ and $\by$~\cite{tong2024invaert}.
The existence of $\bcM_{\by}$ is related to redundancy in the parameterization of the map $f$ (or, in other words, the map $f$ is not injective), as changes in one component of the input $\bv$ can be entirely compensated by adjusting the remaining components, while leaving the output $\by$ unaffected~\cite{wieland2021structural,raue2009structural,kreutz2018easy}.

Structural non-identifiability typically occurs in \emph{partially observed} systems (or observed under the so-called \emph{small data} regime) where $\dim({\by}) < \dim(\bv)$~\cite{wieland2021structural}, and therefore there is not enough information to uniquely determine the input parameters. 
By design, InVAErt networks~\cite{tong2024invaert} leverage input-dependent \emph{latent variables} $\bw\in\bcW$ to compensate for this information gap (see Section~\ref{sec:inVAErt}), where, generally, a higher dimensional $\bcW$ increases the ability to capture complicated manifolds $\bcM_{\by}$.

Second, ill-posedness of an inverse problem can also arise due to data scarcity, model misspecification or presence of noise. All these factors lead to the violation of Hadamard's conditions for well-posedness, that is, the existence of a unique solution that is stable~\cite{hadamard2014lectures,kirsch2011introduction,latz2020well}. In particular, real-life clinical data may not be adequately characterized by a simplified mathematical model or, formally, real measurements may be such that $\by \not\in \textrm{Range}(f)$. This may be due to random noise originating from the finite precision of sensors and acquisition devices used in the clinic, or even human errors in reporting the data, which can push data outside the range of $f$ significantly,  increasing the complexity of inverse problems. Contrary to structural non-identifiability due to the inherent model parameterization, the above issues typically occur when dealing with real data, hence the name of \emph{practical} non-identifiability~\cite{raue2009structural}.
%

Finally, input training data is generated from an assumed prior $\bv \sim p(\bv)$ which induces a distribution in the outputs $f(\bv) = \by \sim p(\by)$. Sampling from and evaluating the output distribution $p(\by)$ are of critical importance when modeling physics-based systems. Therefore inVAErt networks, beside learning the forward $\bv \mapsto \by$ and inverse $\by \mapsto \bv$ maps, are also equipped with a flow-based density estimator (See Section~\ref{sec: cvsim6-str-missing}, Section~\ref{sec: cvsim6-ehr}).
%

\section{Methods}\label{sec:method}

\subsection{The CVSim-6 hemodynamic model}\label{sec:cvsim6}

\begin{figure}[!ht]
    \centering
    \includegraphics[scale=0.35]{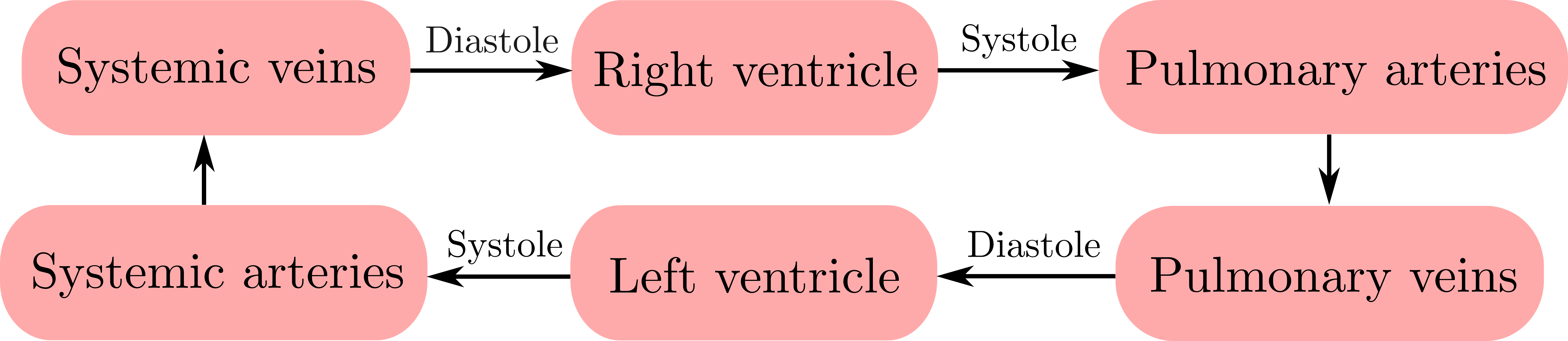}
    \caption{Schematic of compartments in the CVSim-6 model.}
    \label{fig:cvsim6-diag}
\end{figure}

CVSim-6 is a six-compartment 0D model, originally proposed for teaching cardiovascular physiology using computers~\cite{davis1991teaching}. Additional studies related to the use of CVSim-6 for educational tool development, numerical simulation and parameter estimation can be found in~\cite{heldt2010cvsim, mukkamala2000forward, samar2005cardiovascular, su2023privacy}.
It simulates the evolution of blood pressures, flows and volumes through the left and right ventricular chambers, systemic arteries, systemic veins, pulmonary arteries and pulmonary veins as shown in Figure~\ref{fig:cvsim6-diag}. The subscripts $(\cdot)_l,(\cdot)_r,(\cdot)_a,(\cdot)_v,(\cdot)_{pa},(\cdot)_{pv}$ refer to quantities specific to each compartment.

The input parameters $\bv \in \mathbb{R}^{23}$ are reported in Table~\ref{table:cvsim6 input paras}, with reference values listed from~\cite{davis1991teaching}.
Each of the six compartments is characterized by a capacitance (compliance) $C$, a resistance $R$ and an unstressed (zero-pressure filling) volume $V^0$, except for the two ventricles, where prescribed systolic and diastolic capacitances are responsible for ventricular contraction and relaxation. Models for the atrial chamber are not included in CVSim-6. 
Consistent with the compartmental nature of the model and its formulation in terms of a close-loop autonomous system of ODEs, we follow the convention where the outflow resistance of one compartment coincides with the inflow resistance of the subsequent compartment.
Additional parameters include the heart rate $Hr$, transthoracic pressure $P_{th}$ (disregarding any circulatory-respiratory coupling), and the systolic fraction $r_{sys}$ of a heart cycle.

\begin{table}[ht!]
{\small
\begin{center}
\begin{tabular}{@{} l l l l @{}}
\toprule
Num. & Description & Ref. & Unit\\
\midrule
1. & Heart rate ($Hr$) & 72.00 & (bpm)\\
2. & Transthoracic pressure ($P_{th}$) & -4.00 & (mmHg)\\
3. & Systolic ratio per heart cycle ($r_{sys}$) & 0.33 & $-$\\ 
4. & Left ventricular diastolic capacitance ($C_{l,dia}$) & $7.50\cdot 10^{-3}$ & (mL/Barye)\\ 
5. & Left ventricular systolic capacitance ($C_{l,sys}$) & $3.00\cdot 10^{-4}$ & (mL/Barye)\\ 
6. & Arterial capacitance ($C_a$) & $1.20\cdot 10^{-3}$  & (mL/Barye)\\ 
7. & Venous capacitance ($C_v$) & $7.50\cdot 10^{-2}$   & (mL/Barye)\\ 
8. & Right ventricular diastolic capacitance ($C_{r,dia}$) & $1.50\cdot 10^{-2}$ & (mL/Barye)\\ 
9. & Right ventricular systolic capacitance ($C_{r,sys}$) & $9.00\cdot 10^{-4}$ & (mL/Barye)\\ 
10. & Pulmonary arterial capacitance ($C_{pa}$) & $3.23\cdot10^{-3}$ & (mL/Barye)\\ 
11. & Pulmonary venous capacitance ($C_{pv}$) & $6.30\cdot 10^{-3}$ & (mL/Barye)\\ 
12. & Left ventricular input resistance ($R_{l,in}$)  & 13.33  & (Barye$\cdot$s/mL)\\ 
13. & Left ventricular output resistance ($R_{l,out}$) & 8.00 & (Barye$\cdot$s/mL)\\  
14. & Arterial resistance ($R_a$) & 1333.22  &(Barye$\cdot$s/mL)\\ 
15. & Right ventricular input resistance ($R_{r,in}$)  & 66.66  & (Barye$\cdot$s/mL)\\ 
16. & Right ventricular output resistance ($R_{r,out}$) & 4.00  &(Barye$\cdot$s/mL)\\  
17. & Pulmonary venous resistance ($R_{pv}$) & 106.66 & (Barye$\cdot$s/mL)\\ 
18. & Unstressed left ventricular volume ($V_{l}^0$) & 15.00 & (mL) \\
19. & Unstressed arterial volume ($V_{a}^0$) & 715.00 & (mL) \\
20. & Unstressed venous volume ($V_{v}^0$) & 2500.00 & (mL) \\
21. & Unstressed right ventricular volume ($V_{r}^0$) & 15.00 & (mL) \\
22. & Unstressed pulmonary arterial volume ($V_{pa}^0$) & 90.00 & (mL) \\
23. & Unstressed pulmonary venous volume ($V_{pv}^0$) & 490.00 & (mL) \\
\bottomrule
\end{tabular}
\end{center}
\caption{Input parameters ($\bv$) for the CVSim-6 system.}
\label{table:cvsim6 input paras}}
\end{table}

The CVSim-6 system is governed by the six ODEs in~\eqref{equ: cvsim6-ode}, one per compartment, with six pressures as state variables, i.e., $[P_l,P_a,P_v,P_r,P_{pa},P_{pv}]$ whose evolution is governed by a combination of Kirchhoff's law, two-element Windkessel models, and four ideal (no regurgitation) unidirectional heart valves, as detailed in~\cite{davis1991teaching}.
Generalized Ohm equations between pressures and flows are reported in~\eqref{equ: cvsim6-flow}, where $\mathbb{I}_{(\cdot)}$ denotes the indicator function, and the time-dependent ventricular compliances $C_l(t)$, $C_r(t)$ are defined in Appendix~\ref{sec: cvsim6-details}.

Solving the above ODE system also requires proper specification of initial conditions for the pressure in each compartment. This is accomplished through the solution of a linear system of equations that ensures mass conservation at $t=0$. 
For additional information, interested readers may refer to~\cite{davis1991teaching, heldt2010cvsim}, or equations~\eqref{equ: cvsim6-ic} in Appendix~\ref{sec: cvsim6-details}.

\begin{minipage}{0.55\textwidth}
\begin{equation}
    \hspace{-1.0cm}
    \left\{
    \begin{aligned}
        \dot{P}_l(t) & = \displaystyle  \frac{Q_{l,in}(t) - Q_{l,out}(t) - \Big(P_l(t) - P_{th}\Big)\dot{C_l}(t)}{C_l(t)} \\\vspace{-0.35cm}
        \dot{P}_a(t) & = \displaystyle \frac{ Q_{l,out}(t) - Q_a(t) }{C_a} \\\vspace{-0.35cm}
        \dot{P}_v(t) & = \displaystyle \frac{Q_a(t) - Q_{r,in}(t)}{C_v} \\\vspace{-0.35cm}
        \dot{P}_r(t) & = \displaystyle \frac{ Q_{r,in}(t) - Q_{r,out}(t) - \Big(P_r(t) - P_{th}\Big)\dot{C}_r(t) }{ C_r(t)} \\\vspace{-0.35cm}
        \dot{P}_{pa}(t) & = \displaystyle \frac{Q_{r,out}(t) - Q_{pv}(t)}{C_{pa}} \\\vspace{-0.35cm}
        \dot{P}_{pv}(t) & = \displaystyle \frac{Q_{pv}(t) - Q_{l,in}(t)}{C_{pv}}
    \end{aligned}
    \right.\hspace{-0.5cm}
    \label{equ: cvsim6-ode}
\end{equation}
\end{minipage}%
\begin{minipage}{0.44\textwidth}
\begin{equation}
    \left\{
    \begin{aligned}
        Q_{l,in}(t) & = \displaystyle \frac{P_{pv}(t) - P_l(t)}{R_{l,in}}\mathbb{I}_{P_{pv}(t) > P_l(t)} \\\vspace{-0.35cm}
        Q_{l,out}(t) & = \displaystyle \frac{P_{l}(t) - P_a(t)}{R_{l,out}}\mathbb{I}_{P_{l}(t) > P_a(t)} \\\vspace{-0.35cm}
        Q_a(t) &= \displaystyle \frac{P_a(t) - P_v(t)}{R_a} \\\vspace{-0.35cm}
        Q_{r,in}(t) & = \displaystyle \frac{P_{v}(t) - P_r(t)}{R_{r,in}}\mathbb{I}_{P_{v}(t) > P_r(t)} \\\vspace{-0.35cm}
        Q_{r,out}(t) & = \displaystyle \frac{P_{r}(t) - P_{pa}(t)}{R_{r,out}}\mathbb{I}_{P_{r}(t) > P_{pa}(t)}\\\vspace{-0.35cm}
    Q_{pv}(t) &= \displaystyle \frac{P_{pa}(t) - P_{pv}(t)}{R_{pv}}    \end{aligned}
    \right.\hspace{-0.3cm}
    \label{equ: cvsim6-flow}
\end{equation}
\end{minipage}%

\subsection{InVAErt networks}\label{sec:inVAErt}

\begin{figure}[!ht]
    \centering
    \includegraphics[scale=0.6]{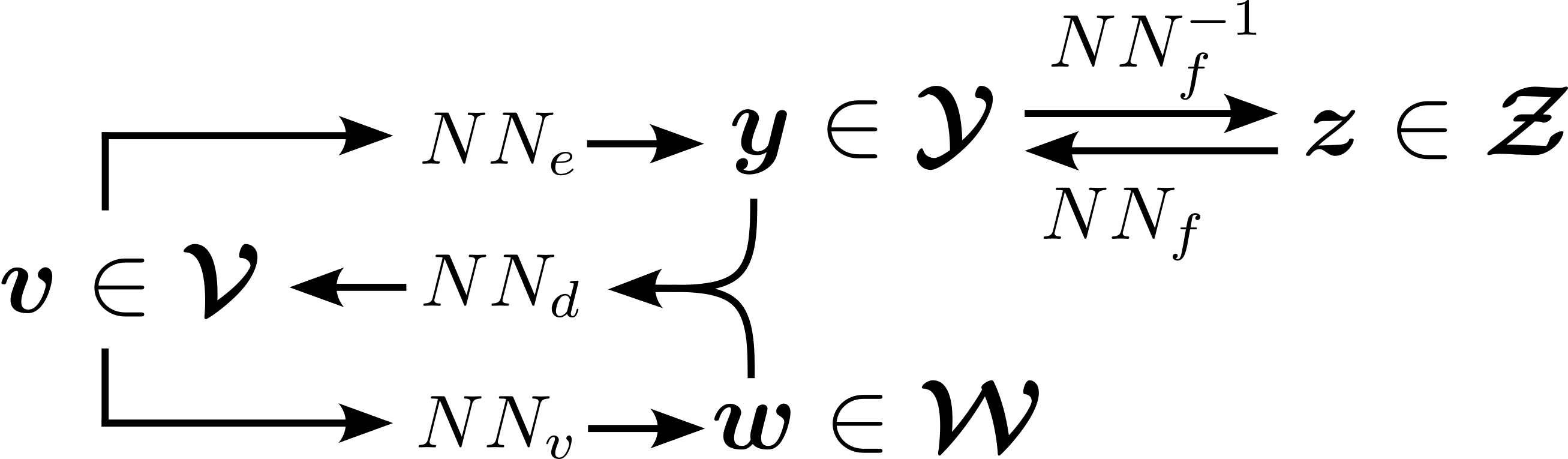}
    \caption{Schematic of all components of an inVAErt network and their interactions.}
    \label{fig:invert-maps}
\end{figure}

InVAErt networks~\cite{tong2024invaert} provide a comprehensive, data-driven methodology for physics-based digital twins (see diagram in Figure~\ref{fig:invert-maps}).
They consist of (1) an emulator $NN_e: \boldsymbol{\mathcal{V}}\mapsto \boldsymbol{\mathcal{Y}}$ as a surrogate for the input-to-output map $f$, (2) a $L$-layer flow-based density estimator $NN_f: \boldsymbol{\mathcal{Z}}\mapsto \boldsymbol{\mathcal{Y}}$ for the output distribution $p(\by)$, (3) a variational encoder $NN_v: \boldsymbol{\mathcal{V}}\mapsto \boldsymbol{\mathcal{W}}$ generating an input-dependent latent space responsible for the lack of bijectivity between inputs and outputs. Finally, (4) a decoder $NN_d: \boldsymbol{\mathcal{Y}} \times \boldsymbol{\mathcal{W}}\mapsto \boldsymbol{\mathcal{V}}$ approximates the \emph{inverse} map from combined outputs and latent variables to model inputs.
Given a dataset $\mathcal{D} = \big\{(\bv_j, \by_j)\big\}_{j=1}^{N}$, optimal weights and biases for each inVAErt network component are obtained by minimizing a loss function of the form
\begin{equation}
\begin{cases}
     \bphi_e^{\text{opt}} &= \displaystyle \argmin\limits_{\bphi_e} \sum\limits_{j=1}^N \Big\|\by_j - NN_e(\bv_j; \bphi_e) \Big\|_2^2 \ ,\\
    \vspace{-0.2cm}\\
    \bphi_f^{\text{opt}} &= \displaystyle \argmax\limits_{\bphi_f} \sum\limits_{j=1}^N \left( \log \pi_0(\bz^{(0)}_j) - \sum\limits_{l=1}^L \log\left|\det{\left(\frac{d\bz_j^{(l)}}{d\bz_j^{(l-1)}}\right)}\right| \right) \ , \\
    \vspace{-0.2cm}\\
    \bphi_v^{\text{opt}}, \bphi_d^{\text{opt}} &= \displaystyle  \argmin\limits_{\bphi_v, \bphi_d} \Bigg[ \beta_d \sum\limits_{j=1}^N \Big\|\bv_j - NN_d(\by_j, \bw_j; \bphi_d)\Big\|_2^2 + \frac{\beta_v}{2} \sum\limits_{j=1}^N\sum\limits_{k=1}^{\dim (\bw)} \Big( \mu_{jk}^2 +\sigma_{jk}^2 - \log(\sigma_{jk}^2) - 1 \Big) \\
    & \displaystyle  + \beta_r \sum\limits_{j=1}^N \Big\|\by_j - NN_e\Big( NN_d(\by_j, \bw_j; \bphi_d); \bphi_e^{\text{opt}} \Big) \Big\|_2^2 \Bigg],
\end{cases}
    \label{equ: nn parameters opt}
\end{equation}
which combines a MSE loss for the emulator $NN_e$, a MLE loss for the normalizing flow density estimator $NN_f$, and joint MSE and KL-divergence losses for the variational encoder $NN_v$ and the decoder $NN_d$.

We utilize a Real-NVP (real-valued non-volume preserving transformation) discrete normalizing flows architecture~\cite{dinh2016density}, which consists of $L$ bijections defined as
$$NN_f(\bz^{(0)}; \bphi_f) = (g_L \circ \cdots \circ g_2 \circ g_1)( \bz^{(0)} ), \quad \bz^{(0)} \sim \pi_0(\bz) \ ,$$
where $\bz^{(l)} = g_{l} (\bz^{(l-1)}), \bz^{(l-1)} = g_{l}^{-1} (\bz^{(l)}), l = 1,\dots,L$, and $\bz^{(L)} \approx \by \sim p(\by)$.
Additional details on approaches for discrete and continuous normalizing flows can be found in~\cite{papamakarios2021normalizing,kobyzev2020normalizing}.

The losses associated with the variational encoder and dense decoder are penalized with constants $\beta_v$, $\beta_d$, and $\beta_r$, respectively. Generation of input-dependent latent variables $\bw \sim p(\bw|\bv)$ follows the classical VAE framework~\cite{kingma2013auto}, i.e.,
\begin{equation}
    NN_v(\bv; \bphi_v) = \big[\bmu, \log \bsigma^2\big]^T, \quad \bw = \bmu + \bepsilon \odot \bsigma, \quad \bepsilon \sim \mathcal{N}(\boldsymbol{0}, \mathbf{I}) \ ,
\end{equation}
and the KL-divergence regularizes the latent space $\boldsymbol{\mathcal{W}}$ by reducing the statistical distance between the posterior distribution $p(\bw|\bv)$ and a standard normal prior $p(\bw)$.
Finally, the output reconstruction loss (sometimes referred to as the \emph{re-evaluation or re-constraining} loss~\cite{tong2024invaert}) uses the previously trained emulator, i.e., $NN_e(\cdot, \bphi_e^{\text{opt}})$, to make sure decoded parameters lead to model outputs that are close to the observations used as inputs to the decoder.

\subsubsection{Noise injection as training data augmentation}\label{sec:inVAErt-training noise}

To tackle challenging inference tasks involving practical non-identifiability, we add artificial Gaussian noise to the labels $\by$ during the training for $NN_f$, $NN_v$ and $NN_d$, as a regularization strategy.
Incorporating training noise as a form of data augmentation has been shown to improve model generalization and reduce overfitting~\cite{wang1999training,bishop1995training}. 
This is helpful in inference tasks when dealing with out-of-distribution data, measurement noise, and model misspecification (see Section~\ref{sec: cvsim6-ehr}). 

Note that the forward map $f$ is assumed to be deterministic in this study, hence no noise injection is performed when training the neural emulator $NN_e$.
However, when training all other inVAErt network components, we consider the output labels corrupted by zero-mean, uncorrelated and heteroskedastic Gaussian noise, i.e.,
\begin{equation}
    \by^{\prime} = \by + \boldsymbol{\eta}, \quad \boldsymbol{\eta} \sim \mathcal{N}\big(\boldsymbol{0}, \text{diag}(\boldsymbol{s}^2)\big), \quad \boldsymbol{s}^2 = [s_1^2, \cdots, s^2_{\dim(\by)}]^T \ ,
\end{equation}
with known standard deviations for each output component.
The random noise is independent to both the input and output data, with a new realization generated at each epoch.
In summary, $\by^{\prime}$ replaces $\by$ in all terms in~\eqref{equ: nn parameters opt}, except for the MSE emulator loss.

\section{Numerical experiments}\label{sec:numerical ex}

In this section, we provide a number of new results for inVAErt networks in the context of computational physiology. 
We first characterize the stiffness of the differential system used throughout this section, showing that lumped parameters models give rise to stiff system of ODEs. We also highlight the most important mechanisms that are responsible for this stiffness. 
We then use synthetically generated data from a baseline set of parameters (the so-called \emph{default} parameter set) to show the ability of inVAErt networks to identify entire manifolds of non-identifiable parameters and provide ways to graphically characterize such manifolds. 
Finally, we show how to select the parameter prior, the amount of noise, and how to deal with missing data, demonstrating the use of inVAErt networks for parametric inversion in the presence of practical non-identifiability.

\subsection{Stiffness analysis}\label{sec: cvsim6-stiffness}

Like other pulsatile cardiovascular models driven by periodic forcing~\cite{marquis2018practical}, the CVSim-6 system exhibits stiffness and may suffer from numerical instability without a careful selection of the integration time step.
To formally study this phenomenon, we first write equations~\eqref{equ: cvsim6-ode} in matrix form
\begin{equation}
    \dot{\mathbf{P}}(t) = \mathbf{A}(t)\,\mathbf{P}(t) + \boldsymbol{b}(t) \ ,
    \label{equ: Ap+b}
\end{equation}
where $\mathbf{P}(t) = [P_l(t), P_a(t), P_v(t), P_r(t), P_{pa}(t), P_{pv}(t)]^T$ represents the state vector, $\mathbf{A}(t)$ denotes the time-dependent coefficient matrix and $\boldsymbol{b}(t)$ is a time-varying forcing vector which quantifies the contribution of the intrathoracic pressure~\cite{davis1991teaching}.

Rather than following a rigorous definition, stiffness in ODE systems is often identified by its effects on numerical solutions~\cite{shi2022polynomial, ashino2000behind,lambert1991numerical}. 
One way to analyze stiffness is through the eigen-decomposition of the RHS Jacobian of~\eqref{equ: Ap+b}, i.e. $\mathbf{A}(t) = \mathbf{Q}_t \boldsymbol{\Lambda}_t \mathbf{Q}_t^{-1}$.
Here, $\boldsymbol{\Lambda}_t \in \mathbb{C}^{6\times 6}$ is a diagonal matrix with entries equal to the eigenvalues $\lambda_1(t), \cdots, \lambda_6(t)$, with $|\text{Re}\big(\lambda_1(t)\big)|\geq \cdots \geq |\text{Re}\big(\lambda_6(t)\big)|$, and $\mathbf{Q}_t \in \mathbb{C}^{6\times 6}$ is the matrix of eigenvectors, i.e., $\mathbf{Q}_t = [\boldsymbol{Q}_1(t),\cdots,\boldsymbol{Q}_6(t)]$.
As $\mathbf{A}(t)$ depends on the state vector $\mathbf{P}(t)$ via indicator functions (see~\eqref{equ: cvsim6-flow}), $\mathbf{A}(t)$ is approximated using a 5-th order implicit Runge-Kutta method (Radau) with adaptive step~\cite{wanner1996solving} as discussed in~\cite{davis1991teaching}.
The simulation runs for 12 cardiac cycles, i.e., $12 \cdot T_{tot} = 12\cdot 60/Hr$, utilizing the \emph{default} set of input parameters listed in Table~\ref{table:cvsim6 input paras}.

\begin{figure}[!ht]
\centering
\includegraphics[scale= 0.5]{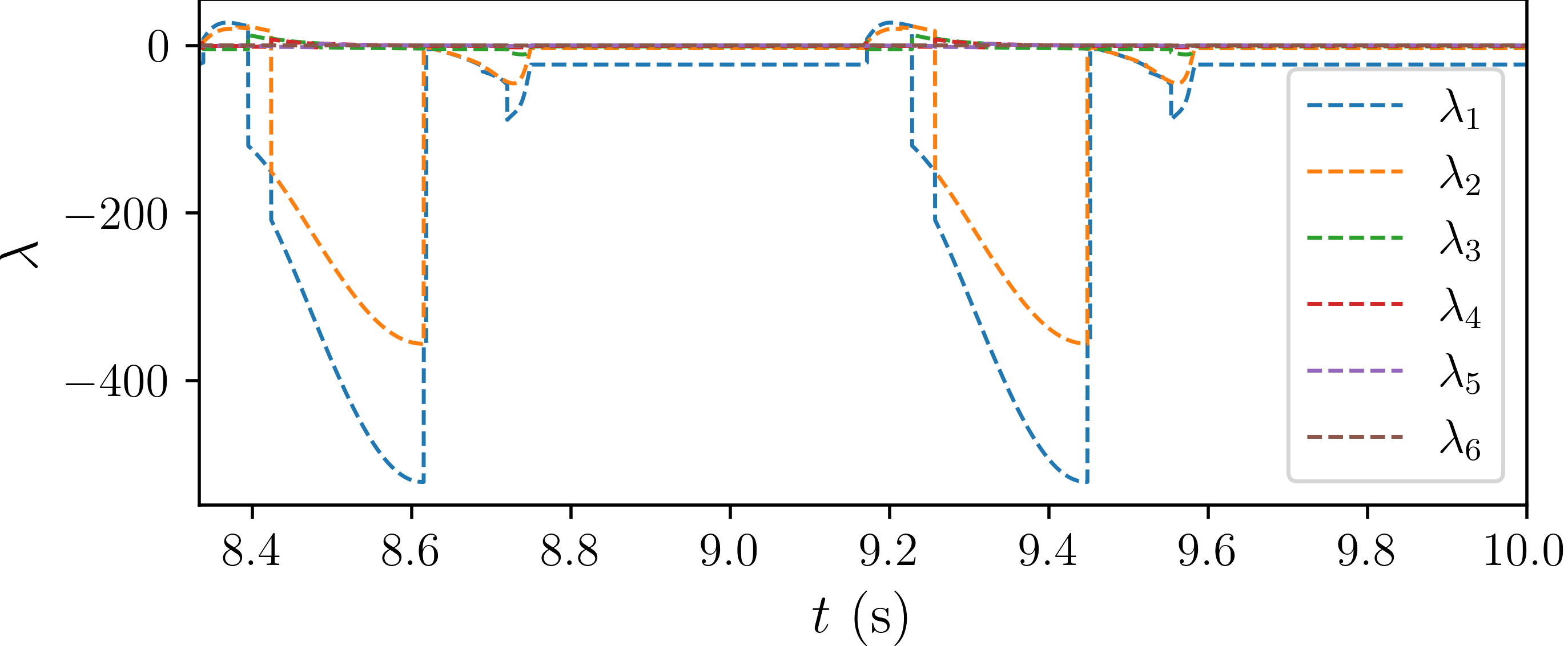}
    \caption{Time history of the eigenvalues of $\mathbf{A}(t)$ plotted over two heart cycles.}
    \label{fig: eig-decay}
\end{figure}

All six eigenvalues of $\mathbf{A}(t)$ are real, with time histories over the last two heart cycles shown in Figure~\ref{fig: eig-decay}. 
The first two modes are associated with the largest magnitudes occurring at the same time during each heart cycle.

At any time $t$, the \emph{stiffness ratio} SR$(t)$ is defined as 
\begin{equation}
    \text{SR}(t) = \frac{\max|\textrm{Re}(\boldsymbol{\Lambda}_t)|}{\min|\textrm{Re}(\boldsymbol{\Lambda}_t)|} \ , 
    \label{equ: sr}
\end{equation}
and can inform on the \emph{degree} of stiffness of a given system.
Specifically, a large SR implies significant difference in the rate of change between solution components, so that a smaller time step is required to capture both fast- and slow-varying dynamics~\cite{vijayarangan2024data, shi2022polynomial, ashino2000behind,lambert1991numerical}.
In addition, it is possible for $\min|\textrm{Re}(\boldsymbol{\Lambda}_t)|$ to be zero for some $t$ (see Figure~\ref{fig: eig-decay}), indicating a solution component remaining constant (e.g., due to a time-varying capacitance \emph{driving} the autonomous system response).
In such a case, the denominator in equation~\eqref{equ: sr} is replaced by the non-zero eigenvalue with the smallest magnitude, plus a tolerance $tol$ ($tol=1\cdot 10^{-14}$ is used here, see Remark~\ref{rmk: sr tol}).

\begin{figure}[!ht]
\begin{subfigure}[b]{0.163\textwidth}
        \centering
         \includegraphics[scale=0.18]{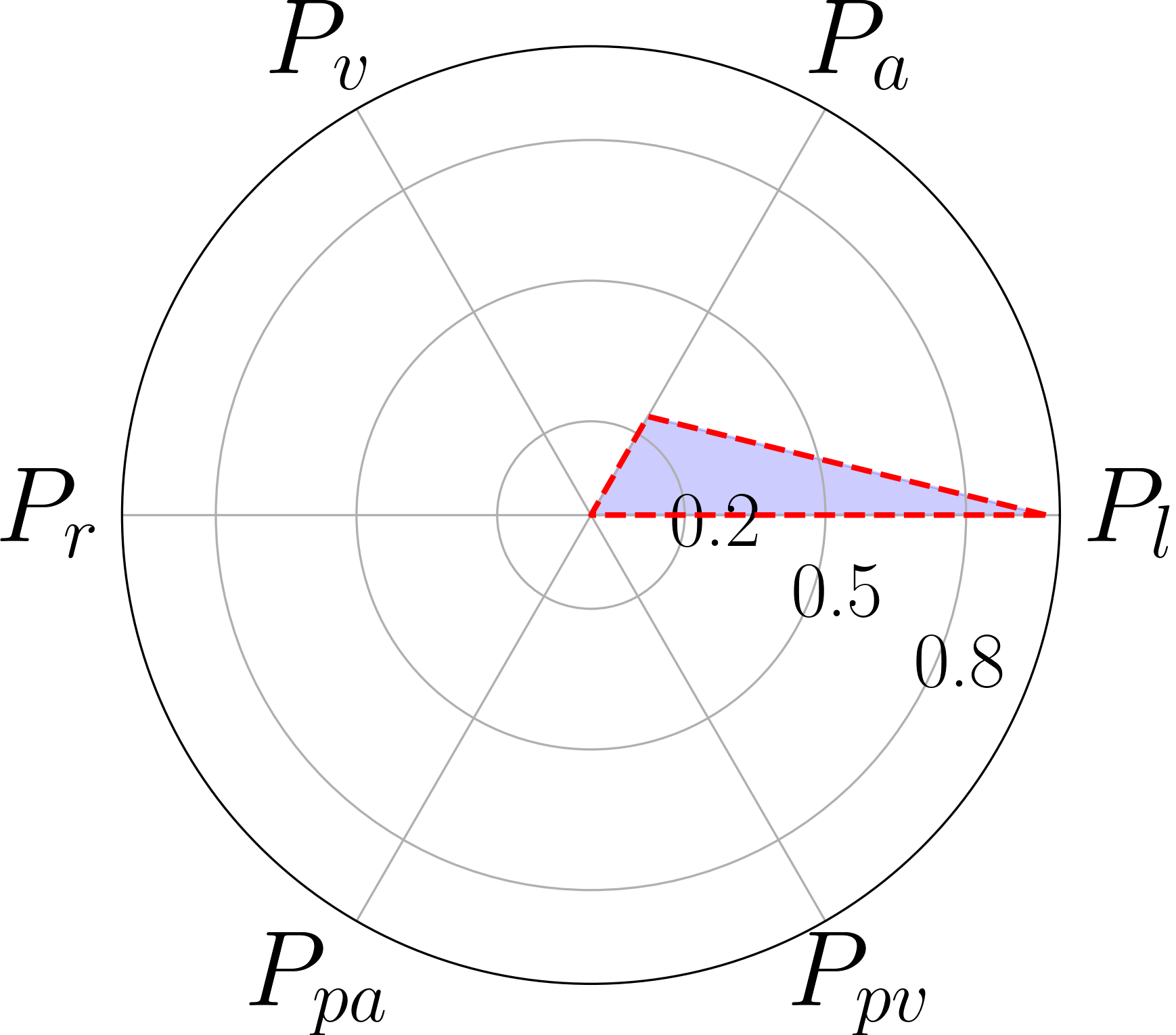}
         \caption{$|\boldsymbol{Q}_1|$.}
\end{subfigure}
\begin{subfigure}[b]{0.163\textwidth}
        \centering
         \includegraphics[scale=0.18]{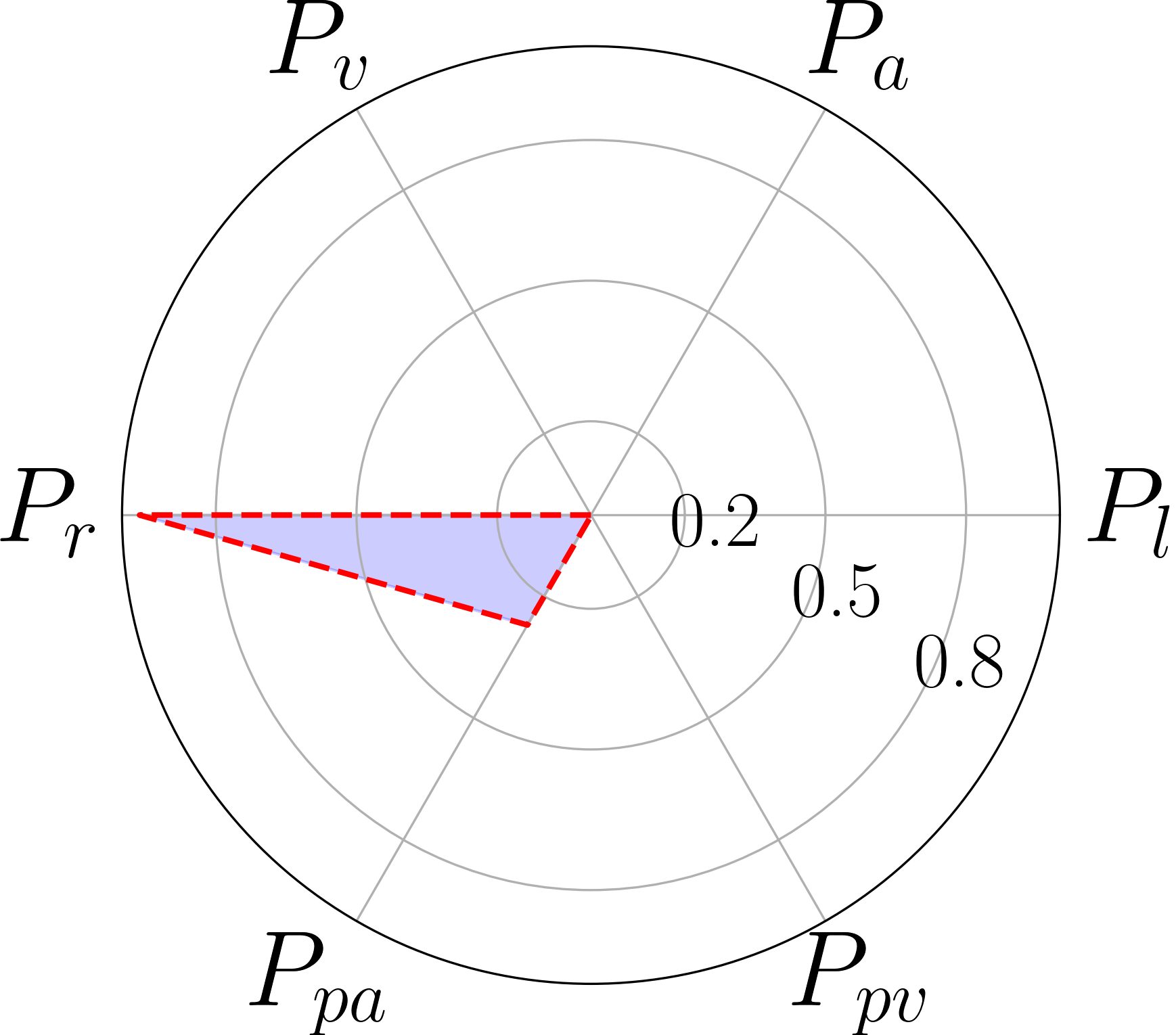}
         \caption{$|\boldsymbol{Q}_2|$.}
\end{subfigure}
\begin{subfigure}[b]{0.163\textwidth}
        \centering
         \includegraphics[scale=0.18]{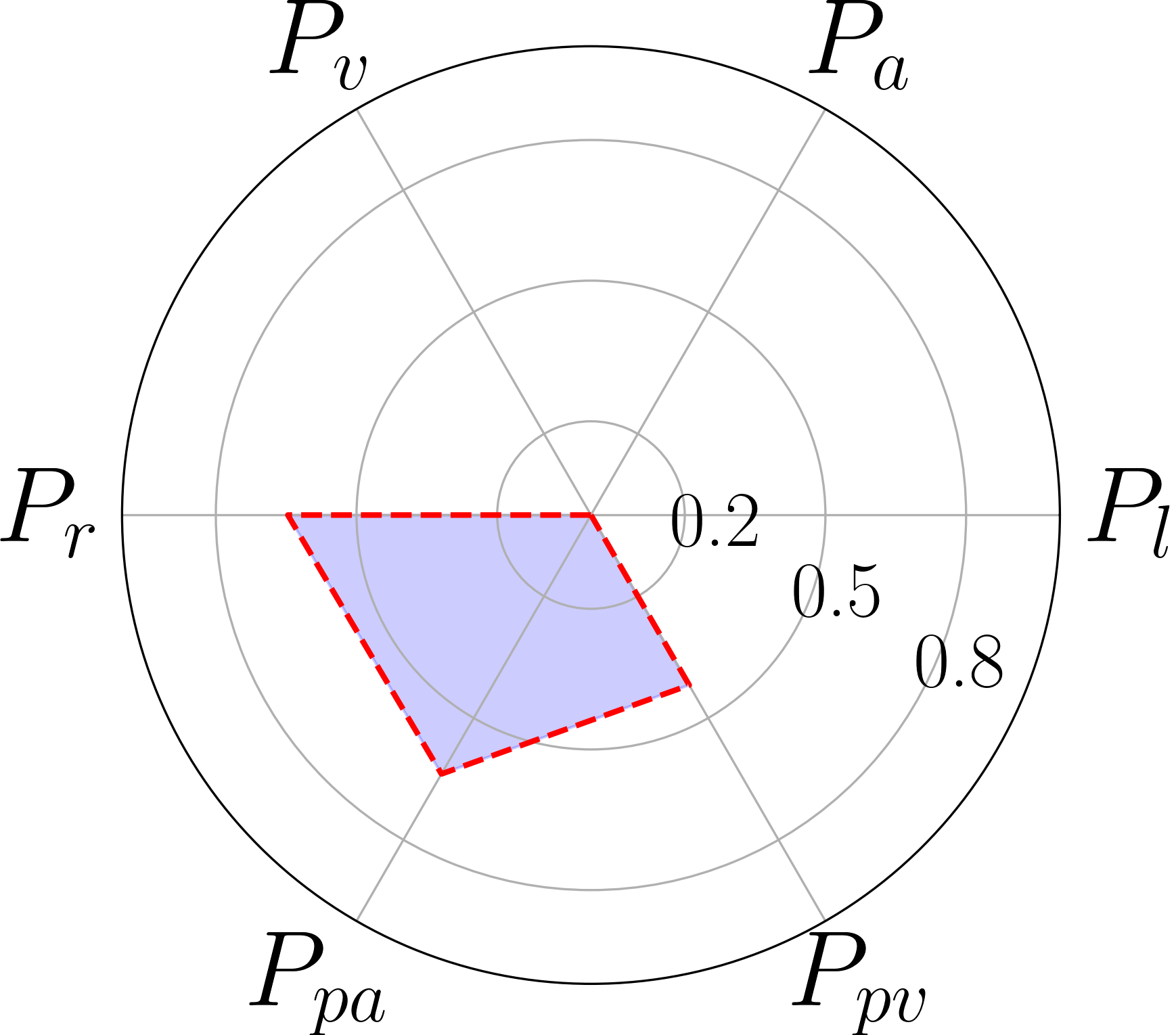}
         \caption{$|\boldsymbol{Q}_3|$.}
\end{subfigure}
\begin{subfigure}[b]{0.163\textwidth}
        \centering
         \includegraphics[scale=0.18]{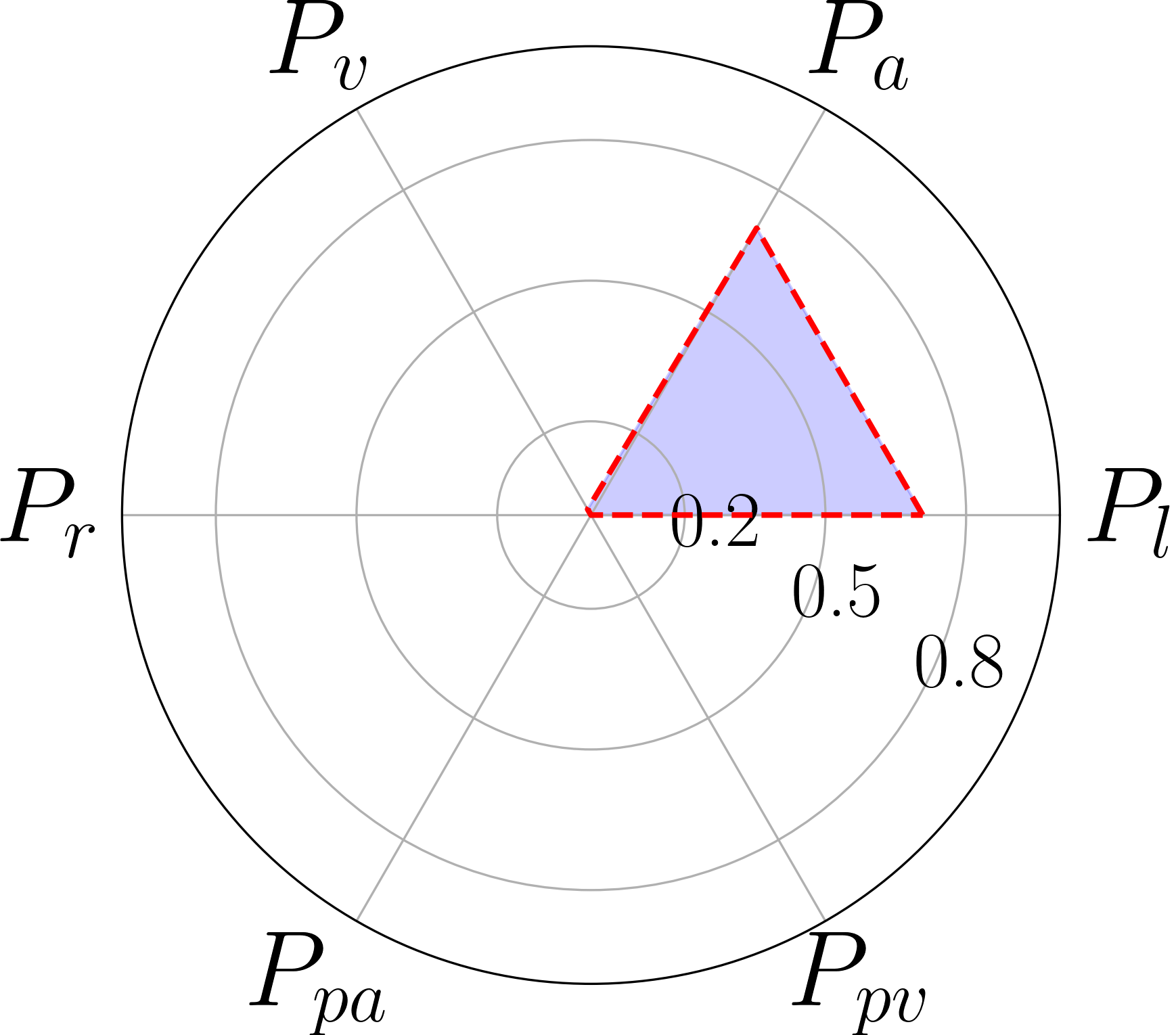}
         \caption{$|\boldsymbol{Q}_4|$.}
\end{subfigure}
\begin{subfigure}[b]{0.163\textwidth}
        \centering
         \includegraphics[scale=0.18]{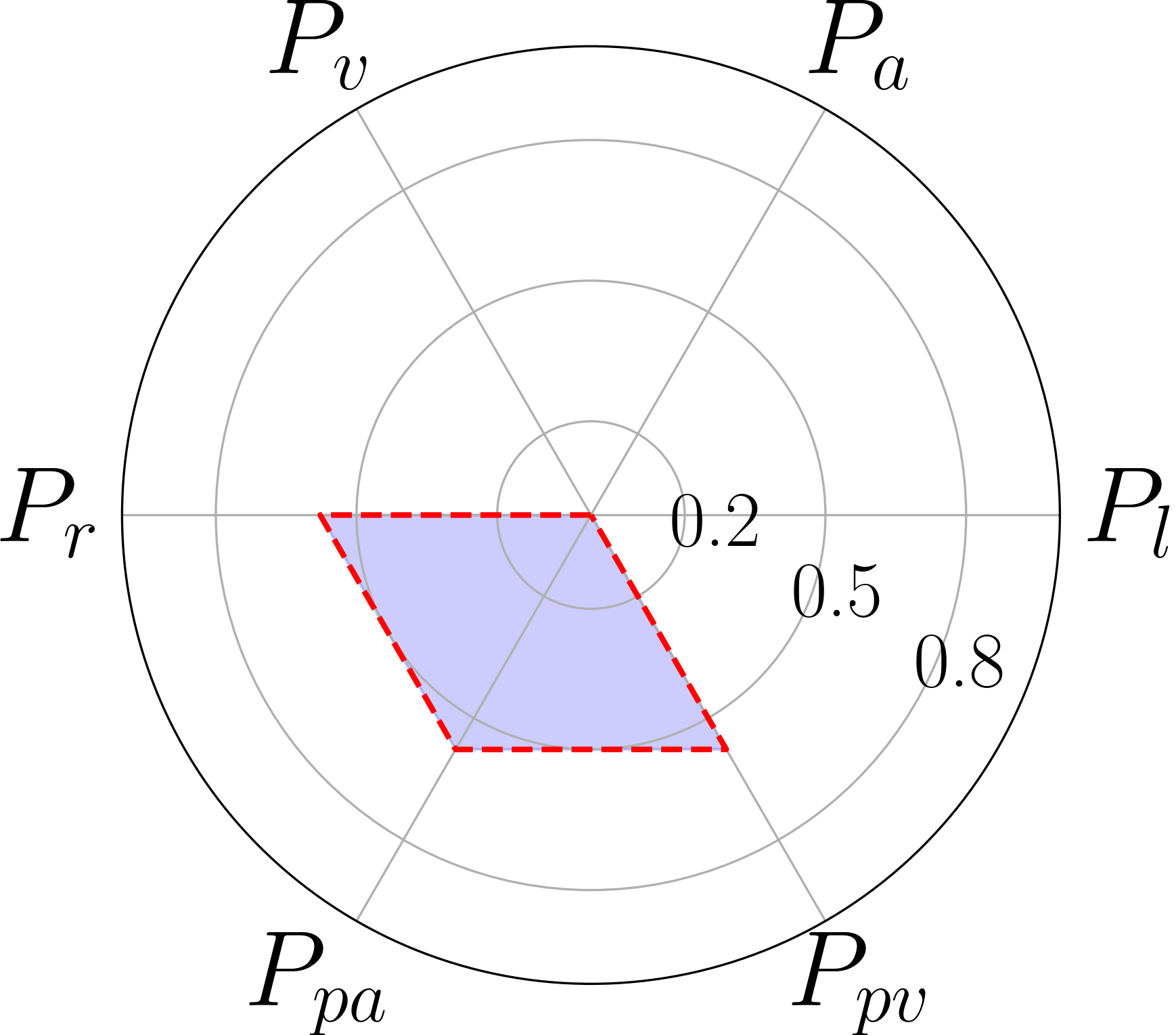}
         \caption{$|\boldsymbol{Q}_5|$.}
\end{subfigure}
\begin{subfigure}[b]{0.163\textwidth}
        \centering
         \includegraphics[scale=0.18]{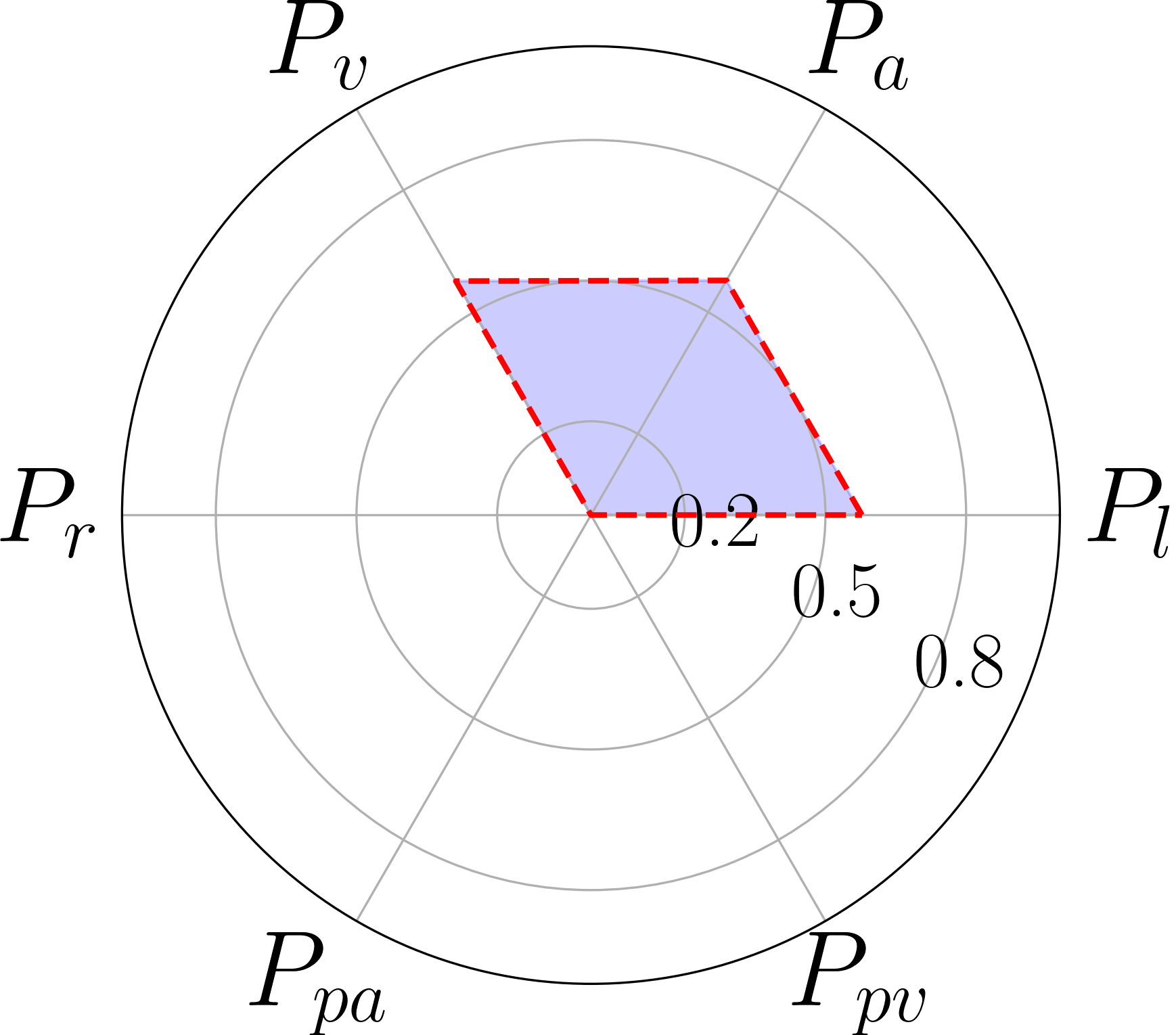}
         \caption{$|\boldsymbol{Q}_6|$.}
\end{subfigure}
\caption{Radar plots of absolute eigenvector components for maximum SR$(t)$ at $t\sim 9.444$ (s). Associated eigenvalues: $\lambda_1 = -520.95$, $\lambda_2 = -355.93$, $\lambda_3=-3.75$, $\lambda_4= -0.51$, $\lambda_5=6.25\cdot 10^{-4}$, $\lambda_6=2.91\cdot 10^{-5}$.}
\label{fig:eigvec-SRmax}
\end{figure}

\begin{figure}[!ht]
\begin{subfigure}[b]{0.163\textwidth}
        \centering
         \includegraphics[scale=0.18]{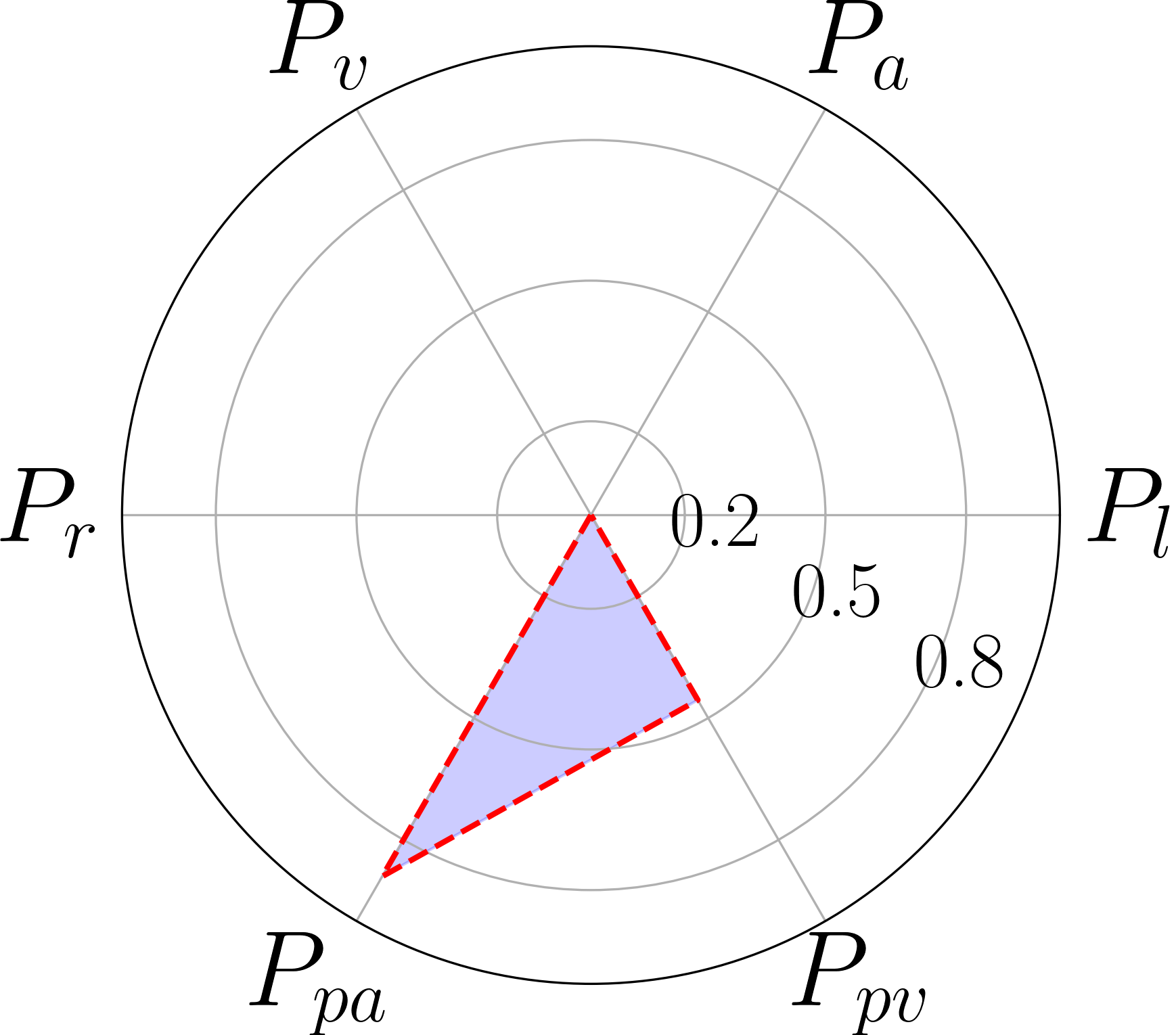}
         \caption{$|\boldsymbol{Q}_1|$.}
\end{subfigure}
\begin{subfigure}[b]{0.163\textwidth}
        \centering
         \includegraphics[scale=0.18]{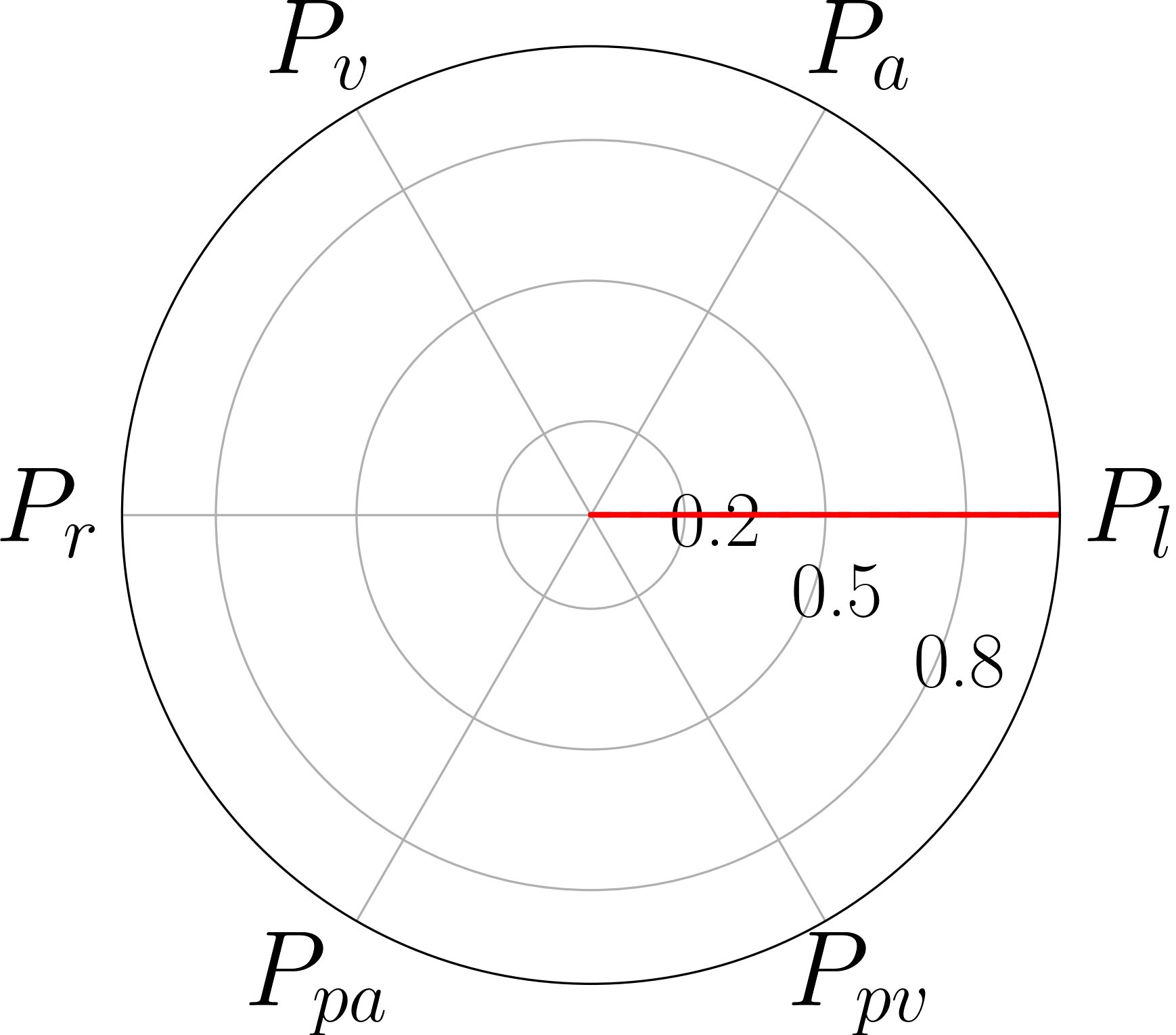}
         \caption{$|\boldsymbol{Q}_2|$.}
\end{subfigure}
\begin{subfigure}[b]{0.163\textwidth}
        \centering
         \includegraphics[scale=0.18]{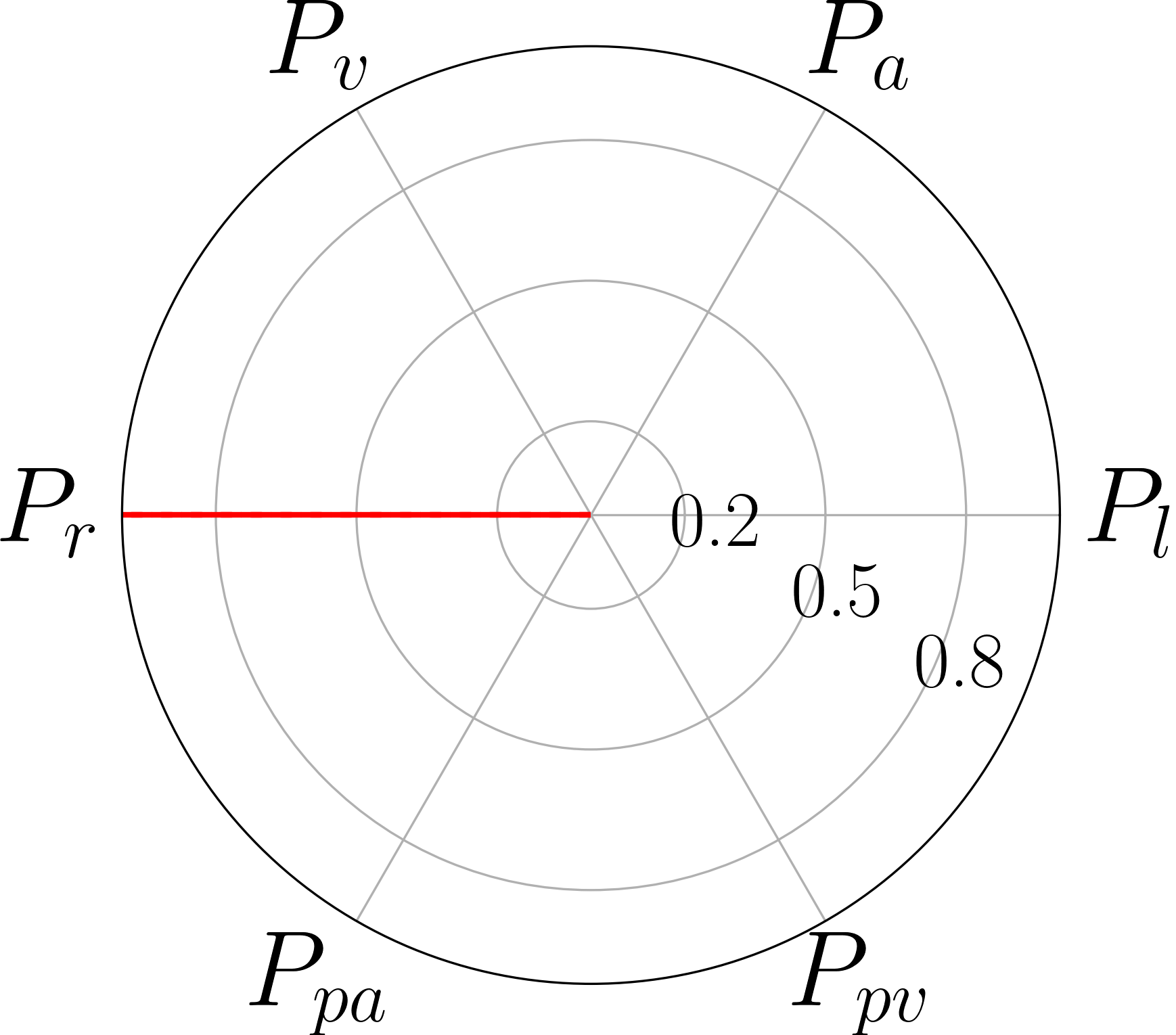}
         \caption{$|\boldsymbol{Q}_3|$.}
\end{subfigure}
\begin{subfigure}[b]{0.163\textwidth}
        \centering
         \includegraphics[scale=0.18]{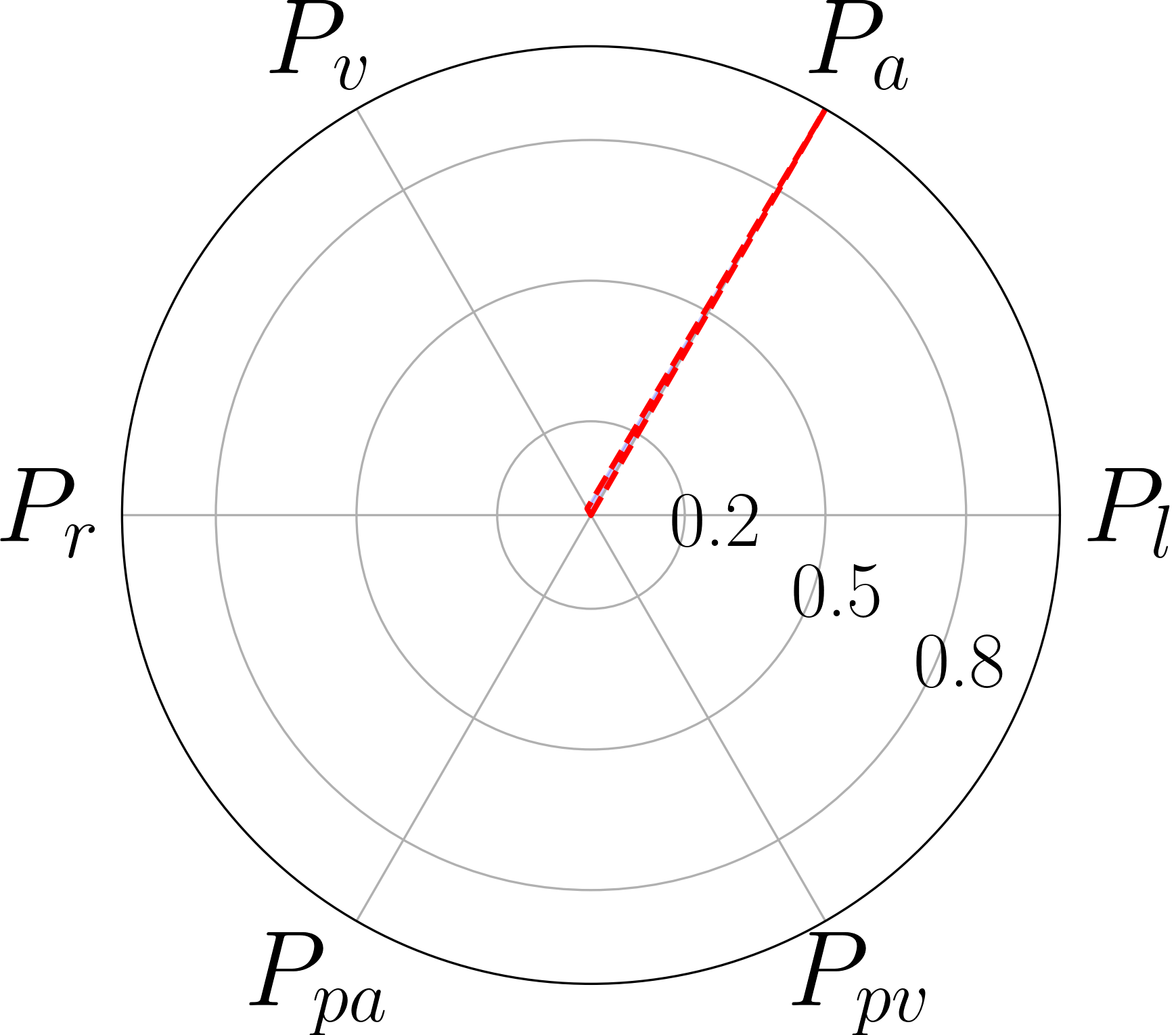}
         \caption{$|\boldsymbol{Q}_4|$.}
\end{subfigure}
\begin{subfigure}[b]{0.163\textwidth}
        \centering
         \includegraphics[scale=0.18]{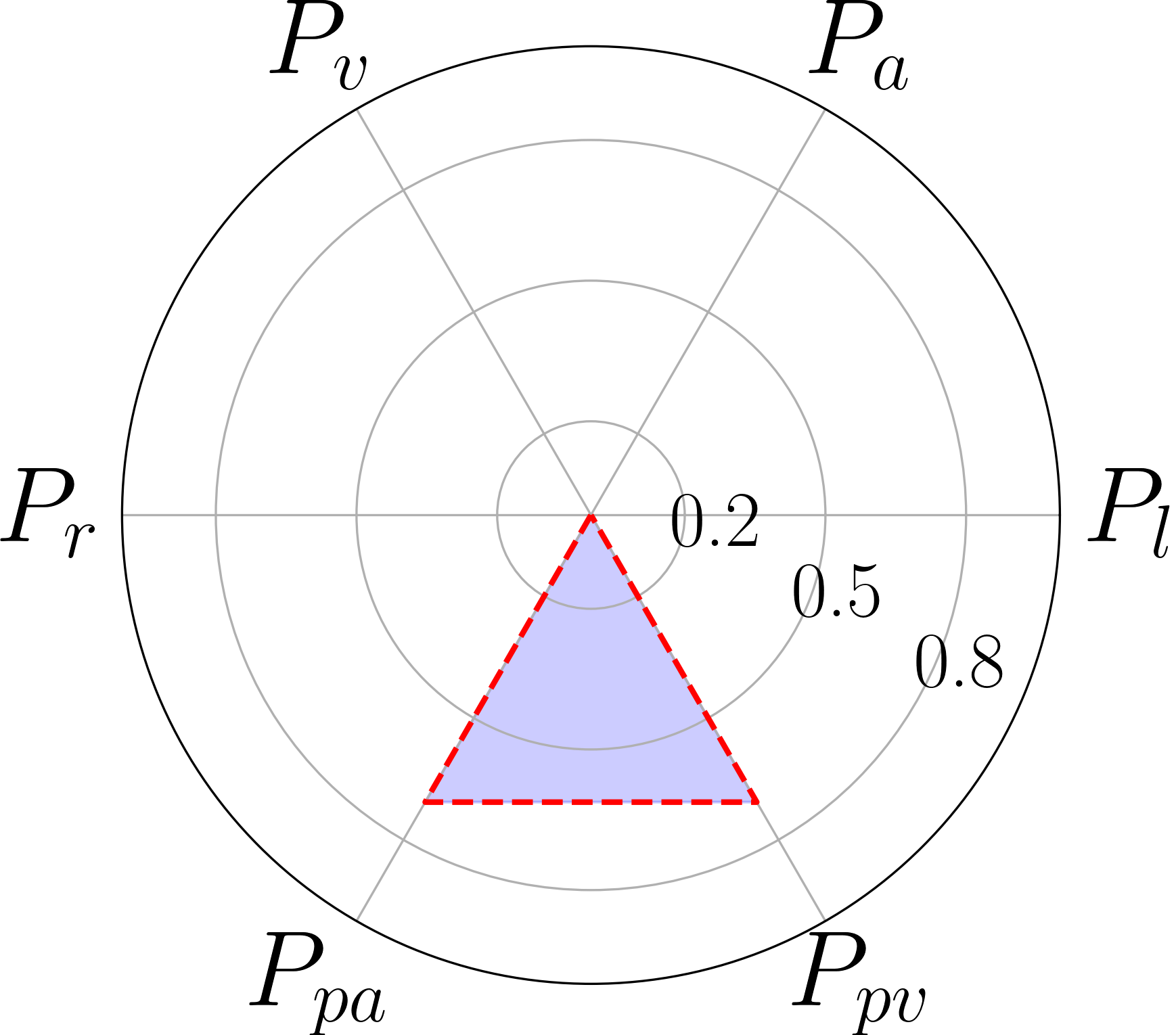}
         \caption{$|\boldsymbol{Q}_5|$.}
\end{subfigure}
\begin{subfigure}[b]{0.163\textwidth}
        \centering
         \includegraphics[scale=0.18]{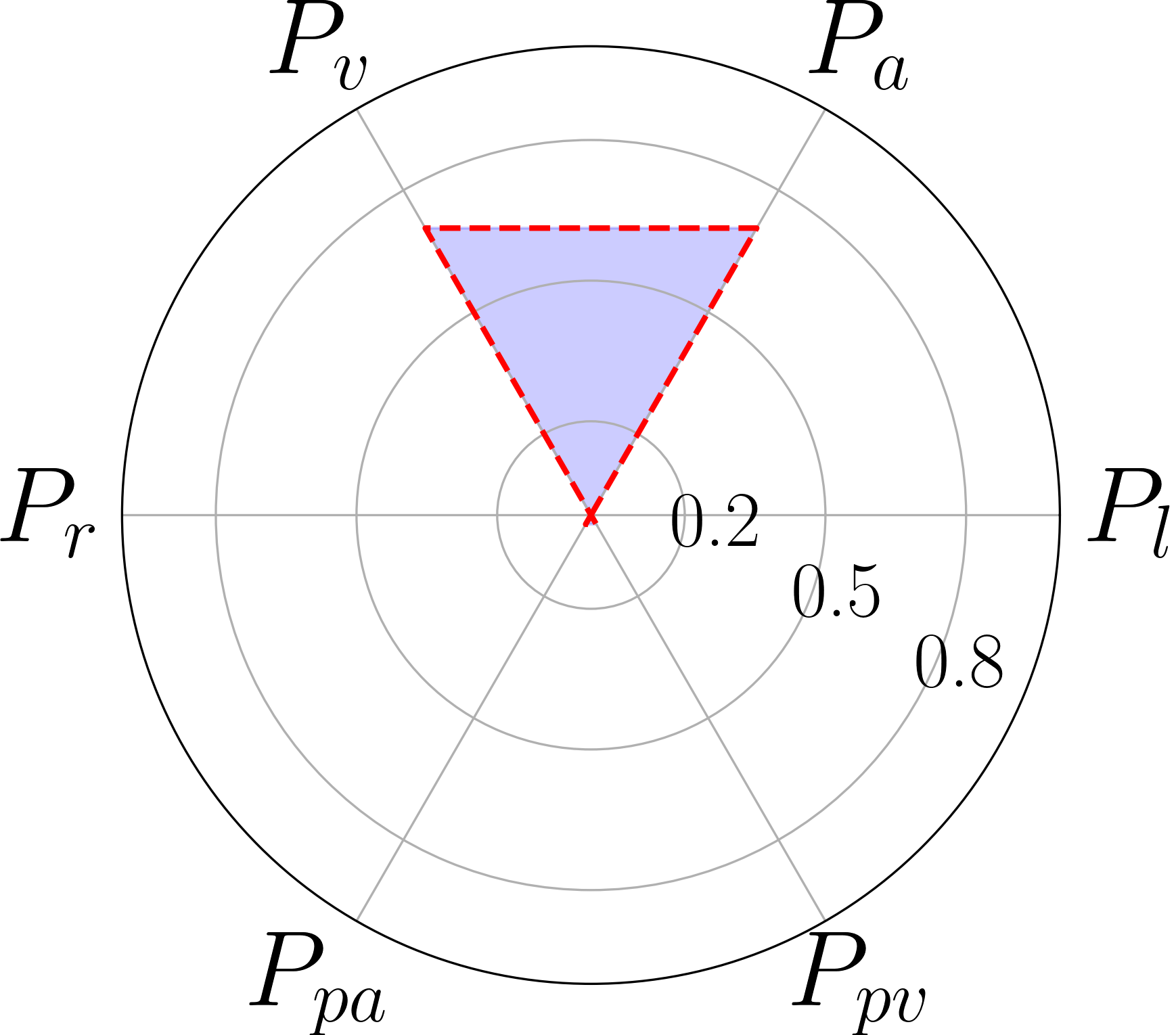}
         \caption{$|\boldsymbol{Q}_6|$.}
\end{subfigure}
\caption{Radar plots of absolute eigenvector components for minimum SR$(t)$ at $t\sim 8.618$ (s). Associated eigenvalues: $\lambda_1 = -4.40$, $\lambda_2 = -1.81$, $\lambda_3=-1.77$, $\lambda_4= -0.64$, $\lambda_5=2.22\cdot 10^{-16}$, $\lambda_6=3.47\cdot 10^{-18}$.}
\label{fig:eigvec-SRmin}
\end{figure}

In Figure~\ref{fig:eigvec-SRmax} and~\ref{fig:eigvec-SRmin}, we also show the absolute components of each eigenvector and report the corresponding eigenvalue for times at which the stiffness ratio SR$(t)$ reaches its minimum and maximum value, respectively, within the last two heart cycles.

\begin{figure}[!ht]
\begin{subfigure}[b]{0.498\textwidth}
        \centering
         \includegraphics[scale=0.48]{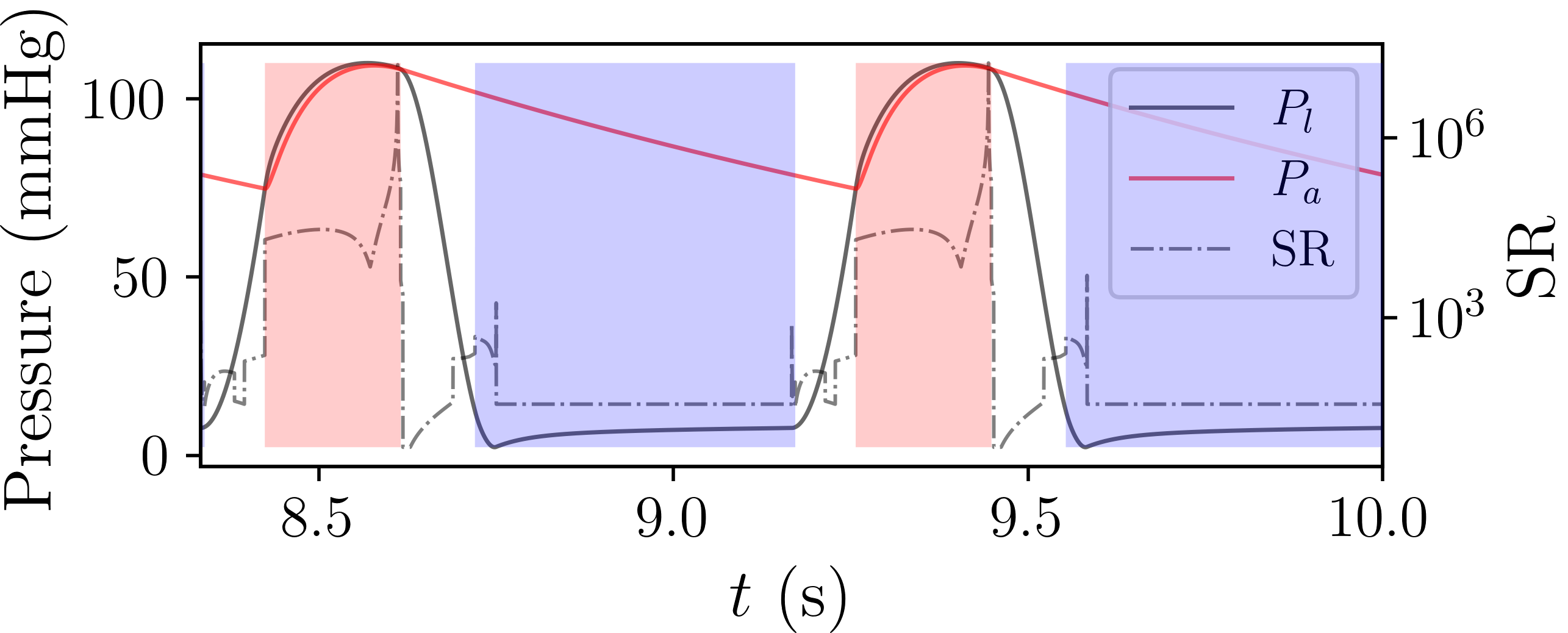}
\end{subfigure}
\begin{subfigure}[b]{0.498\textwidth}
        \centering
         \includegraphics[scale=0.48]{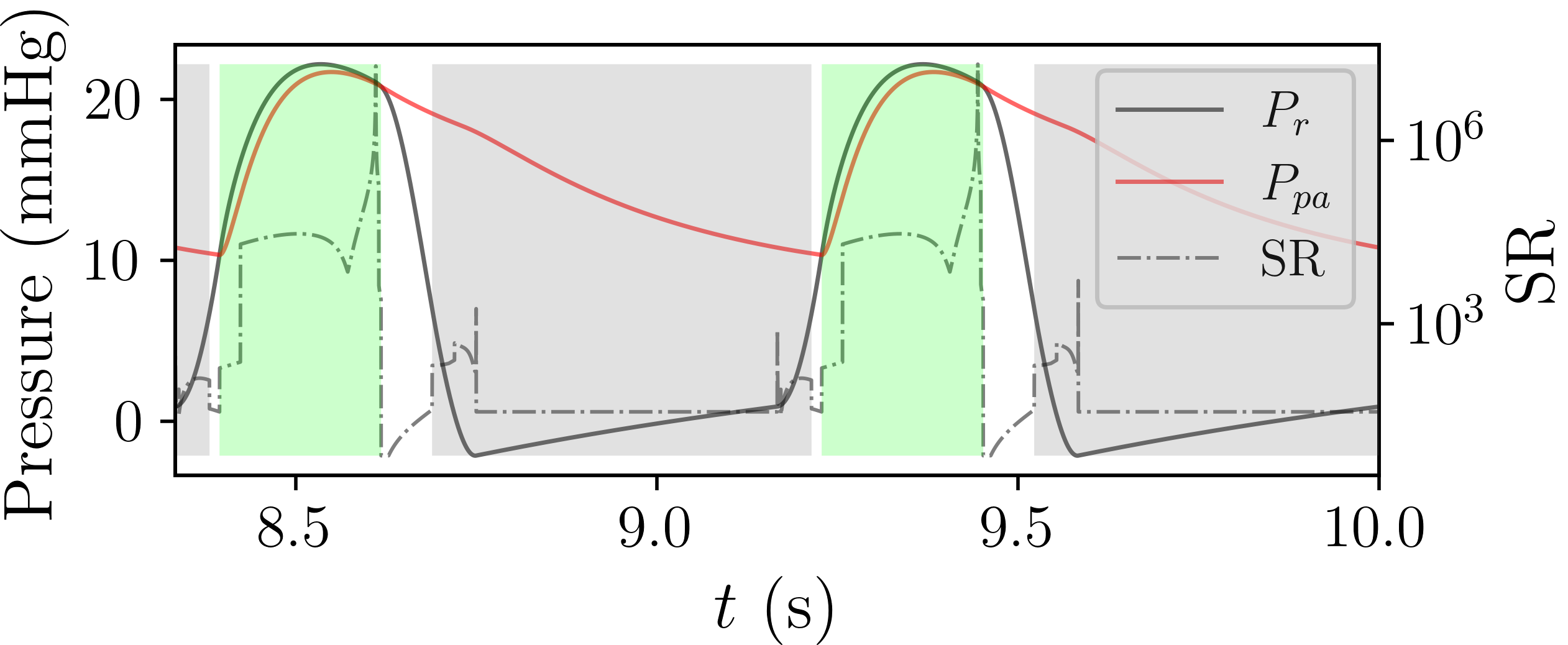}
\end{subfigure} \\
\vspace{-0.25cm} \\
\begin{subfigure}[b]{0.498\textwidth}
        \centering
         \includegraphics[scale=0.48]{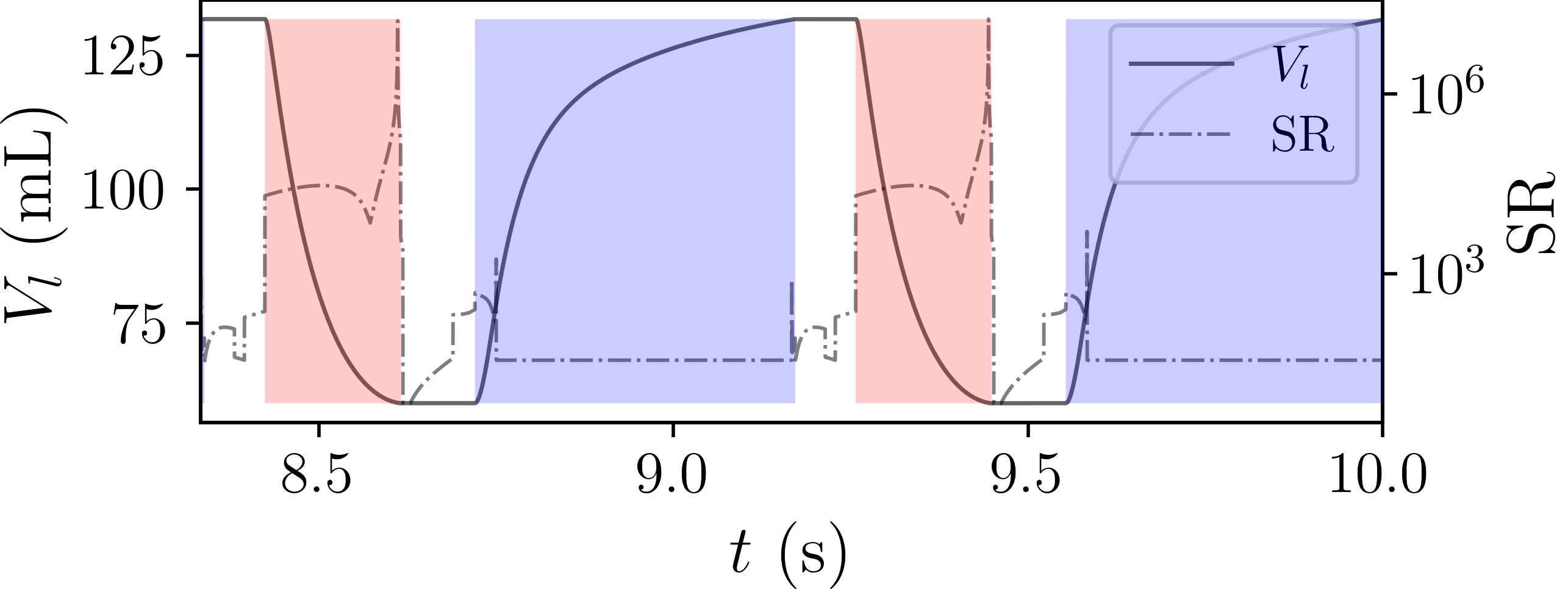}
\end{subfigure}
\begin{subfigure}[b]{0.498\textwidth}
        \centering
         \includegraphics[scale=0.48]{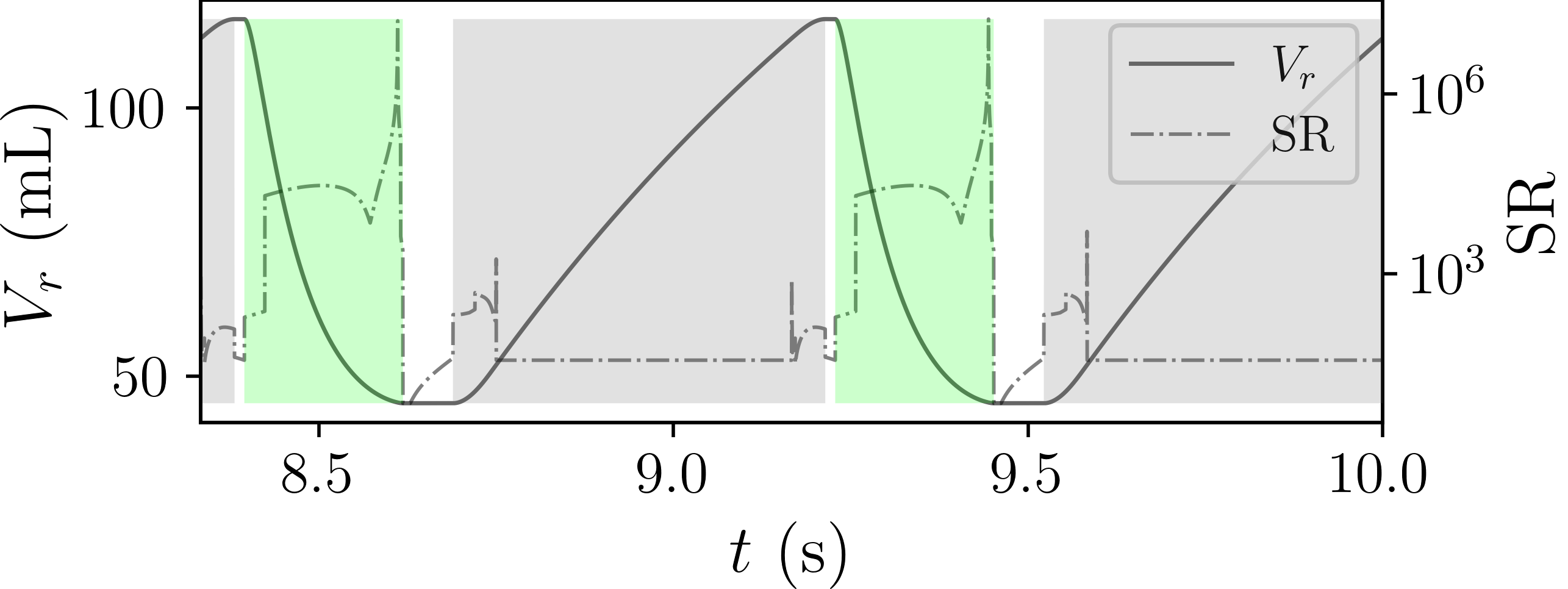}
\end{subfigure} \\
\vspace{-0.25cm} \\
\begin{subfigure}[b]{0.498\textwidth}
        \centering
         \includegraphics[scale=0.48]{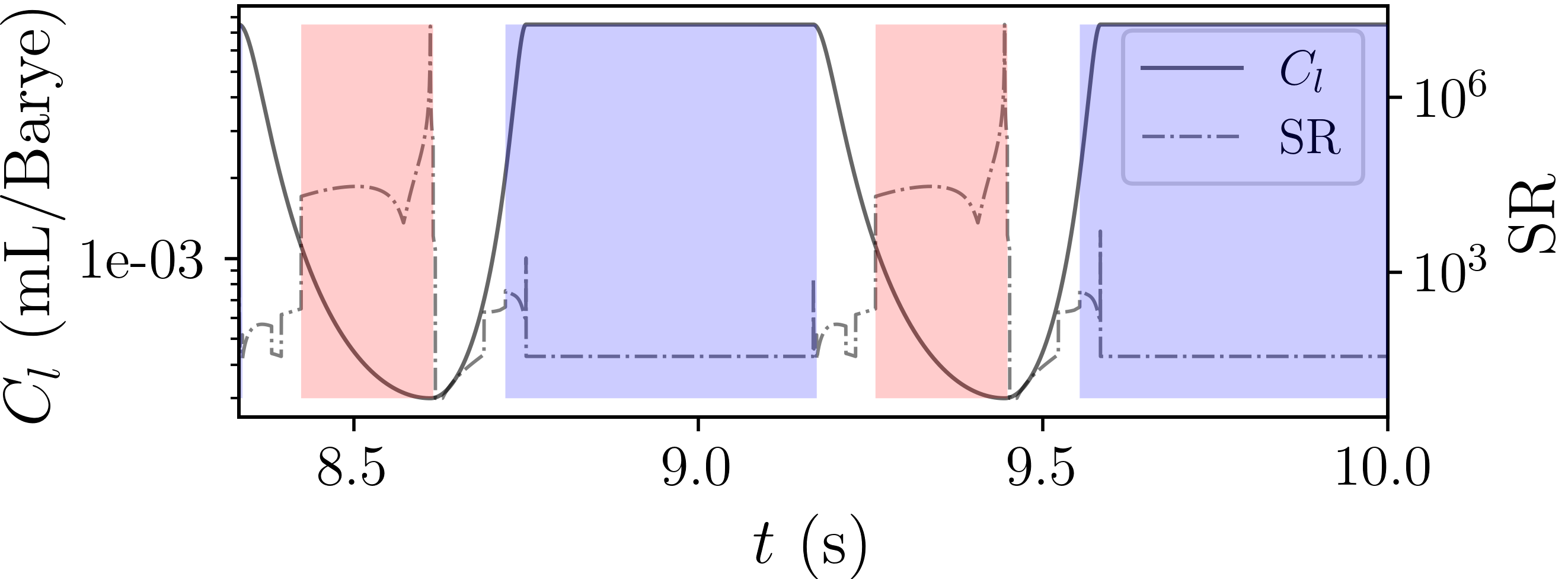}
\end{subfigure} 
\begin{subfigure}[b]{0.498\textwidth}
        \centering
         \includegraphics[scale=0.48]{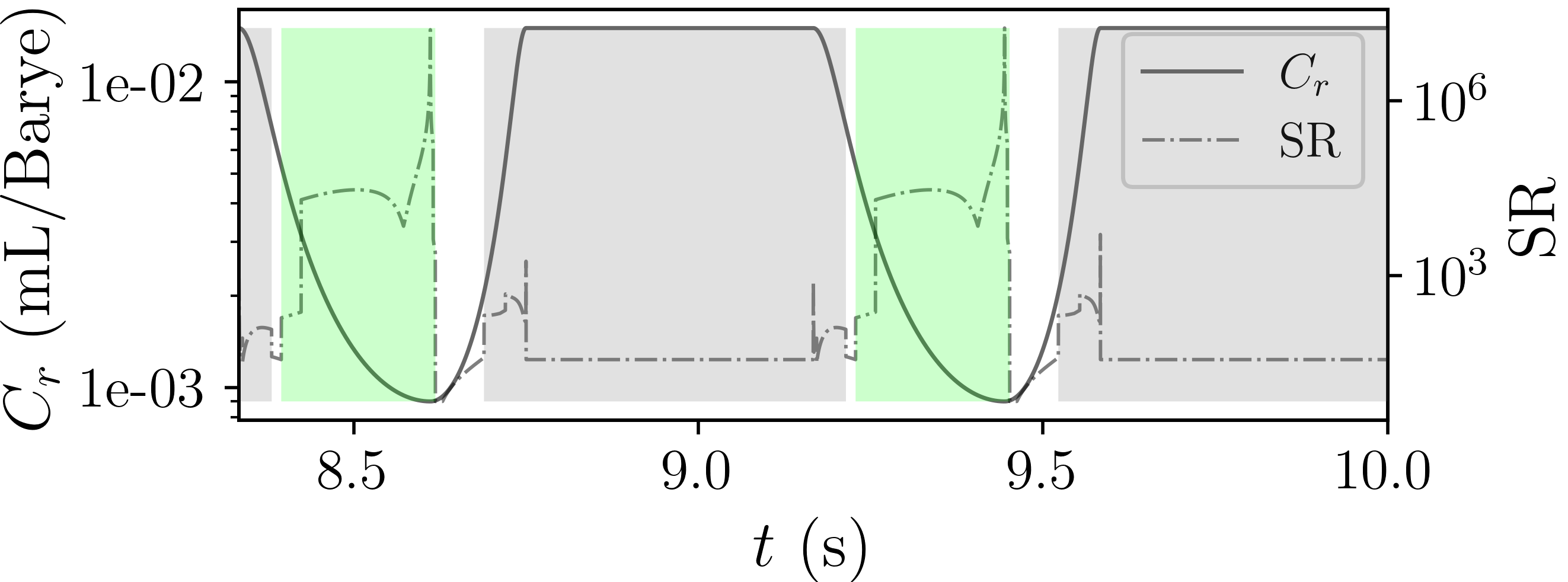}
\end{subfigure}\\
\vspace{-0.43cm} \\
\caption{Variation of stiffness ratio SR$(t)$ over two heart cycles. The CVSim-6  dynamics is also shown in terms of ventricular volumes ($V_l$, $V_{r}$), ventricular pressures ($P_l$, $P_r$), arterial pressure ($P_a$), pulmonary arterial pressure ($P_{pa}$), and time-varying ventricular capacitance ($C_l, C_r$).
The opening and closing intervals for the mitral (blue), aortic (red), tricuspid (gray), and pulmonary (green) valves are depicted using shaded colors.}
\label{fig: VP-eigSR}
\end{figure}

The maximum stiffness is observed when the magnitudes of $\lambda_1$ and $\lambda_2$ are maximized (see Figure~\ref{fig: eig-decay} and Figure~\ref{fig:eigvec-SRmax}), with intrinsic timescale constants equal to~\cite{vijayarangan2024data}
$$
\tau_1 = 1/|\lambda_1| = 0.0019 \ (\text{s}), \quad \tau_2 = 1/|\lambda_2| = 0.0028 \ (\text{s}).
$$
These values agree well with $RC$-constants computed at the outflow of the left and right ventricles, respectively (see Table 4.2 of~\cite{davis1991teaching}), defined as $RC_{j,k} = R_{j, k}C_{j}C_{k}/(C_{j} + C_{k})$ for neighboring compartments $j$ and $k$~\cite{davis1991teaching}.
The largest difference between solution components occurs due to ventricular-arterial coupling (see $|\boldsymbol{Q}_5|$ and $|\boldsymbol{Q}_6|$ in Figure~\ref{fig:eigvec-SRmax}), and is observed right before the closure of the aortic and pulmonary valve in late systole (see Figure~\ref{fig: VP-eigSR}).

Figure~\ref{fig: VP-eigSR} also shows how the stiffness ratio SR$(t)$ varies with respect to CVSim-6 system dynamics, specifically right and left ventricular pressure, volume, time-varying compliance, and valve opening times.
SR increases significantly during systole, peaking just before the aortic and pulmonary valves close. It then drops sharply to its minimum at the onset of isovolumetric relaxation, and remains low throughout the rest of diastole.

\begin{figure}[!ht]
    \centering
    \begin{subfigure}[b]{0.49\textwidth}
        \centering
         \includegraphics[scale=0.56]{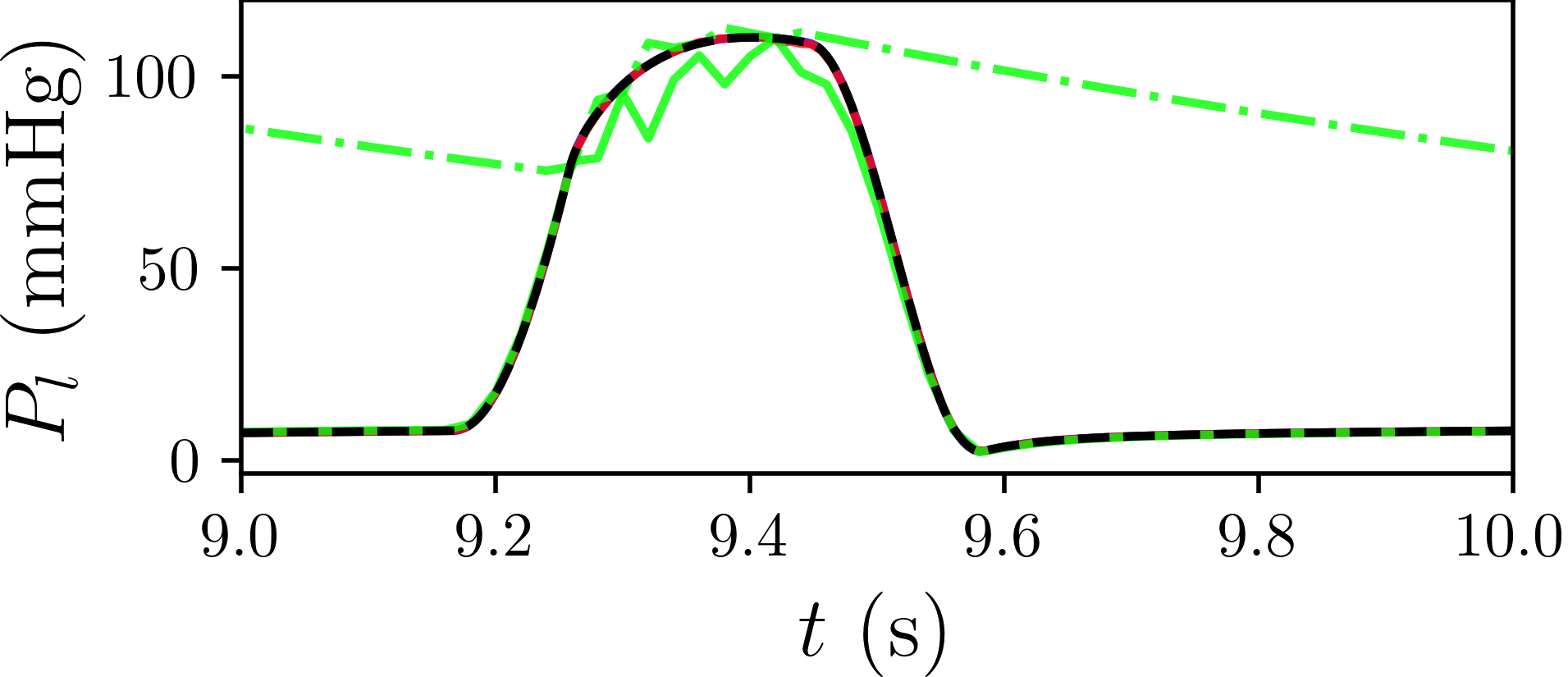}
    \end{subfigure} 
    \begin{subfigure}[b]{0.49\textwidth}
        \centering
         \includegraphics[scale=0.56]{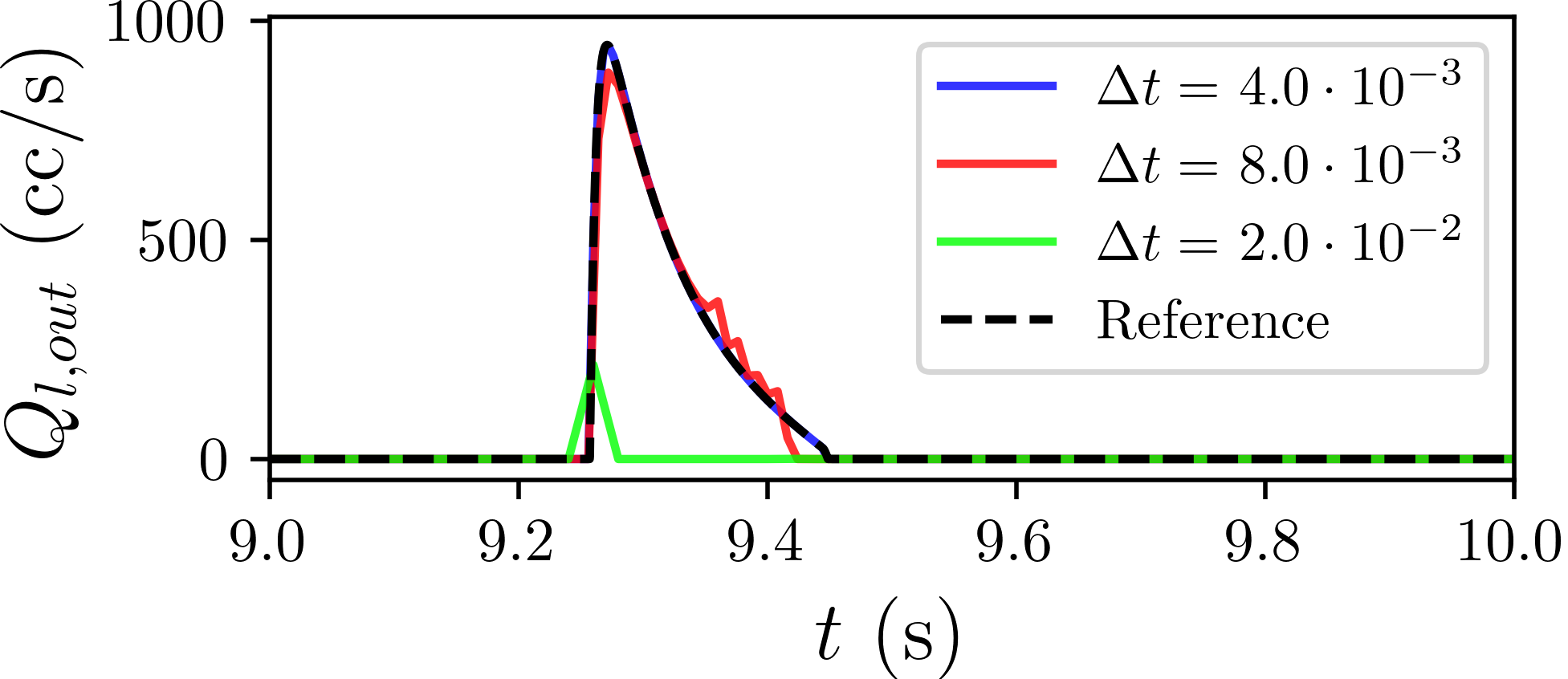}
    \end{subfigure}\\
    \vspace{0.15cm}
    \begin{subfigure}[b]{0.49\textwidth}
        \centering
         \includegraphics[scale=0.56]{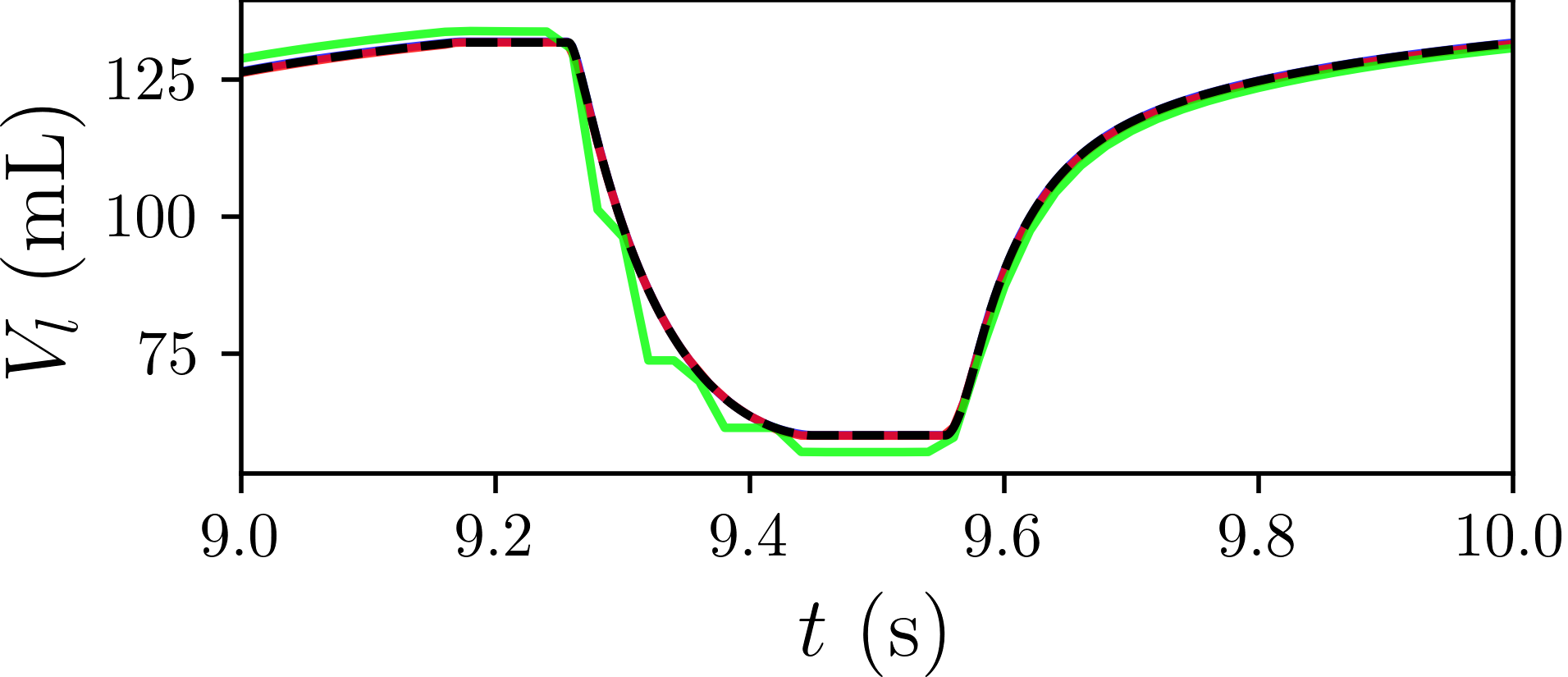}
    \end{subfigure}
    \begin{subfigure}[b]{0.49\textwidth}
        \centering
         \includegraphics[scale=0.56]{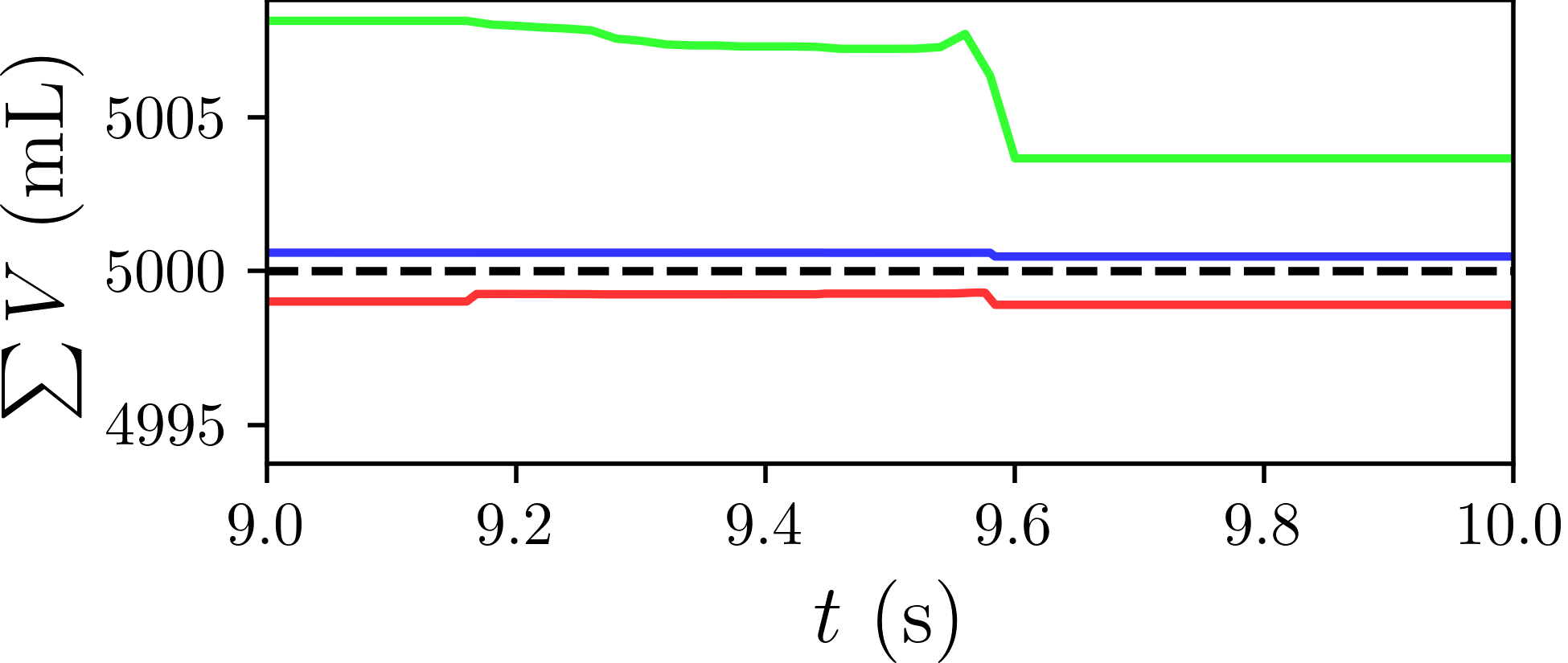}
    \end{subfigure}
    \caption{Comparison between numerical approximations for left ventricular pressure ($P_l$), outflow ($Q_{l,out}$), volume ($V_l$) and total stressed volume ($\sum V$). 
    We compare an explicit 4-th order Runge-Kutta method with three different time step sizes equal to $\Delta t= 4\cdot 10^{-3}$ (s), $8\cdot 10^{-3}$ (s), and $2\cdot 10^{-2}$ (s). The reference approach is an implicit 5-th order Runge-Kutta method (Radau~\cite{wanner1996solving}) with adaptive step size. 
    The systemic arterial pressure solution $P_a$ when $\Delta t= 2\cdot 10^{-2}$ (s) is also super-imposed in the top-left figure as the green dash-dotted line.}
    \label{fig:cvsim6-RK4-explicit-implicit}
\end{figure}

Finally, Figure~\ref{fig:cvsim6-RK4-explicit-implicit} shows how the use of an inappropriate time step size can trigger numerical instability and a non-physiological response. Specifically, when $\Delta t = 2\cdot 10^{-2}$ (s), an inaccurate approximation of the left ventricular pressure $P_l$ makes it lower than the systemic arterial pressure $P_a$ during systole (where a high stiffness ratio is observed), nearly causing the aortic valve to close.
The figure also shows that an appropriate time step and numerical algorithms are needed to ensure conservation of total stressed volume $\sum V$ (see also Appendix~\ref{sec: cvsim6-details}).

In conclusion, stiffness in ODE systems can pose challenges in both data-driven forward emulation and inverse modeling, since simulated data of poor quality may result from situations where the time step size or the numerical algorithm are not carefully selected, as shown in Figure~\ref{fig:cvsim6-RK4-explicit-implicit}. Also, an inadequate numerical solution can inject practical non-identifiability into the inference problem, even if no measurement noise is present in the data acquisition process.

\begin{remark}
    \label{rmk: sr tol}
    The stiffness of an ODE system can be overestimated by SR. For example, if $|\lambda_{\max}| = 0.1$ and $|\lambda_{\min}| = 1\cdot 10^{-11}$, SR $= 1\cdot 10^{10}$ while the fastest timescale $\tau = 1/|\lambda_{\max}|$ remains fairly large. 
    Besides, computation of small eigenvalues is affected by round off error, numerical solver accuracy, etc.
    An alternative criterion to measure stiffness of an ODE system is through the Lipschitz constant~\cite{shi2022polynomial,lambert1991numerical}, which is bounded below by the spectral radius, i.e., the magnitude of the largest eigenvalue for the RHS Jacobian.
    Therefore, as discussed above, we set a certain tolerance for the smallest eigenvalue to compute the stiffness ratio, but also specifically pay attention to the largest eigenvalue alone when characterizing stiffness.
\end{remark}

\begin{remark}
    As a result of assuming perfect synchronization between the left and right ventricular time-varying capacitance in Equation~\eqref{equ: driversfunction}, the duration of systole and diastole are expected to be identical for the two ventricles.
    However, due to differences in ventricular pressure and volume, and numerical accuracy, a small delay can occur between heart cycles as seen from the left and right ventricles (e.g., see Figure~\ref{fig: VP-eigSR}).
    This delay may depend on the specific input parameter combination, but in general is very limited. Therefore, for simplicity, we do not take it into consideration when calculating the right ventricular end-diastolic pressure $P_{r,edp}$. 
\end{remark}

\subsection{Output selection}\label{sec: cvsim6-output}

The selection of output attributes $\by \in \mathbb{R}^{16}$ for the CVSim-6 system (see Table~\ref{table:cvsim6 output paras}) is based on the EHR dataset~\cite{harrod2021predictive} to be studied in Section~\ref{sec: cvsim6-ehr}. 
This dataset includes patient-specific measurements of heart rate, systolic and diastolic blood pressures, mean venous pressure, left ventricular volumes, cardiac output, and vascular resistance.
Corresponding model outputs are computed from the CVSim-6 solutions over the final three of 12 simulated heart cycles (we refer to these three cycles as $\mathbf{T}_p$ with size $n_p$) due to a fully converged periodic response.

Note that the pulmonary wedge pressure $P_w$ and central venous pressure $P_{cvp}$ are indirect measures for the left and right atrial pressures~\cite{klabunde2011cardiovascular}, respectively, despite the absence of atrial compartments in CVSim-6.
Furthermore, the stroke volume (SV) is calculated as the difference between diastolic and systolic volumes in the left ventricle, i.e., $\text{SV} = V_{l,dia} - V_{l,sys}$ and from which, the cardiac output can be alternatively estimated as $\text{CO} = \text{SV} \times Hr$ (see, e.g.,~\cite{klabunde2011cardiovascular}).
Finally, to account for measurement noise, we use standard deviations for each output component as reported in Table~\ref{table:cvsim6 output paras}, based on values suggested in the literature~\cite{harrod2021predictive,gordon1983reproducibility,maceira2006normalized,yared2011pulmonary}.

\begin{table}[ht!]
{\small
\begin{center}
\begin{tabular}{@{} l l l l l @{}}
\toprule
Num. & Description & Calculation & Unit & Std.\\
\midrule
1. & Heart rate ($Hr$)  & Copied from input  & (bpm) & 3.0\\ 
\vspace{-0.25cm}\\
2. & Systolic blood pressure ($P_{a,sys}$)  & $\max_{t\in \mathbf{T}_p} P_a(t)$  & (mmHg) & 1.5\\ 
\vspace{-0.25cm}\\
3. & Diastolic blood pressure ($P_{a,dia}$)  & $\min_{t\in \mathbf{T}_p} P_a(t)$  & (mmHg) & 1.5\\ 
\vspace{-0.25cm}\\
4. & Right ventricular systolic pressure ($P_{r,sys}$) & $\max_{t\in \mathbf{T}_p} P_r(t)$ & (mmHg) & 1.0\\
\vspace{-0.25cm}\\
5. & Right ventricular diastolic pressure ($P_{r,dia}$) & $\min_{t\in \mathbf{T}_p} P_r(t)$ & (mmHg) & 1.0\\
\vspace{-0.25cm}\\
6. & Pulmonary arterial systolic pressure ($P_{pa,sys}$) & $\max_{t\in \mathbf{T}_p} P_{pa}(t)$ & (mmHg) & 1.0\\
\vspace{-0.25cm}\\
7. & Pulmonary arterial diastolic pressure ($P_{pa,dia}$) & $\min_{t\in \mathbf{T}_p} P_{pa}(t)$ & (mmHg) & 1.0\\
\vspace{-0.25cm}\\
8. & Right ventricular end diastolic pressure ($P_{r,edp}$) & $P_r\big(\mathbf{T}_p[n_p]\big)$ & (mmHg) & 1.0\\
\vspace{-0.25cm}\\
9. & Pulmonary wedge pressure ($P_{w}$) & $\frac{1}{|\mathbf{T}_p|}\int_{\mathbf{T}_p} P_{pv}(t) \ dt$ & (mmHg) & 1.0\\
\vspace{-0.25cm}\\
10. & Central venous pressure ($P_{cvp}$) & $\frac{1}{|\mathbf{T}_p|}\int_{\mathbf{T}_p} P_{v}(t) \ dt$ & (mmHg) & 0.5\\
\vspace{-0.25cm}\\
11. & Systolic left ventricular volume ($V_{l,sys}$) & $\min_{t\in \mathbf{T}_p} V_{l}(t)$ & (mL) & 10.0\\
\vspace{-0.25cm} \\
12. & Diastolic left ventricular volume ($V_{l,dia}$) & $\max_{t\in \mathbf{T}_p} V_{l}(t)$ & (mL) & 20.0\\
\vspace{-0.25cm} \\
13. & Left ventricular ejection fraction (LVEF) &  $(V_{l,dia} - V_{l,sys})/{V_{l,dia}}$ & $-$ & 0.02\\
\vspace{-0.25cm} \\
14. & Cardiac output (CO) & $\frac{1}{|\mathbf{T}_p|}\int_{\mathbf{T}_p} Q_{a}(t) \ dt$ & (L/min) & 0.2 \\
\vspace{-0.25cm}\\
15. & Systemic vascular resistance (SVR) &  $\big(\frac{1}{|\mathbf{T}_p|}\int_{\mathbf{T}_p} P_{a}(t) \ dt - P_{cvp}\big)/\text{CO}$ &  (dyn$\cdot$s/cm$^5$) & 50.0 \\
\vspace{-0.25cm}\\
16. & Pulmonary vascular resistance (PVR) & $\big(\frac{1}{|\mathbf{T}_p|}\int_{\mathbf{T}_p} P_{pa}(t) \ dt - P_{w}\big)/\text{CO}$ & (dyn$\cdot$s/cm$^5$) & 5.0\\
\bottomrule
\end{tabular}
\end{center}
\caption{CVSim-6 system outputs ($\by$) and corresponding uncertainty.}
\label{table:cvsim6 output paras}}
\end{table}

\subsection{Structural non-identifiability analysis with synthetic data}\label{sec: cvsim6-str}

We first investigate structural non-identifiability in CVSim-6 using noiseless synthetic data.
A set of $N$=54,000 input-output pairs is generated by numerically solving the CVSim-6 system, and then subdivided in 37,500 training data, 12,500 online testing data, and 4,000 offline validation data.
Input parameters are randomly selected from a uniform prior, centered around the default parameter combination (see Table~\ref{table:cvsim6 input paras}).
Specifically, the compliance parameters are allowed to vary by up to $\pm 50\%$, and other input parameters can fluctuate up to $\pm 30\%$, relative to their default value.
For brevity, we summarize results of the forward emulation and output density estimation in Appendix~\ref{sec: cvsim6-emu-nf} and focus on the inverse problem in this section. The choice of hyperparameters and other training details can be found in Appendix~\ref{sec: NN-details-cvsim6-str}.

We first validate the ability of an inVAErt network to solve amortized inverse problems by checking if a given output $\by$ can be reconstructed via the decoded parameters $\widehat{\bv}$.
To do so, we extract all 4,000 labels from the validation dataset and sample an equal number of latent variables from the standard normal, i.e., $\bw\sim \mathcal{N}(\boldsymbol{0}, \mathbf{I})$.
Then, concatenated outputs and latent space samples $[\by, \bw]^T$ are fed to the trained decoder $NN_d$, resulting in 4,000 inverse predictions $\widehat{\bv}$.
Finally, we evaluate this parameters through the exact CVSim-6 simulator, i.e., $\widehat{\by} = f(\widehat{\bv})$, and record the absolute reconstruction difference $|\by - \widehat{\by}|$ of each output component in Table~\ref{table:inverse-recons-component}.
Errors are very limited on average. In particular, we observe a maximum absolute error in the systolic arterial pressure equal to 3.4 mmHg, due to either an output $\by$ or a latent variable $\bw$ belonging to the tail of the respective distribution.

\begin{table}[ht!]
{\small
\begin{center}
\begin{tabular}{@{} l c c c c c @{}}
\toprule
& $Hr$ (bpm) & $P_{a,sys}$ (mmHg) & $P_{a,dia}$ (mmHg) & $P_{r,sys}$ (mmHg) & $P_{r,dia}$ (mmHg) \\
\midrule
Average & 4.62$\cdot 10^{-2}$ & 3.37$\cdot 10^{-1}$ & 2.50$\cdot 10^{-1}$ & 8.16$\cdot 10^{-2}$ & 1.12$\cdot 10^{-2}$ \\
Max     & 4.06$\cdot 10^{-1}$ & 3.34$\cdot 10^{+0}$ & 2.42$\cdot 10^{+0}$ & 1.60$\cdot 10^{+0}$ & 1.41$\cdot 10^{-1}$ \\
Std     & 4.36$\cdot 10^{-2}$ & 3.32$\cdot 10^{-1}$ & 2.36$\cdot 10^{-1}$ & 8.27$\cdot 10^{-2}$ & 1.20$\cdot 10^{-2}$ \\
\bottomrule
\end{tabular}
\bigskip

\begin{tabular}{@{} l c c c c c @{}}
\toprule
 & $P_{pa,sys}$ (mmHg) & $P_{pa,dia}$ (mmHg) & $P_{r,edp}$ (mmHg) & $P_{w}$ (mmHg) & $P_{cvp}$ (mmHg) \\
\midrule
Average & 8.06$\cdot 10^{-2}$ & 5.12$\cdot 10^{-2}$ & 1.82$\cdot 10^{-2}$ & 4.89$\cdot 10^{-2}$ & 2.72$\cdot 10^{-2}$ \\
Max     & 1.58$\cdot 10^{+0}$ & 1.08$\cdot 10^{+0}$ & 2.78$\cdot 10^{-1}$ & 1.12$\cdot 10^{+0}$ & 3.26$\cdot 10^{-1}$ \\
Std     & 8.16$\cdot 10^{-2}$ & 5.48$\cdot 10^{-2}$ & 1.93$\cdot 10^{-2}$ & 5.20$\cdot 10^{-2}$ & 2.71$\cdot 10^{-2}$ \\
\bottomrule
\end{tabular}
\bigskip

\begin{tabular}{@{} l c c c c c c @{}}
\toprule
 & $V_{l,sys}$ (mL) & $V_{l,dia}$ (mL) & LVEF (-) & CO (L/min) & SVR (dyn$\cdot$s/cm$^5$) & PVR (dyn$\cdot$s/cm$^5$) \\
\midrule
Average & 1.49$\cdot 10^{-1}$ & 3.29$\cdot 10^{-1}$ & 5.04$\cdot 10^{-4}$ & 1.41$\cdot 10^{-2}$ & 1.02$\cdot 10^{+0}$ & 8.79$\cdot 10^{-2}$ \\
Max     & 1.79$\cdot 10^{+0}$ & 4.67$\cdot 10^{+0}$ & 6.95$\cdot 10^{-3}$ & 2.22$\cdot 10^{-1}$ & 1.20$\cdot 10^{+1}$ & 6.54$\cdot 10^{-1}$ \\
Std     & 1.57$\cdot 10^{-1}$ & 3.41$\cdot 10^{-1}$ & 5.09$\cdot 10^{-4}$ & 1.44$\cdot 10^{-2}$ & 9.32$\cdot 10^{-1}$ & 8.07$\cdot 10^{-2}$ \\
\bottomrule
\end{tabular}
\bigskip
\caption{Statistics of the absolute reconstruction error across the validation dataset of size 4,000.}
\label{table:inverse-recons-component}
\end{center}}
\end{table}

\subsubsection{Going beyond point estimates}\label{sec: cvsim6-str-go-beyond-ps}

In the previous section, each output $\by$ was associated with a single latent variable realization $\bw$, thereby generating \emph{a single} inverse problem solution.
However, structurally non-identifiable systems are characterized by the possibility that infinitely many inputs $\bv \in \bcM_{\by} \subset \bcV$ may map to a common output $\by$, i.e., $f(\bv) = \by$ for any $\bv\in \bcM_{\by}$.
For an inVAErt network, discovering multiple solutions from the non-identifiable manifold $\bcM_{\by}$ is as easy as drawing more latent samples from $\bcW$, and perform one decoder evaluation for each of these samples.
To test this, we first select a representative output from the trained density estimator $NN_f$, denoted as $\by^*$, with its components listed in Table~\ref{table:inverse-id-check-point}.

\begin{table}[ht!]
{\small
\begin{center}
\begin{tabular}{@{} c c c c c c @{}}
\toprule
$Hr$ (bpm) & $P_{a,sys}$ (mmHg) & $P_{a,dia}$ (mmHg) & $P_{r,sys}$ (mmHg) & $P_{r,dia}$ (mmHg) & $P_{pa,sys}$ (mmHg) \\
\midrule
72.91 & 142.78 & 111.55 & 22.78 & -1.51 & 22.37\\
\bottomrule
\end{tabular}
\bigskip

\begin{tabular}{@{} c c c c c @{}}
\toprule
$P_{pa,dia}$ (mmHg) & $P_{r,edp}$ (mmHg) & $P_{w}$ (mmHg) & $P_{cvp}$ (mmHg) & $V_{l,sys}$ (mL) \\
\midrule
13.13 & 1.61 & 12.56 & 7.27 & 76.81 \\
\bottomrule
\end{tabular}
\bigskip

\begin{tabular}{@{} c c c c c @{}}
\toprule
$V_{l,dia}$ (mL) & LVEF (-) & CO (L/min) & SVR (dyn$\cdot$s/cm$^5$) & PVR (dyn$\cdot$s/cm$^5$) \\
\midrule
153.61 & 0.50 & 5.60 & 1725.48 & 74.53 \\
\bottomrule
\end{tabular}
\bigskip
\end{center}}
\caption{Selected CVSim-6 output $\by^*$ sampled from $NN_f$.}
\label{table:inverse-id-check-point}
\end{table}

\begin{figure}[!ht]
\begin{subfigure}[b]{0.495\textwidth}
        \centering
         \includegraphics[scale=0.43]{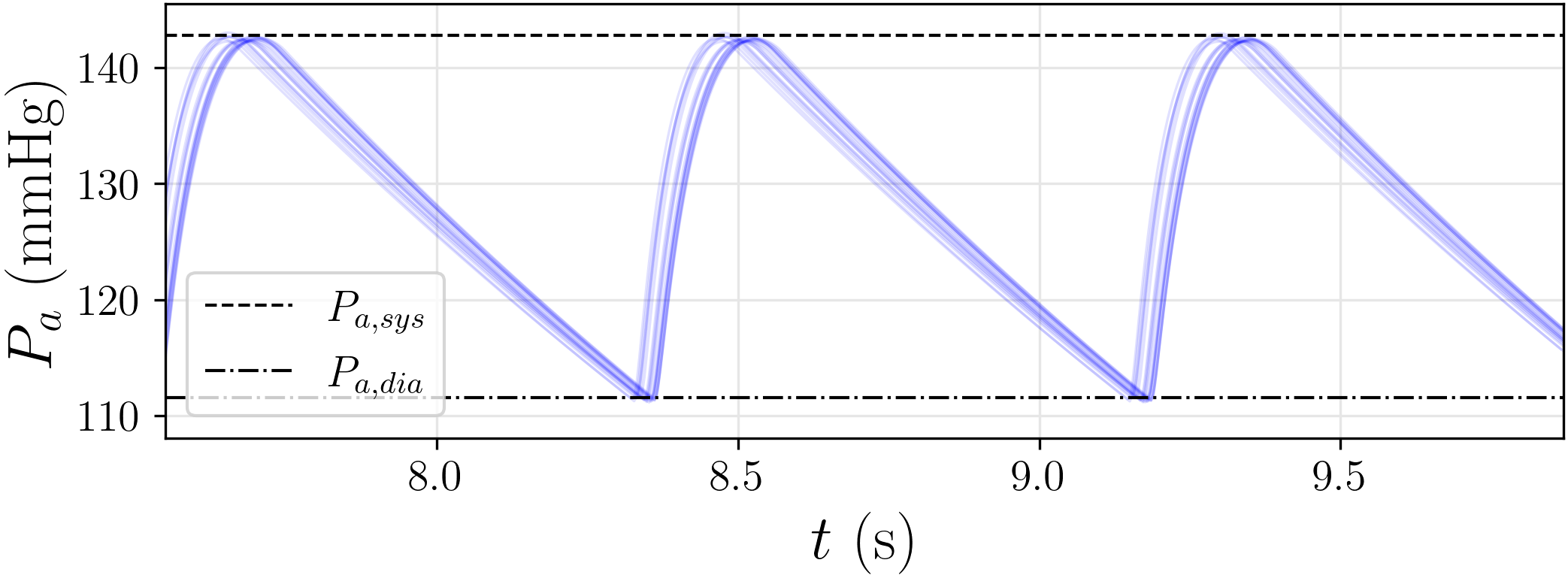}
     \end{subfigure}
    \begin{subfigure}[b]{0.495\textwidth}
        \centering
         \includegraphics[scale=0.43]{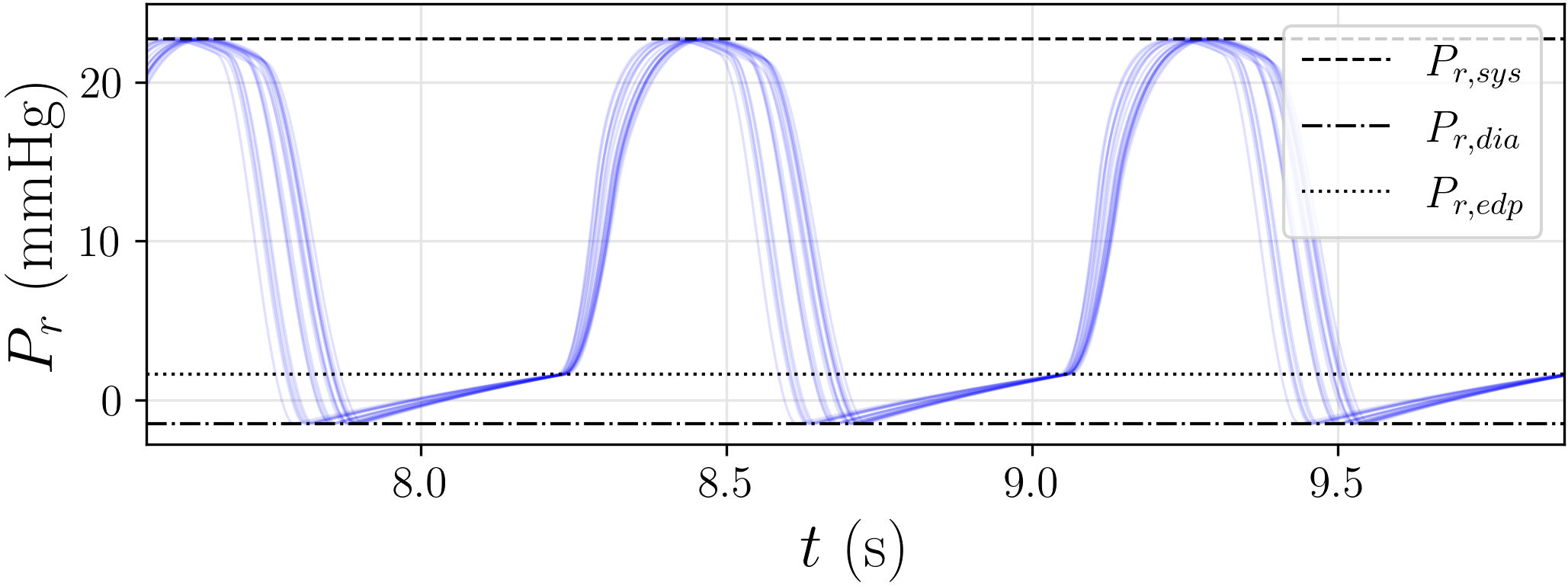}
    \end{subfigure}\\
    
    \begin{subfigure}[b]{0.495\textwidth}
        \centering
         \includegraphics[scale=0.43]{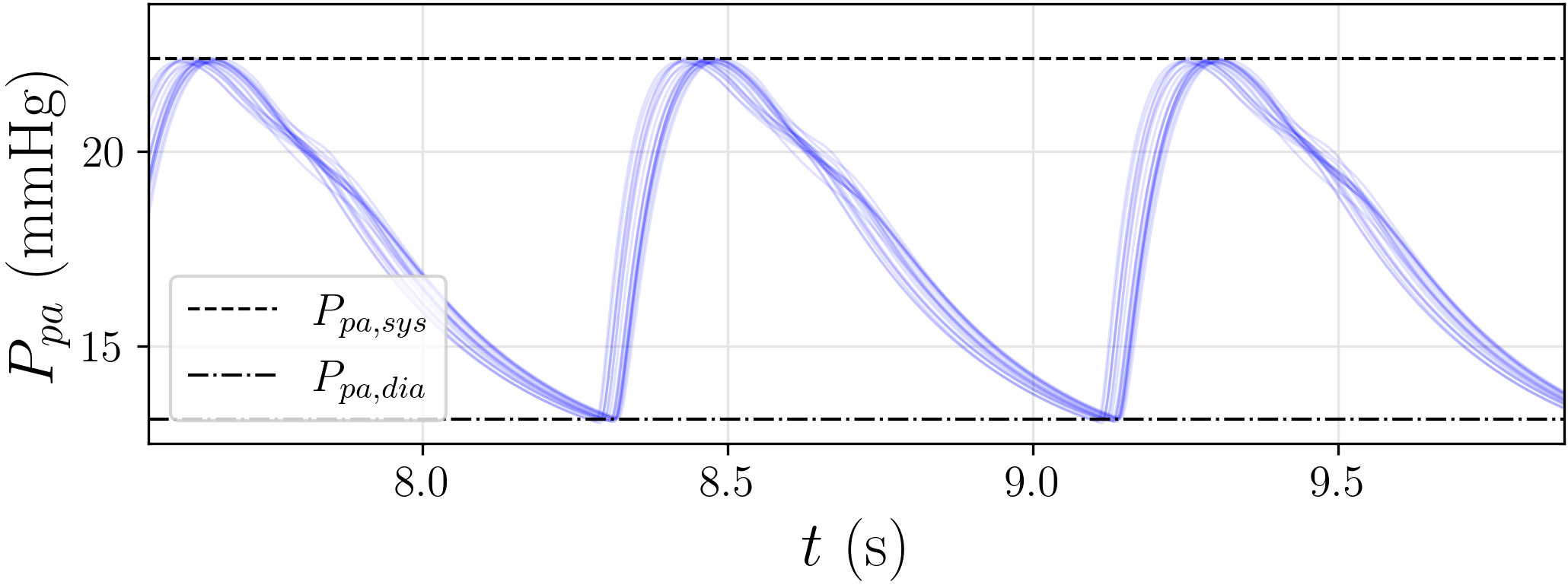}
     \end{subfigure}
    \begin{subfigure}[b]{0.495\textwidth}
        \centering
         \includegraphics[scale=0.43]{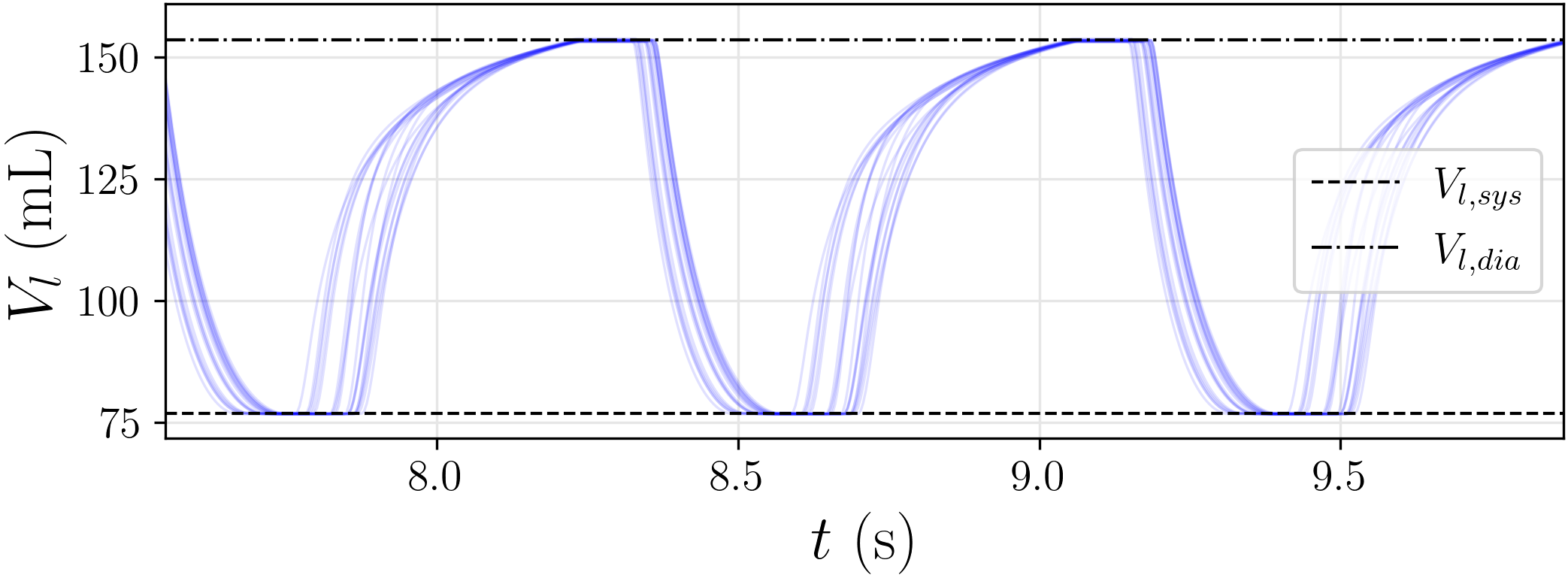}
    \end{subfigure}
    \caption{Trajectories of systemic arterial pressure $P_a$, right ventricular pressure $P_r$, pulmonary arterial pressure $P_{pa}$, and left ventricular volume $V_l$ corresponding to different latent space samples, computed from the observation $\by^*$ (see Table~\ref{table:inverse-id-check-point}). To better visualize each reconstruction, we plot only 20/100 pressure and volume time traces, together with systolic/diastolic targets in $\by^*$ using horizontal lines.}
    \label{fig: pv-check-str-id}
\end{figure}

For this output, we draw 100 samples $\bw$ from the latent space and process each $[\by^*, \bw]^T$ through $NN_d$, generating 100 predictions of $\widehat{\bv}$.  
To confirm these predictions are indeed associated with the non-identifiable manifold $\bcM_{\by^*}$, we use them to reconstruct $\by^*$ via the exact numerical simulator $f$, achieving a maximum relative error of approximately 0.139\%.
When compared to the components in $\by^*$ as presented in Figure~\ref{fig: pv-check-str-id}, the pressure/volume trajectories computed through $\widehat{\bv}$, are able to accurately reproduce both systolic and diastolic targets.

\subsubsection{Analysis of non-identifiable manifolds}\label{sec: cvsim6-str-manifold}

\begin{figure}[!ht]
\begin{subfigure}[b]{0.999\textwidth}
        \centering
         \includegraphics[scale=0.38]{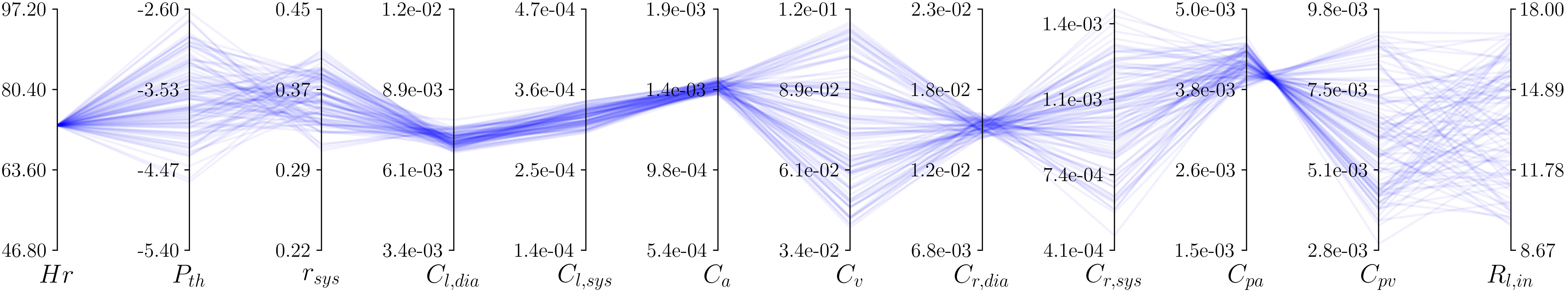}
     \end{subfigure}\\
     \vspace{-0.15cm}\\
\begin{subfigure}[b]{0.999\textwidth}
        \centering
         \includegraphics[scale=0.38]{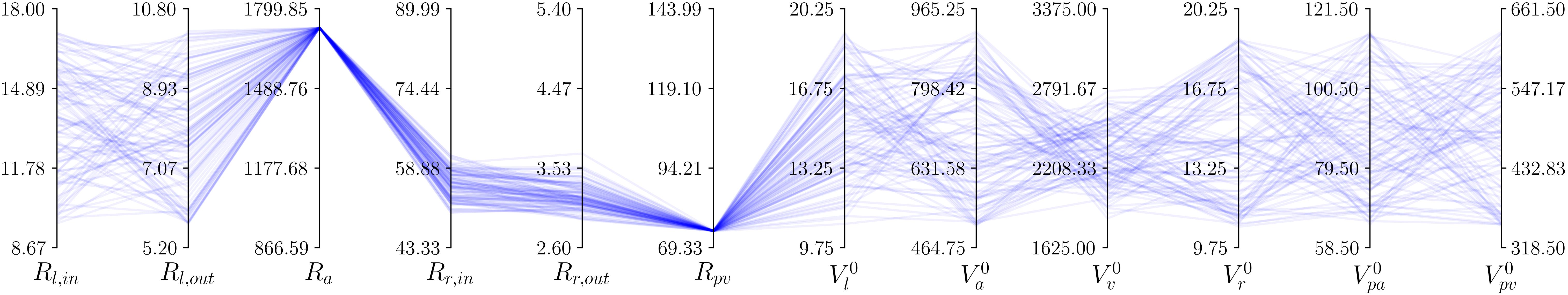}
     \end{subfigure}
     \caption{Parallel coordinate plot of the non-identifiable manifold $\bcM_{\by^*}$, obtained by inverting the fixed observation $\by^*$ (Table~\ref{table:inverse-id-check-point}), together with 100 latent variable realizations $\bw$ drawn from a standard normal. Plot limits for each parameter are $\pm 5\%$ larger than the bounds defined in Section~\ref{sec: cvsim6-str}, with units omitted for brevity (see Table~\ref{table:cvsim6 input paras}). The parameter $R_{l,in}$ is repeated in both rows to include all connections.}
    \label{fig: varcheck-str-id}
\end{figure}

A parallel coordinate plot of the manifold $\bcM_{\by^*}$ is presented in Figure~\ref{fig: varcheck-str-id}.
First, a significant variability can be observed in some of the input components, confirming that the CVSim-6 system is structurally non-identifiable.
From the picture, parameters such as $C_v, C_{r,sys}, C_{pv}$, $R_{l,in}, R_{l,out}$ and most of the unstressed volumes appear to be mostly responsible for the lack of identifiability, as their predicted values cover almost entirely their prior ranges.
However, input parameters such as $R_a$, $R_{pv}$ show very limited variations. This suggests the non-identifiable manifold $\bcM_{\by^*}$ embedded in the 23-dimensional input space $\bcV$ may be characterized by a smaller intrinsic dimension, and, in turn, that some of the parameters are actually identifiable.

\begin{figure}[!ht]
\begin{subfigure}[b]{0.495\textwidth}
        \centering
         \includegraphics[scale=0.4]{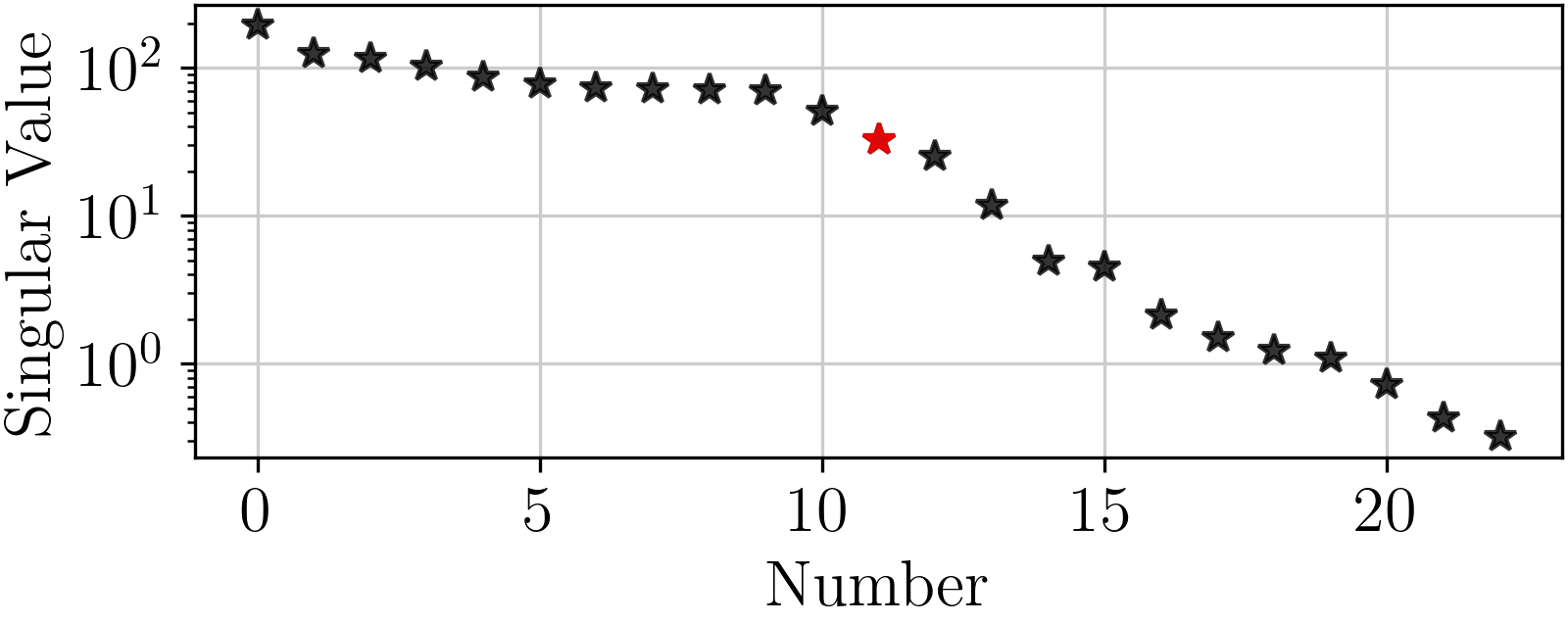}
     \end{subfigure}
    \begin{subfigure}[b]{0.495\textwidth}
        \centering
         \includegraphics[scale=0.4]{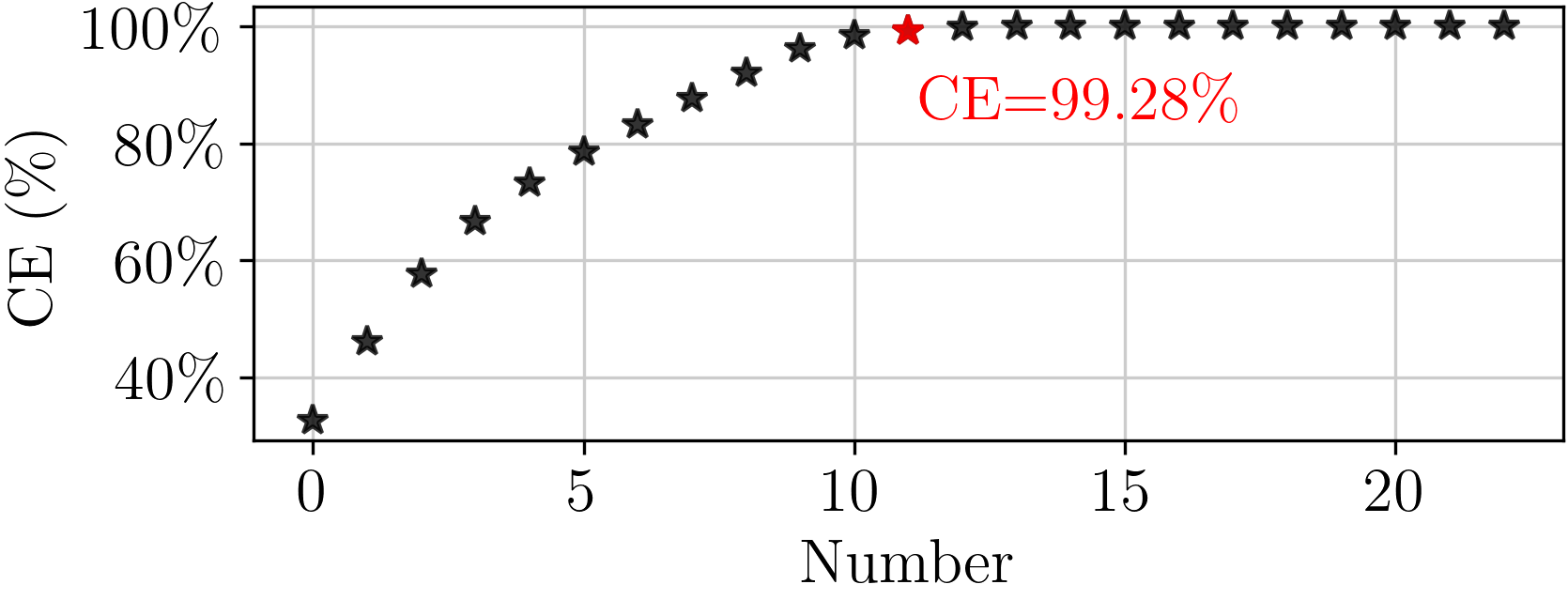}
    \end{subfigure}
    \caption{Singular value spectrum and cumulative energy distribution of 5,000 input parameters belonging to $\bcM_{\by^*}$.}
    \label{fig: svd-str-id}
\end{figure}

To discover the lower dimensional structure embedded in $\bcM_{\by^*}$, we compute the singular values of the data matrix sampled from $\bcM_{\by^*}$, and the associated \emph{cumulative energy} (CE)~\cite{leskovec2020mining, brunton2022data} (See Figure~\ref{fig: svd-str-id}).
CE measures how much data variance is restored by keeping the first $n$ principal components. It is the ratio of the sum of the first $n$ squared singular values, to the sum of the squares of all singular values.
As illustrated in Figure~\ref{fig: svd-str-id}, 12 modes are sufficient to reconstruct 99.28\% of the total \emph{energy}~\cite{leskovec2020mining}.

\subsubsection{Missing data analysis}\label{sec: cvsim6-str-missing}

Consider now the situation where the components $\{ P_{a, dia}, P_{r,dia}, P_{pa,dia}, P_w, V_{l,sys}, \text{SVR}, \text{PVR}\}$ of $\by^*$ are missing, e.g., unobserved. 
We wonder whether one can infer these missing elements through the pre-trained density estimator $NN_f$, and then solve the inverse problem using the decoder network $NN_d$.
We answer this question using Algorithm~\ref{alg: missing data inference}.
The ranking mechanism in Algorithm~\ref{alg: missing data inference} promotes certain combinations of components that are most frequently observed in the solution of the CVSim-6 system, providing an approach for physics-based multiple imputation (see, e.g.,~\cite{sterne2009multiple}).

\begin{algorithm*}[ht!]
\caption{Inversion with missing data}
\label{alg: missing data inference}
\begin{algorithmic}
\State \textbf{Inputs:}
\begin{itemize}
    \item An observation $\boldsymbol{y}^*$ with missing components.
    \item The set $\chi$ with the indices of the missing components in $\boldsymbol{y}^*$, and its complement $\chi^c$.
    \item Optimally trained density estimator $NN_f$ and decoder $NN_d$.
\end{itemize}
\State \textbf{Output:}
\begin{itemize}
    \item Inverse predictions: $\widehat{\bv}$.
\end{itemize}
\State \textbf{Steps:}
\begin{enumerate}
    \item Generate a large sample set $\mathbf{Y}_f \in \mathbb{R}^{M\times \dim(\by)}$ of size $M$ from $NN_f$. 
    \item Replace $\mathbf{Y}_f$'s components with those not missing in $\by^*$, i.e., $\mathbf{Y}_f[:, \chi^c] = \by^*[\chi^c]$.
    \item Evaluate the log-density for each modified sample in $\mathbf{Y}_f$ (denoted as $\by^c$) using normalizing flow
    \begin{equation*}\label{equ:change_var_nf}
    \log p(\by^c) = \log \pi_0(\bz) + \log \left|\det \frac{d \by^c}{d \bz}\right|^{-1}, \quad \bz = NN_f^{-1}(\by^c; \bphi_f^{\text{opt}}) \ .
    \end{equation*}
    \item Re-order $\mathbf{Y}_f$ by ranking the log-density of each $\by^c$ and take the most likely samples.
    \item For each selected $\by^c$, solve the inverse problem
    \begin{equation*}
        \widehat{\bv} = NN_d([\by^c, \bw]^T; \bphi_d^{\text{opt}}), \quad \bw \sim \boldsymbol{\mathcal{W}}.
    \end{equation*}
\end{enumerate}
\end{algorithmic}
\end{algorithm*}
%
%
%
%

\begin{figure}[!ht]
\begin{subfigure}[b]{0.495\textwidth}
        \centering
         \includegraphics[scale=0.43]{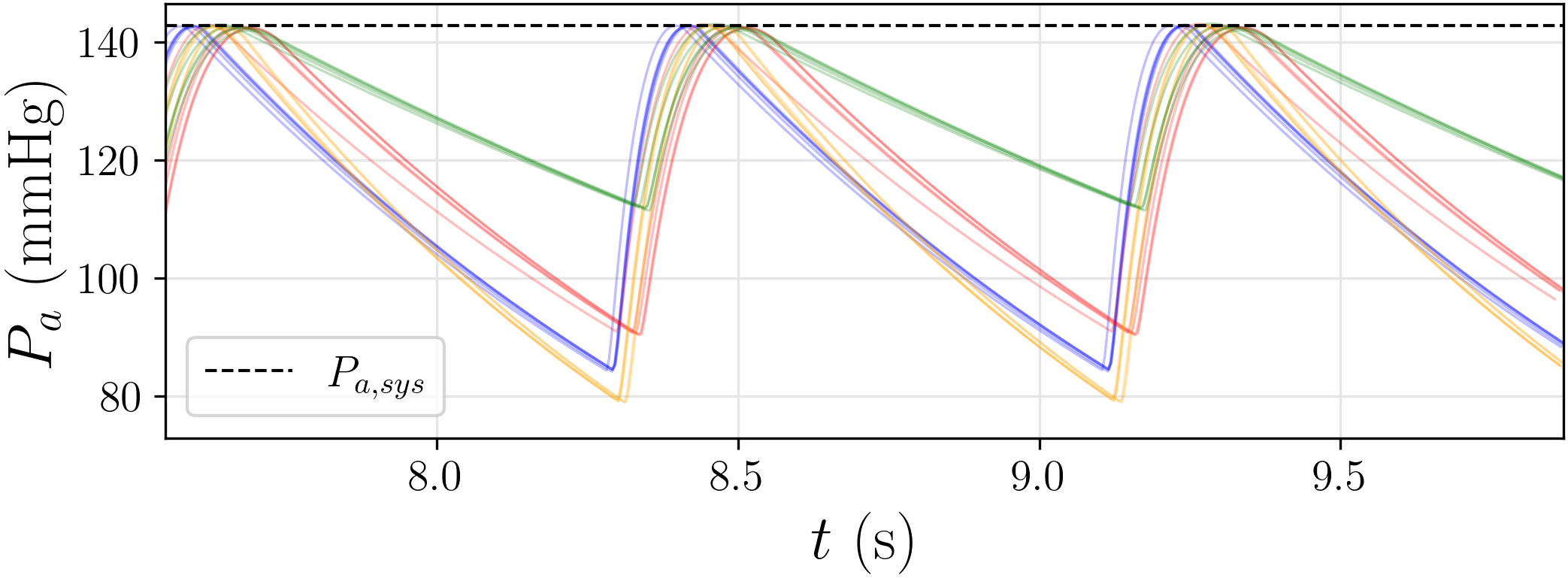}
     \end{subfigure}
    \begin{subfigure}[b]{0.495\textwidth}
        \centering
         \includegraphics[scale=0.43]{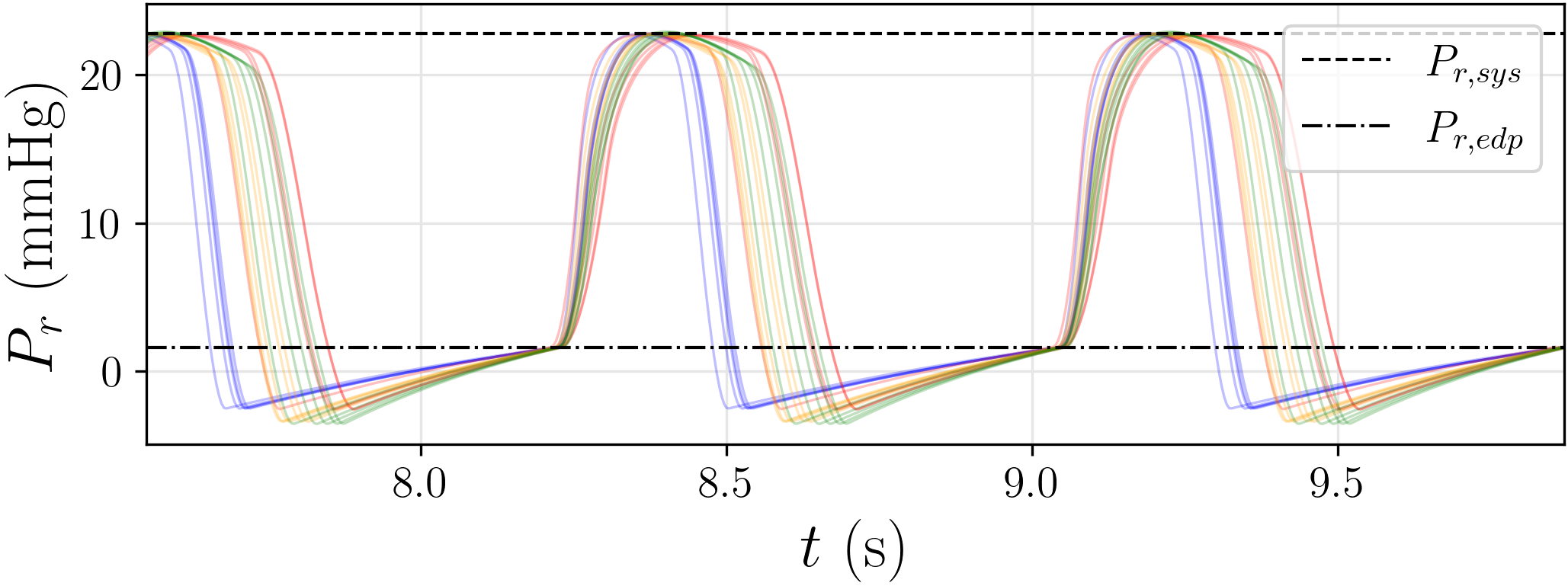}
    \end{subfigure}\\
    
    \begin{subfigure}[b]{0.495\textwidth}
        \centering
         \includegraphics[scale=0.43]{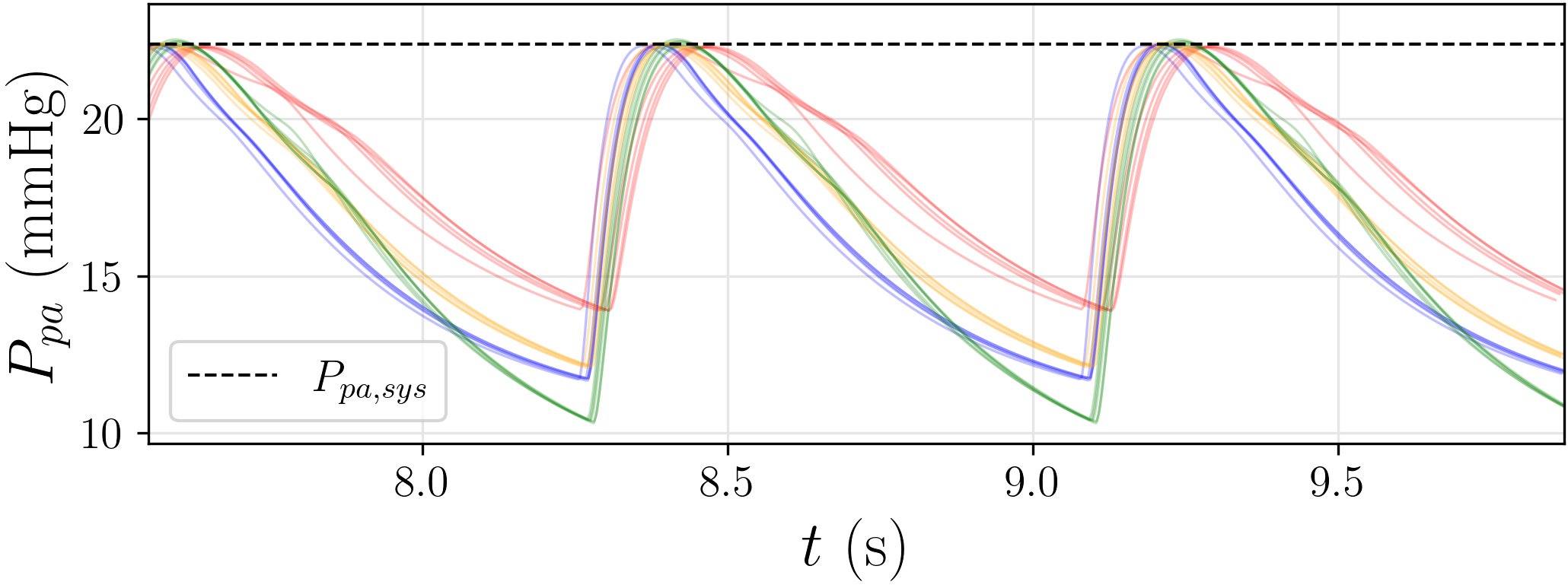}
     \end{subfigure}
    \begin{subfigure}[b]{0.495\textwidth}
        \centering
         \includegraphics[scale=0.43]{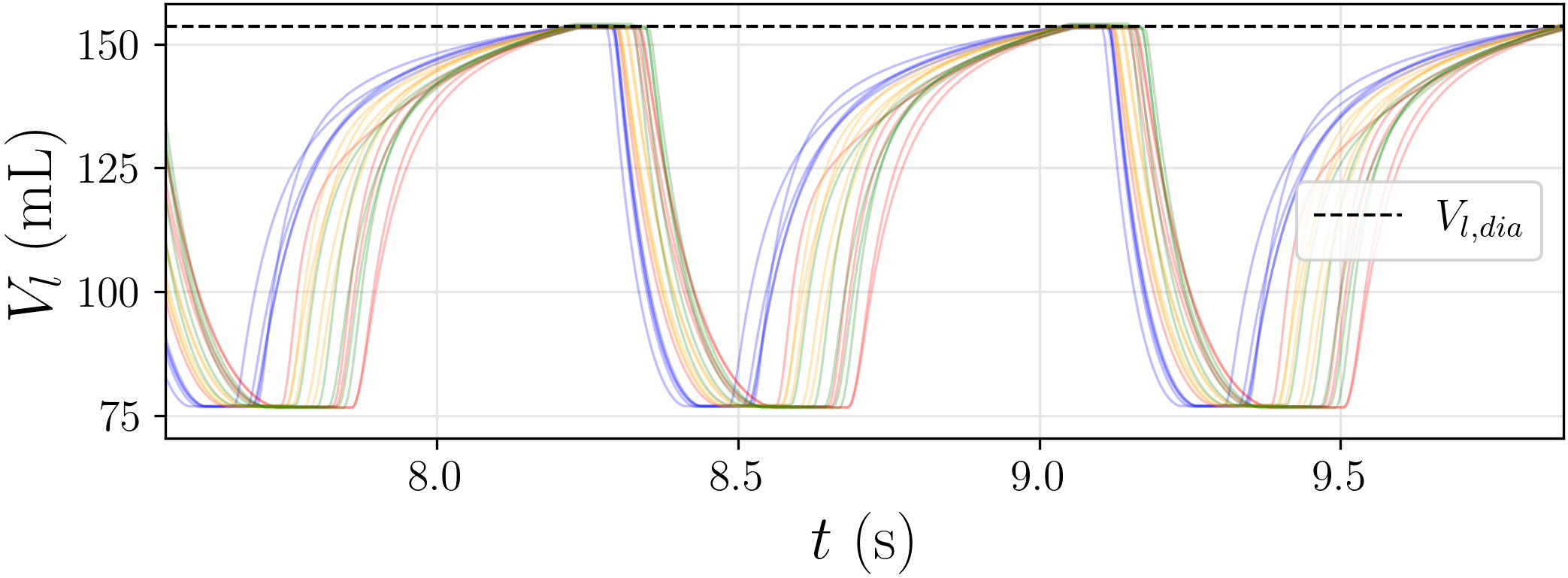}
    \end{subfigure}
    \caption{Trajectories of systemic arterial pressure $P_a$, right ventricular pressure $P_r$, pulmonary arterial pressure $P_{pa}$, and left ventricular volume $V_l$ obtained from $\widehat{\bv}$ under missing data.
    Four potential outputs are extracted based on the density ranking, with five latent samples of $\bw$ utilized for the inverse problem in each group. 
    Different colors are applied to each group, and the target (non missing) systolic/diastolic components of $\by^*$ are plotted using horizontal lines.}
    \label{fig: inverse-missing-curves}
\end{figure}

To show the accuracy of this approach, we remove the previously mentioned components (assumed missing) from $\by^*$, and extract the top four most probable outputs using Algorithm~\ref{alg: missing data inference}.
Next, for each output, we solve the inverse problem by sampling 5 latent variables $\bw$ from the standard normal, resulting in 4 groups of inverse predictions of $\widehat{\bv}$.
We then evaluate these parameters with the exact CVSim-6 simulator $f$, and plot the resulting pressure/volume curves in Figure~\ref{fig: inverse-missing-curves}, using a different color for each group.

First, the systolic and diastolic extrema for the trajectories in Figure~\ref{fig: inverse-missing-curves} align closely with the corresponding targets (non missing) in $\by^*$, indicating a high reconstruction accuracy (maximum relative reconstruction error based on the $4\times5=20$ predictions is about $0.408\%$).
Second, the missing diastolic components for the systemic and pulmonary arterial pressures in $\by^*$ are free to change between values with high likelihood. In contrast, the right ventricular diastolic pressures show only minor variations, and the left ventricular systolic volumes remain relatively consistent across each group. This behavior is specific for the selected $\by^*$ and may change with other observations and different missing components.
%

\subsection{Inverse problem with real clinical measurements}\label{sec: cvsim6-ehr}

We now consider data from a real EHR dataset provided through a Google's ATAP research project~\cite{harrod2021predictive, su2023privacy}, containing clinical measurements from 84 fully anonymized adult patients.
Each patient has up to 26 features but none of the patients has data for all the features. While the missing components vary among the cohort, attributes like heart rate, or systolic and diastolic blood pressures are available for almost all patients, since they can be easily measured.
As mentioned in Section~\ref{sec: cvsim6-output}, 16 out of these 26 clinical targets are selected as outputs for the CVSim-6 system in this paper. 

Unlike studies with synthetic data in the previous sections, inferring input parameters through these real patient-specific data is a challenging task, mainly due to (1) potential measurement noise, (2) missing data in each measurement, and  (3) simplifying assumptions and model misspecification in CVSim-6 with respect to the real cardiovascular physiology. 
In other words, an EHR measurement may be not belong to the CVSim-6 model range (a.k.a. out-of-distribution data).
As a result, the complexity of the inversion tasks may increase significantly due to combined structural and practical non-identifiability. 

To handle these challenges, we made three main adjustments. First, we expanded the bounds for the prior used to selected training inputs (compared to Section~\ref{sec: cvsim6-str}) such that the corresponding CVSim-6 outputs would cover most of the clinical measurements in the EHR dataset. 
Specifically, we allowed the heart rate to vary from $-20\%$ to $+60\%$, and each compliance and resistances to vary from $-80\%$ to $+60\%$, and the remaining parameters to vary $\pm 30\%$, with respect to their default values reported in Table~\ref{table:cvsim6 input paras}.
Second, we added noise during network training to enhance generalization. The noise for each output component is modeled as a zero-mean Gaussian, with standard deviation from the literature (see Table~\ref{table:cvsim6 output paras}). Third, we adjusted the neural network hyperparameters to avoid over-fitting.
Similar to Section~\ref{sec: cvsim6-str}, a total of 54,000 input-output pairs are numerically generated according to the new prior. Of these, 4,000 data will be used for offline validation, and the rest is divided into the training and online testing set with a 3-to-1 ratio.

For the EHR dataset, we slightly modify the notation, where the clinical data for the $q$-th patient $(1\leq q \leq 84)$ is denoted by $\by^{\text{EHR},(q)}$.
A diagram for the inference task where $N_w$ latent variable samples are used for each EHR measurement is illustrated in equation~\eqref{equ:ehr-inference} for the CVSim-6 model with input dimension $\dim(\bv) = 23$ and output dimension $\dim(\by) = 16$.
\begin{equation}
    \by^{\text{EHR},(q)} \ {\underset{\text{Algorithm~\ref{alg: missing data inference}}}{\xrightarrow{\hspace*{1cm}}}} \ \widehat{\boldsymbol{V}}^{\text{EHR},(q)} \in \mathbb{R}^{N_w\times 23} \ {\underset{NN_e}{\xrightarrow{\hspace*{0.5cm}}}} \ \widehat{\boldsymbol{Y}}^{\text{EHR},(q)} \in \mathbb{R}^{N_w\times 16} \ .
    \label{equ:ehr-inference}
\end{equation}
Note that the emulator $NN_e$ is employed to replace the exact CVSim-6 simulator during the output reconstruction process from $\widehat{\boldsymbol{V}}^{\text{EHR}}$ to $\widehat{\boldsymbol{Y}}^{\text{EHR}}$.
There are two reasons for this. First the emulator is significantly faster and second the exact numerical integrator may crash due to an input parameter combination resulting from an extremely poor inversion task, while the network model can still make a prediction without causing the code to terminate.
Moreover, to further justify the use of an emulator, the histogram for the $l_2$ norm of the prediction error on the validation dataset of size 4,000 is shown in Figure~\ref{fig:emulator-hist-ehr} below.
\begin{figure}[!ht]
    \centering
    \includegraphics[scale=0.35]{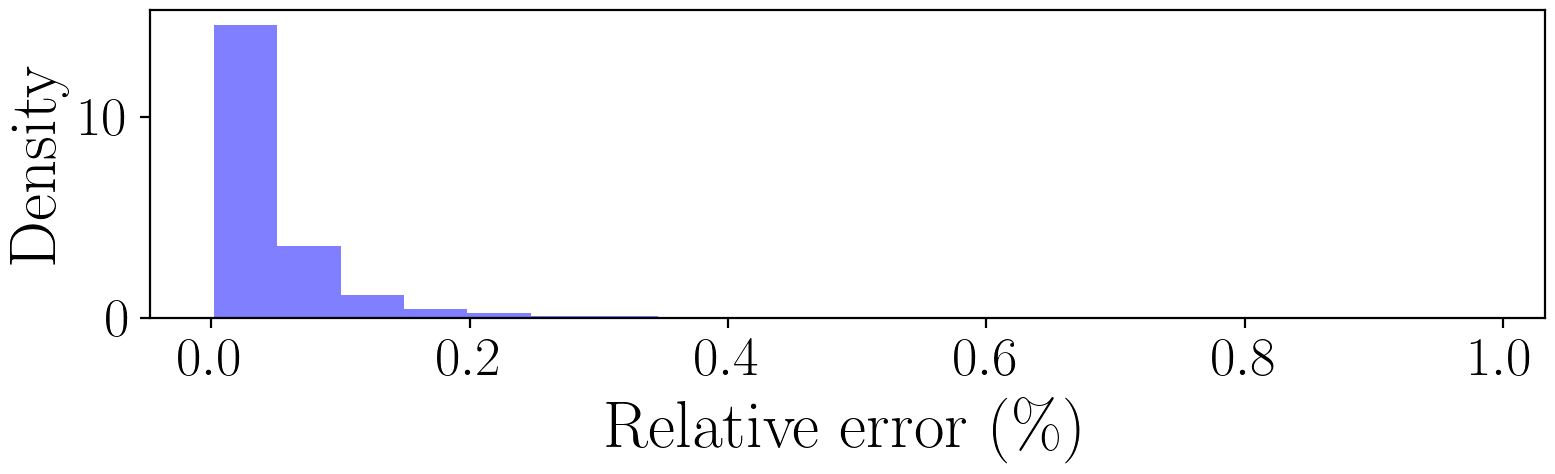}
    \caption{$l_2$ relative error distribution for the emulator. The maximum error on 4,000 validation examples is equal to 0.984\%.}
    \label{fig:emulator-hist-ehr}
\end{figure}

\subsubsection{Inversion results}\label{sec: cvsim6-ehr-training-noise}

In this section, we quantify how much noise should be added to the CVSim-6 outputs so an inVAErt network could achieve sufficient generalization on real data in general, and our EHR dataset in particular.
To do so, we introduce a scale factor $\delta$ which amplifies the values of standard deviation in Table~\ref{table:cvsim6 output paras}. 
This parameter quantifies the amount of noise added during training where, for example, $\delta=0.0$ corresponds to no noise added, and $\delta = 2.0$ is used to denote a zero-mean Gaussian noise with standard deviation twice as the value reported in Table~\ref{table:cvsim6 output paras}.

Additionally, an error criterion to quantify the accuracy of the inversion task on a per output component basis is defined in equation~\eqref{equ: error-noise-level} below 
\begin{equation}
   e_k = \frac{1}{|\mathcal{P}_k|} \cdot \frac{1}{N_w} \sum_{q \in \mathcal{P}_k} \sum_{i=1}^{N_w}\big|y^{\text{EHR},(q)}_k - \widehat{Y}_{ik}^{\text{EHR},(q)}\big|, \quad k=1:16 \ ,
   \label{equ: error-noise-level}
\end{equation}
where $\mathcal{P}_k$ contains the IDs of the patients whose $k$-th attribute is not missing, and $|\mathcal{P}_k|$ is the cardinality of this set, reported under the ``Counts'' column in Table~\ref{table:ehr-error-global}.
Inversion is performed with four noise levels ($\delta = 0.25, 0.5, 1.0, 2.0$), with results reported in Table~\ref{table:ehr-error-global}. The lowest errors are obtained for $\delta = 0.5$.
Note that we use four networks with the same hyperparameters, one for each noise intensity (see Section~\ref{sec: NN-details-cvsim6-EHR}), even though improved results may be obtained by fine tuning hyperparameters on training dataset for each noise level. 
Additionally, note that inversion is only performed for patients having more than 10 measurements, that is, 46 out of 84 patients. 
Also note that, generally speaking, the fewer missing components, the harder the inversion task becomes. This is because EHR measurements usually fall on the tail of the distribution learned from the synthetic CVSim-6 outputs, due to model misspecification and noise corruption. 
%

\begin{table}[ht!]
{\small
\begin{center}
\begin{tabular}{@{}l l c c c c c r@{}}
\toprule
Num. & Quantity  & $\delta=0.25$ & $\delta=0.5$ & $\delta=1.0$ & $\delta=2.0$ & Unit & Counts\\
\midrule
1. & $Hr$  & 0.602 & \b{0.440} & 0.843 & 2.636  & (bpm) & 46/46\\
\vspace{-0.25cm}\\
2. & $P_{a,sys}$   & 5.204  & \b{2.284} & 3.379 & 2.738 & (mmHg) & 46/46\\
\vspace{-0.25cm}\\
3. & $P_{a,dia}$  & 3.614 & \b{1.605} & 2.045 & 2.627 & (mmHg) & 46/46\\
\vspace{-0.25cm}\\
4. & $P_{r,sys}$  & 1.669 & \b{1.323} & 1.469 & 1.476 & (mmHg) & 44/46\\
\vspace{-0.25cm}\\
5. & $P_{r,dia}$  & \b{1.170} & 1.684 & 2.022 & 3.383 & (mmHg) & 11/46\\
\vspace{-0.25cm}\\
6. & $P_{pa,sys}$  & 1.907 & \b{1.227} & 1.606  & 1.386 & (mmHg) & 46/46\\
\vspace{-0.25cm}\\
7. & $P_{pa,dia}$  & 1.285 & \b{1.062} & 1.065 & 1.290 & (mmHg) & 46/46\\
\vspace{-0.25cm}\\
8. & $P_{r,edp}$  & \b{0.540} & 0.599 & 0.786 & 2.787  & (mmHg) & 44/46\\
\vspace{-0.25cm}\\
9. & $P_{w}$   & 0.887 & \b{0.740} & 1.258 & 1.200 & (mmHg) & 46/46\\
\vspace{-0.25cm}\\
10. & $P_{cvp}$  & 1.194  & \b{0.371} & 0.733 & 0.821 & (mmHg) & 3/46\\
\vspace{-0.25cm}\\
11. &$V_{l,sys}$  & \texttt{None} & \texttt{None} & \texttt{None} & \texttt{None} & (mL) & 0/46\\
\vspace{-0.25cm}\\
12. &$V_{l,dia}$  & 50.960 & \b{47.391} & 103.402 & 119.440 & (mL) & 3/46\\
\vspace{-0.25cm}\\
13. &LVEF  & 0.012  & \b{0.008} & 0.015 & 0.024 & (-) & 45/46\\
\vspace{-0.25cm}\\
14. &CO  & 0.200  & \b{0.121} & 0.281 & 0.274 & (L/min) & 46/46\\
\vspace{-0.25cm}\\
15. &SVR   & 28.755 & \b{24.728} & 25.896 & 41.233 & (dyn$\cdot$s/cm$^5$) & 46/46\\
\vspace{-0.25cm}\\
16. &PVR & 1.924 & \b{1.812} & 2.010 & 3.795 & (dyn$\cdot$s/cm$^5$) & 46/46\\
\bottomrule
\end{tabular}
\end{center}}
\caption{Reconstruction error per output component ($N_w = 100$) for the EHR missing data inversion, using the criterion defined in equation~\eqref{equ: error-noise-level}, for noise levels: $\delta = 0.25, 0.5, 1.0, 2.0$. 
The minimum error across all 4 cases is highlighted using bold fonts. Only patients with more than 10 available measurements (i.e., 46 out of 84) are included. 
The total number of patients for which a given attribute was measured is also reported in the last column. Finally, the left ventricular systolic volume $V_{l,sys}$ is never measured on the selected patients.}
\label{table:ehr-error-global}
\end{table}

\begin{figure}[!ht]
\begin{subfigure}[b]{0.498\textwidth}
        \centering
         \includegraphics[scale=0.53]{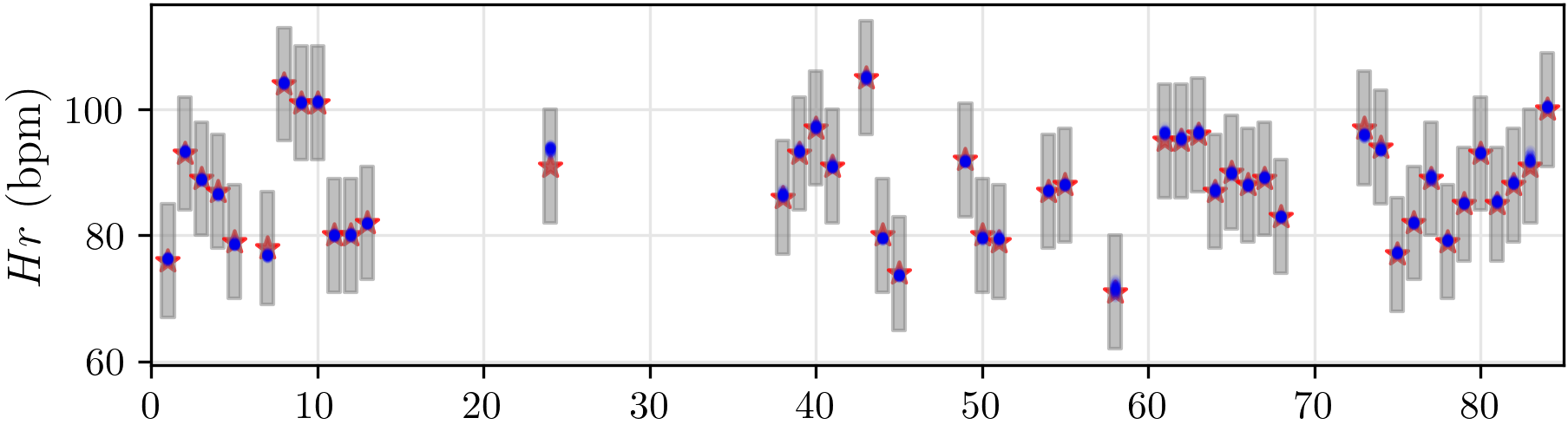}
    \end{subfigure}
\begin{subfigure}[b]{0.498\textwidth}
        \centering
         \includegraphics[scale=0.53]{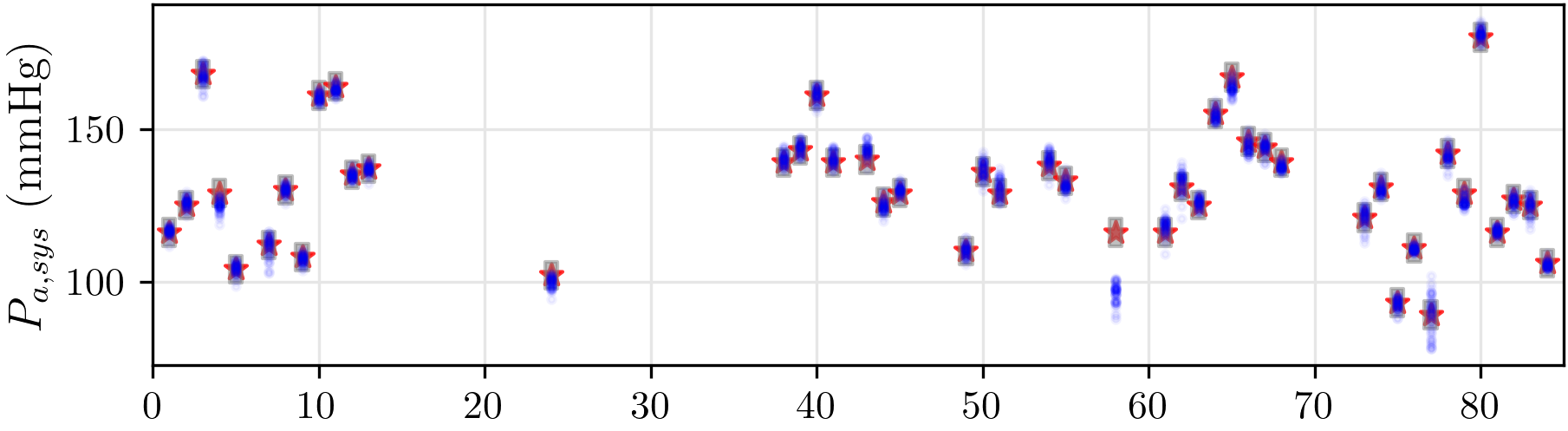}
    \end{subfigure}\\
    
\vspace{-0.2cm}
    \begin{subfigure}[b]{0.498\textwidth}
        \centering
         \includegraphics[scale=0.53]{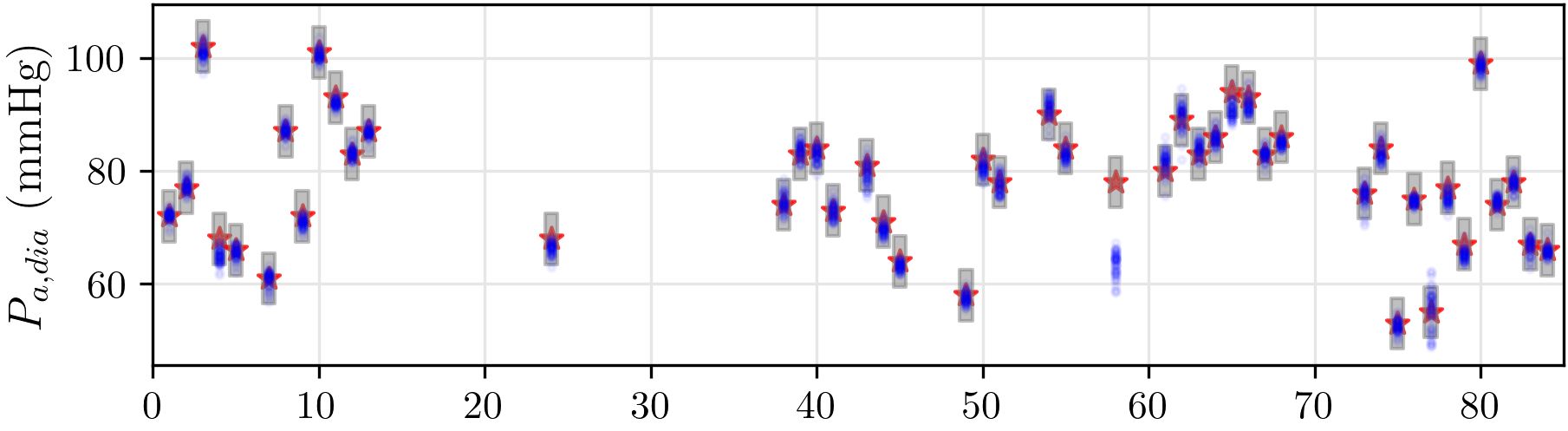}
    \end{subfigure}
\begin{subfigure}[b]{0.498\textwidth}
        \centering
         \includegraphics[scale=0.53]{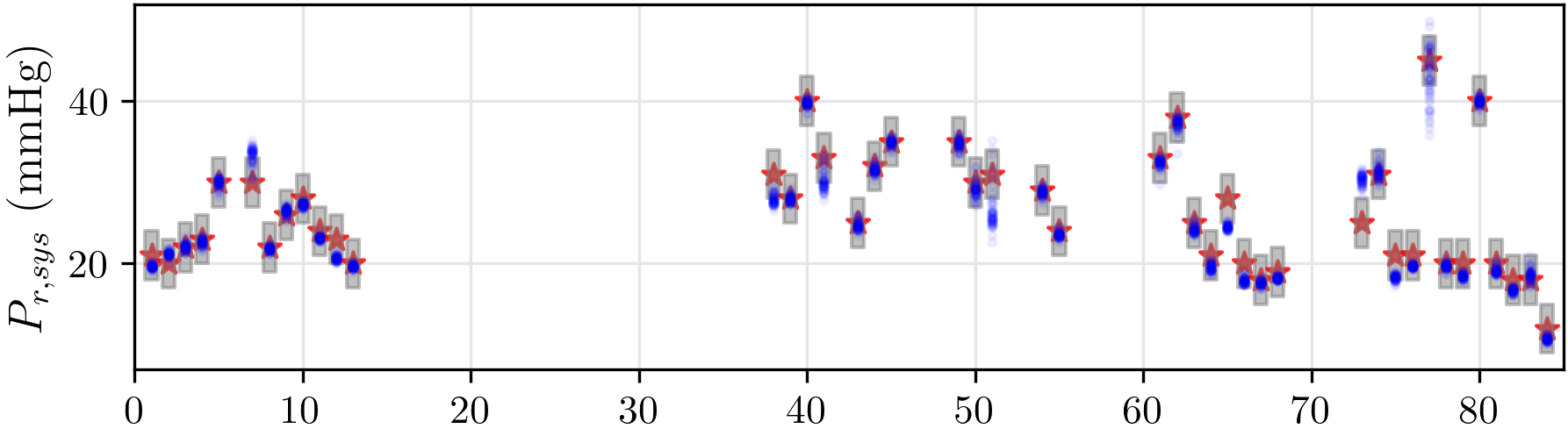}
    \end{subfigure}\\
    
    \vspace{-0.2cm}
    \begin{subfigure}[b]{0.498\textwidth}
        \centering
         \includegraphics[scale=0.53]{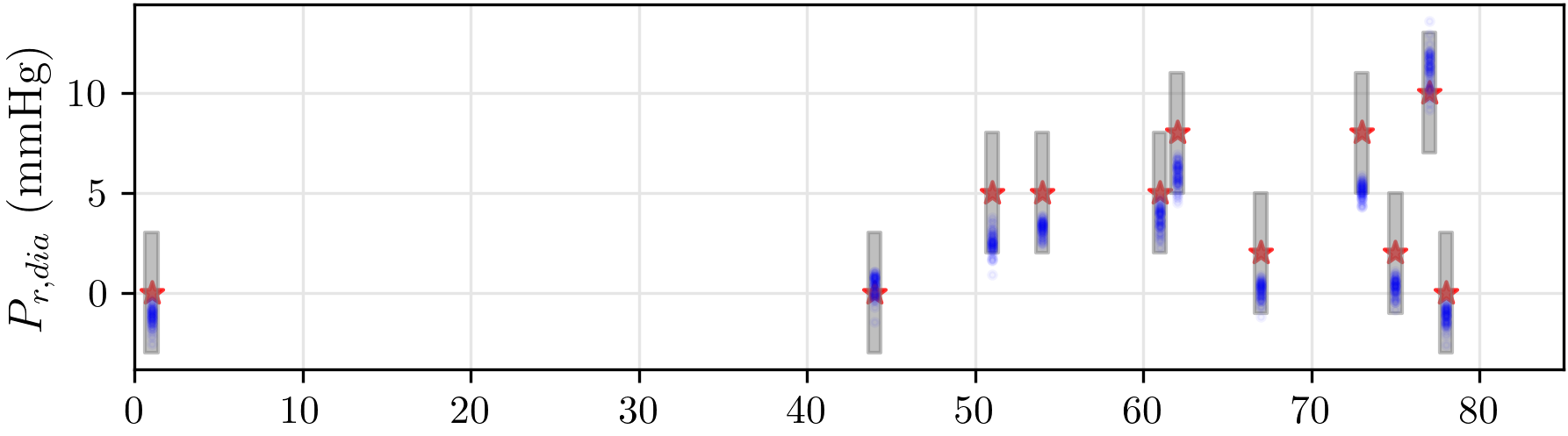}
    \end{subfigure}
\begin{subfigure}[b]{0.498\textwidth}
        \centering
         \includegraphics[scale=0.53]{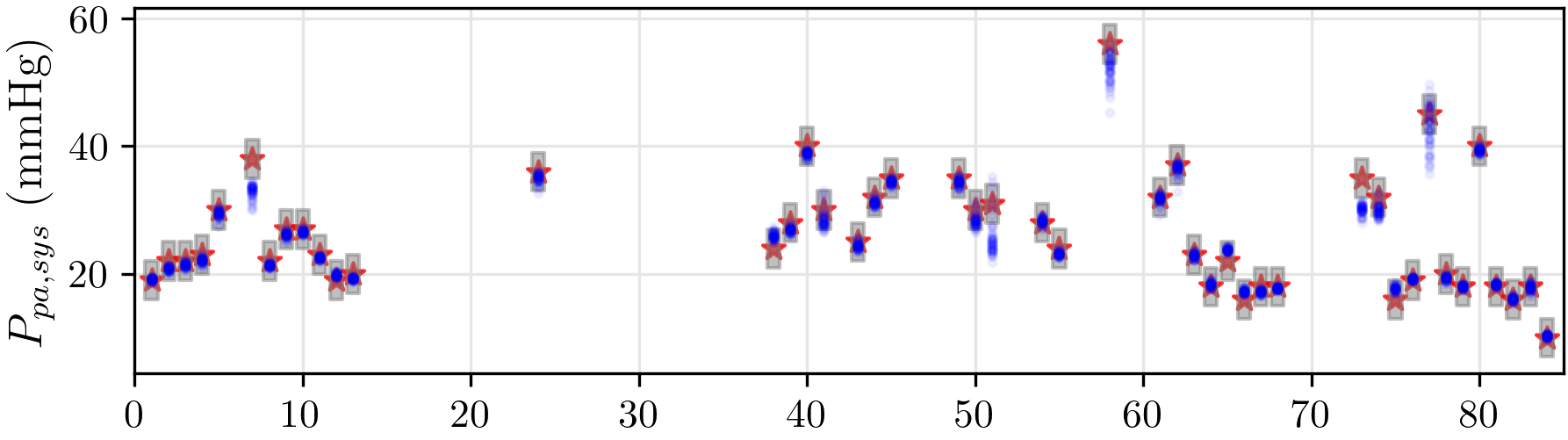}
    \end{subfigure}\\

    \vspace{-0.2cm}
    \begin{subfigure}[b]{0.498\textwidth}
        \centering
         \includegraphics[scale=0.53]{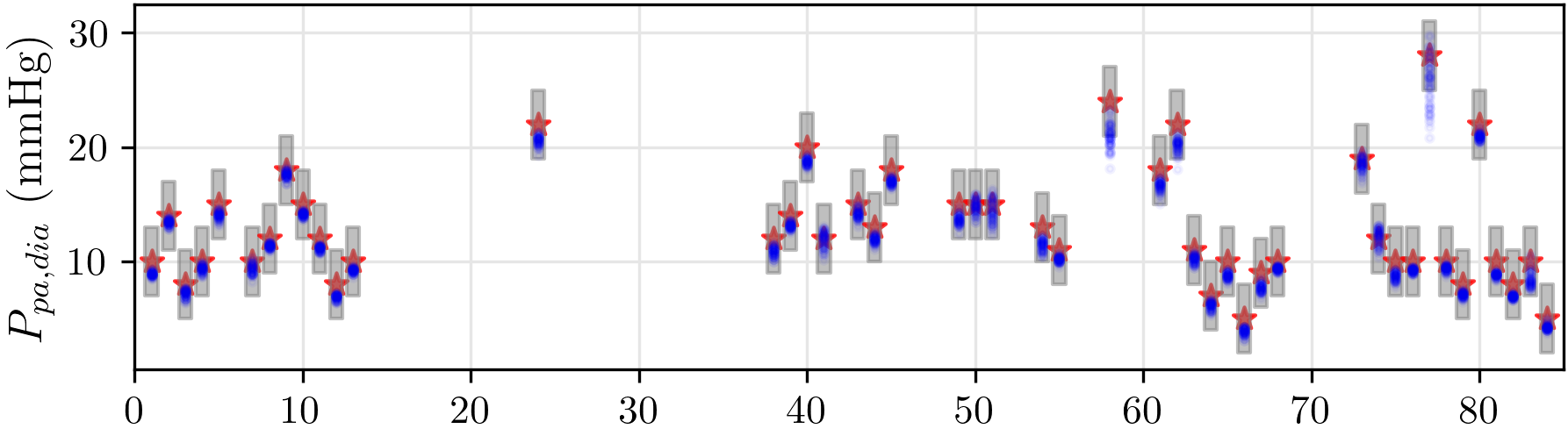}
    \end{subfigure}
\begin{subfigure}[b]{0.498\textwidth}
        \centering
         \includegraphics[scale=0.53]{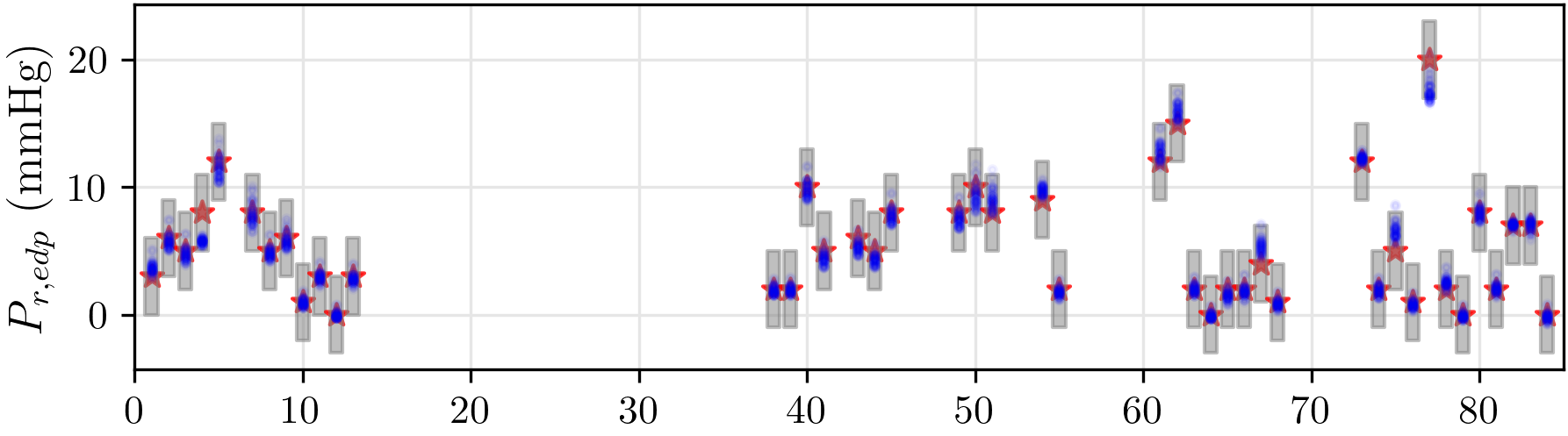}
    \end{subfigure}\\
    
    \vspace{-0.2cm}
    \begin{subfigure}[b]{0.498\textwidth}
        \centering
         \includegraphics[scale=0.53]{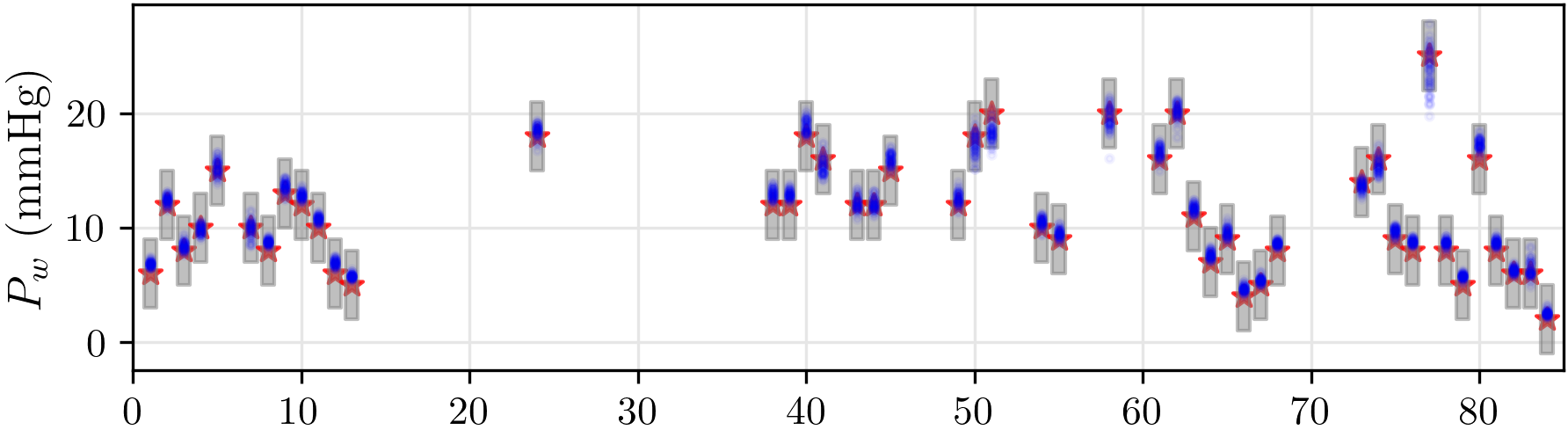}
    \end{subfigure}
\begin{subfigure}[b]{0.498\textwidth}
        \centering
         \includegraphics[scale=0.53]{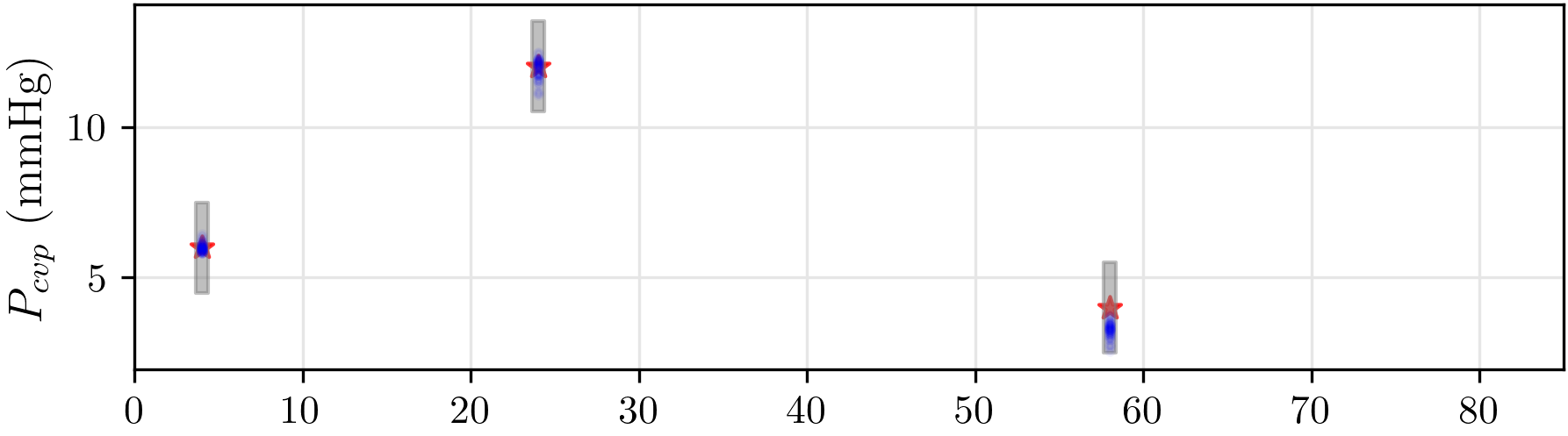}
    \end{subfigure}\\

    \vspace{-0.2cm}
    \begin{subfigure}[b]{0.498\textwidth}
        \centering
         \includegraphics[scale=0.53]{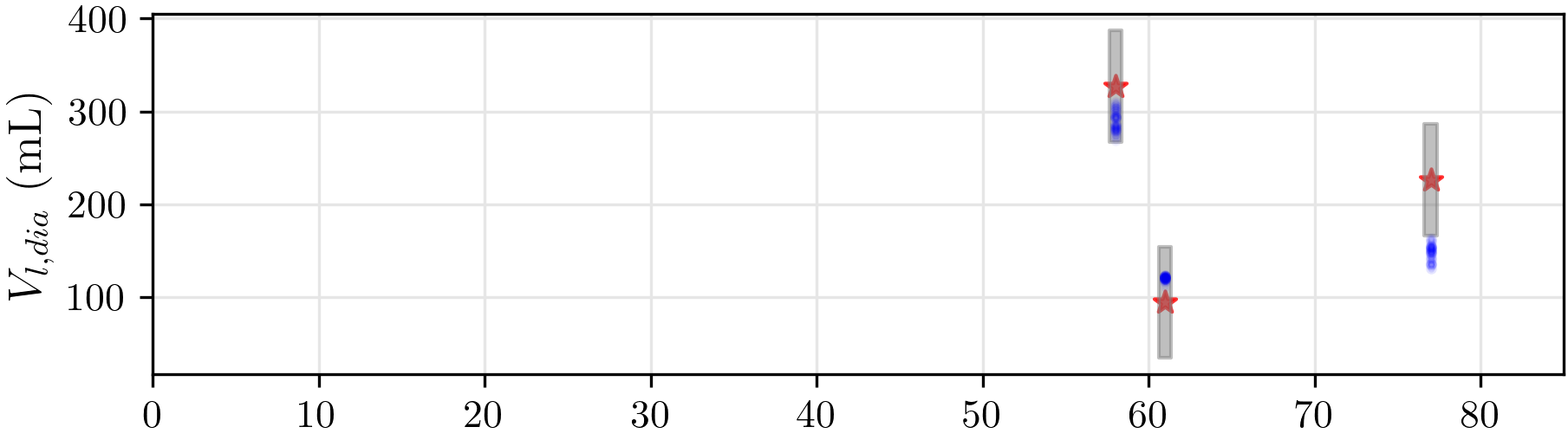}
    \end{subfigure}
\begin{subfigure}[b]{0.498\textwidth}
        \centering
         \includegraphics[scale=0.53]{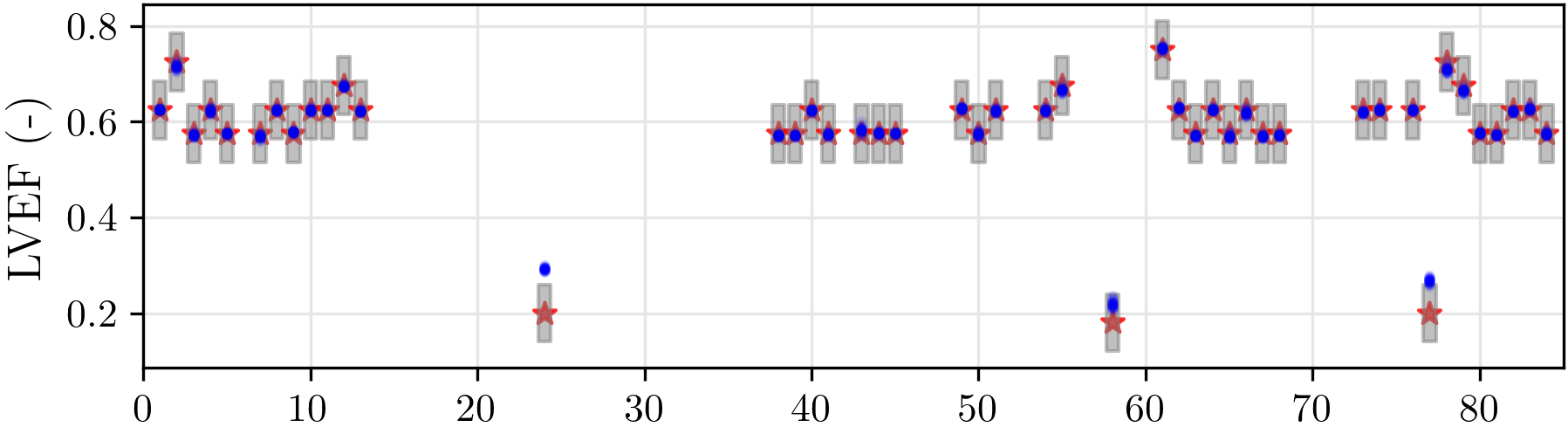}
    \end{subfigure}\\

    \vspace{-0.2cm}
    \begin{subfigure}[b]{0.498\textwidth}
        \centering
         \includegraphics[scale=0.53]{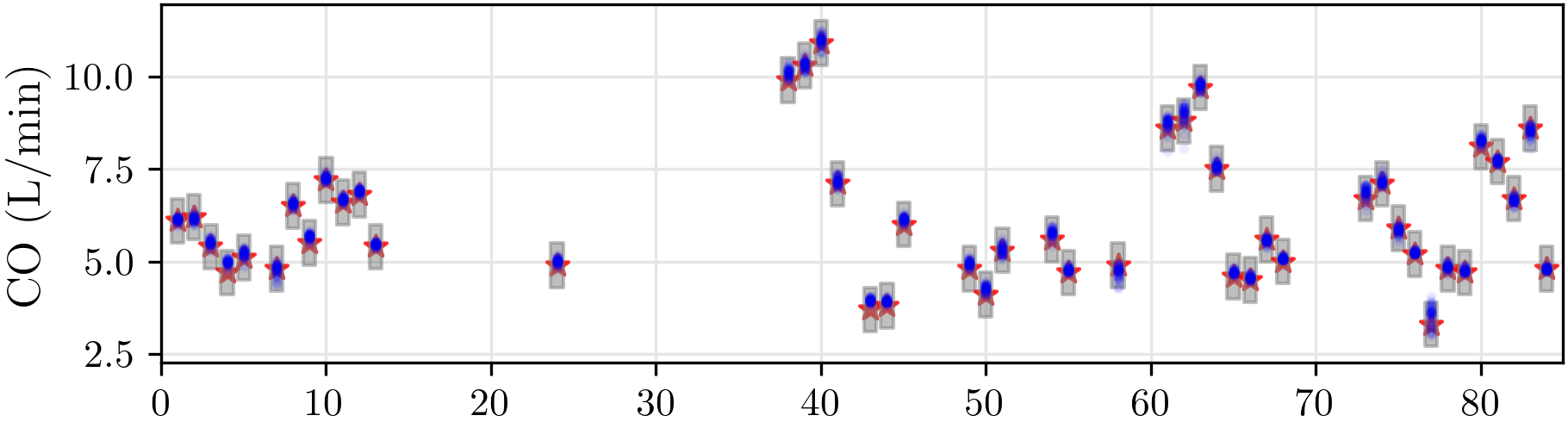}
    \end{subfigure}
\begin{subfigure}[b]{0.498\textwidth}
        \centering
         \includegraphics[scale=0.53]{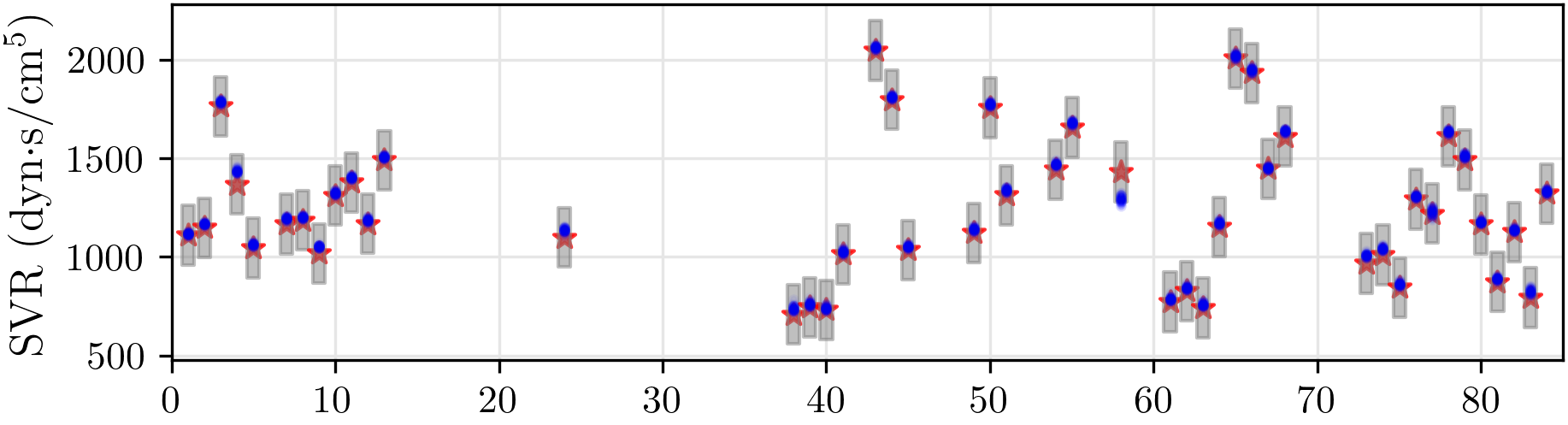}
    \end{subfigure}\\

    \vspace{-0.2cm}
    \begin{subfigure}[b]{0.498\textwidth}
        \centering
         \includegraphics[scale=0.53]{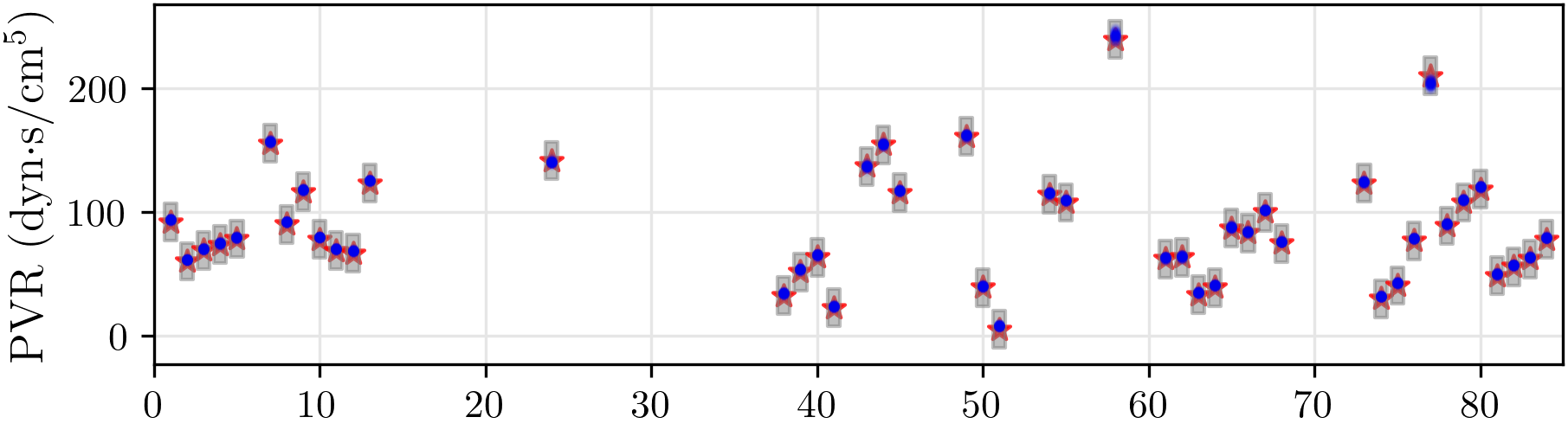}
    \end{subfigure}  
    \caption{Patient-specific predictions and uncertainty quantification for each CVSim-6 output component when the noise intensity factor is set as $\delta = 0.5$. 
    The red star is the EHR measurement of each patient and the blue dots are the network predictions, with each dot associated to a different latent variable realization. 
    The gray interval is obtained from the EHR measurement $\pm 3\sigma$, using the standard deviations reported in Table~\ref{table:cvsim6 output paras}. 
    The left ventricular systolic volume, $V_{l,sys}$, is missing for all patients with more than 10 measurements, hence the corresponding plot is omitted.}\label{tab:delta_results}
\end{figure}

The patient-specific prediction is validated with respect to each of the CVSim-6 output attribute and shown in Figure~\ref{tab:delta_results}, with respect to the noise level $\delta = 0.5$. A bounding interval of each EHR measurement $\pm 3\sigma$, using the reported uncertainty listed in Table~\ref{table:cvsim6 output paras}, is also superimposed to quantify the prediction accuracy.
From Figure~\ref{tab:delta_results}, we notice most of the input predictions can lead to an output close to the corresponding EHR measurement, although inaccuracy and large uncertainty exist in a few cases. 
Among all the output targets, the arterial systolic pressure, i.e., $P_{a,sys}$, is the most uncertain attribute in our experiment, despite when $\delta = 0.5$, the lowest prediction error is obtained.
%
%
%

\section{Discussion and conclusions}\label{sec:conclusion}

In this paper, we propose a new approach for model synthesis in stiff dynamical systems, and demonstrate its potential using both synthetic and real data with missing components from an EHR dataset. 
Once an inVAErt network is trained, it can be used as an \emph{enhanced digital twin} that can instantaneously (i.e., at the cost of a single network evaluation) perform a number of tasks including emulation, amortized inversion, output generation and missing data imputation.
In addition, a unique feature of inVAErt networks is their ability to solve ill-posed inverse problems with multiple solutions, by computing an entire manifold of parameters whose model outputs correspond to the same clinical observation. 
We believe that a paradigm where \emph{all} the solutions of an inverse problem are provided to the analyst (with an accuracy that only depends on the amount of training data provided to the system) is superior to any form of regularization. 

Moreover, besides using a latent space to handle structural non-identifiability, for the first time, noise from a closed-form model is added during training of an inVAErt network to alleviate practical non-identifiability in the solution of inverse problems.  
We focus on an application in computational physiology, using the CVSim-6 six-compartment lumped parameter model. We analyze this system in detail, showing that its stiffness is due to the presence of ideal unidirectional valves modeled through indicator functions, and ventricular-arterial coupling during late systole inducing very small RC constants. 

The results confirm the ability of inVAErt networks to leverage stiff dynamical systems in order to accurately perform a number of tasks involving not only synthetic data, but also real data with missing components.

There are many possible directions for future work, such as focusing on applications involving higher fidelity models such as one- or three-dimensional cardiovascular models, a broader family of latent space distributions, misspecified models, and beyond.

\section*{Acknowledgements}\label{sec:ack}

GGT and DES were supported by a NSF CAREER award \#1942662 (PI DES), a NSF CDS\&E award \#2104831 (University of Notre Dame PI DES) and used computational resources provided through the Center for Research Computing at the University of Notre Dame. CSL was partially supported by an Open Seed Fund between Pontificia Universidad Cat\'olica de Chile and the University of Notre Dame, by the grant Fondecyt 1211643, and by Centro Nacional de Inteligencia Artificial CENIA, FB210017, BASAL, ANID.

\bibliographystyle{abbrv}

\appendix

\section{Additional details on the CVSim-6 model formulation}\label{sec: cvsim6-details}

The following physiological quantities can be computed from the input parameters listed in Table~\ref{table:cvsim6 input paras}:

\begin{equation*}
    \begin{cases}
     \text{Total unstressed blood volume:} \ & V_{tot}^0 = V_{l}^0 + V_{a}^0 + V_{v}^0 + V_{r}^0 + V_{pa}^0 + V_{pv}^0 \ \text{(mL)} \ , \\
     \text{Duration of cardiac cycle:} \ & T_{tot} = 60 / Hr \ \text{(s)} \ , \\
     \text{Systolic time per cardiac cycle:} \ & T_{sys} = T_{tot} \cdot r_{sys} \ \text{(s)} \ , \\
     \text{Diastolic time per cardiac cycle:} \ & T_{dia} = T_{tot} - T_{sys} \ \text{(s)}\ , \\
    \end{cases}
\end{equation*}

and we assume a constant total blood volume $V_{tot} = 5000$ (mL), representative of an adult male weighing 70 kg~\cite{davis1991teaching}. To begin the in-silico simulation of the CVSim-6 system and preserve mass conservation at $t=0$, the initial pressures in the ODE system~\eqref{equ: cvsim6-ode} are determined by solving the following set of equations~\cite{davis1991teaching}

\begin{equation}
    \begin{cases} 
    C_{l,dia} (P_{l,dia}^0 - P_{th}) & - \ C_{l,sys} (P_{l,sys}^0 - P_{th})  \\
    \vspace{-0.35cm}\\ 
    &=  C_{r,dia} (P_{r,dia}^0 - P_{th}) - C_{r,sys} (P_{r,sys}^0 - P_{th})\ ,\\
    \vspace{-0.35cm}\\ 
    &= \displaystyle T_{sys} \frac{P^0_{l,sys} - P_a^0}{R_{l,out}} \ ,\\
    \vspace{-0.35cm}\\ 
    &= \displaystyle T_{tot} \frac{P^0_{a} - P_v^0}{R_{a}} \ ,\\
    \vspace{-0.35cm}\\ 
    &= \displaystyle T_{dia} \frac{P^0_{v} - P_{r,dia}^0}{R_{r, in}} \ ,\\
    \vspace{-0.35cm}\\ 
    &= \displaystyle T_{sys} \frac{P^0_{r,sys} - P_{pa}^0}{R_{r, out}} \ ,\\
     \vspace{-0.35cm}\\ 
    &= \displaystyle T_{tot} \frac{P^0_{pa} - P_{pv}^0}{R_{pv}} \ ,\\
    \vspace{-0.35cm}\\ 
    &= \displaystyle T_{dia} \frac{P^0_{pv} - P_{l,dia}^0}{R_{l,in}} \ ,\\
    \vspace{-0.35cm}\\ 
    V_{tot} - V_{tot}^0 &= C_{l,dia}(P^0_{l,dia} - P_{th}) + C_a(P_a^0 - \frac{1}{3}P_{th} ) + C_vP_v^0  \\
    & + \ C_{r,dia} (P_{r,dia}^0 - P_{th}) + C_{pa} (P_{pa}^0 - P_{th}) + C_{pv} (P_{pv}^0 - P_{th})  \ ,
    \end{cases}
    \label{equ: cvsim6-ic}
\end{equation}

where $P^0_{(\cdot)} = P_{(\cdot)}(t = 0)$ and the ventricular pressures are initialized at their diastolic values, namely, 
$$P_l(t=0) = P_{l,dia}^0, \  P_r(t=0) = P_{r,dia}^0 \ .$$

The time-varying capacitance of the two ventricles utilized in equations~\eqref{equ: cvsim6-ode} are defined through their reciprocals, i.e., through the elastance,
$$E_l(t) = 1/C_l(t), \ E_r(t) = 1/C_r(t) \ ,$$
via the following \emph{activation functions}~\cite{heldt2010cvsim, davis1991teaching} (these formulae are obtained from CVSim, an open-source code repository within the PhysioNet project~\cite{goldberger2000physiobank}).

\begin{equation}
    \begin{aligned}
        E_l(t) & = \begin{cases}
        \displaystyle \frac{1}{C_{l,dia}} ,& \ t = 0 \quad \textrm{or} \quad t> 1.5\cdot T_{sys}\ ,\\
        \vspace{-0.3cm}\\
         \displaystyle  \frac{1}{2}\big(\frac{1}{C_{l,sys}} - \frac{1}{C_{l,dia}}\big)\big(1-\cos(\pi\frac{t}{T_{sys}})\big) + \frac{1}{C_{l,dia}}  ,& \ 0 < t \leq T_{sys} \ ,\\
          \vspace{-0.3cm}\\
         \displaystyle  \frac{1}{2}\big(\frac{1}{C_{l,sys}} - \frac{1}{C_{l,dia}}\big)\big(1+\cos(2\pi\frac{t - T_{sys}}{T_{sys}})\big) + \frac{1}{C_{l,dia}}  ,& \ T_{sys} < t \leq 1.5\cdot T_{sys} \ . \\
        \end{cases} \\
        \vspace{-0.3cm} \\
        E_r(t) & = \begin{cases}
           \displaystyle  \frac{1}{C_{r,dia}} ,& \ t = 0 \quad \textrm{or} \quad t> 1.5\cdot T_{sys} \ ,\\
            \vspace{-0.3cm}\\
         \displaystyle  \frac{1}{2}\big(\frac{1}{C_{r,sys}} - \frac{1}{C_{r,dia}}\big)\big(1-\cos(\pi\frac{t}{T_{sys}})\big) + \frac{1}{C_{r,dia}}  ,& \ 0 < t \leq T_{sys}\ , \\
          \vspace{-0.3cm}\\
         \displaystyle  \frac{1}{2}\big(\frac{1}{C_{r,sys}} - \frac{1}{C_{r,dia}}\big)\big(1+\cos(2\pi\frac{t - T_{sys}}{T_{sys}})\big) + \frac{1}{C_{r,dia}}  ,& \ T_{sys} < t \leq 1.5\cdot T_{sys} \ .\\
        \end{cases} 
        \end{aligned}
        \label{equ: driversfunction}
\end{equation}

At each step, the stressed volume at each compartment is calculated by the linear pressure-volume relationship~\cite{davis1991teaching}

\begin{equation}
\begin{cases}
    V_{l}(t) &= V_{l}^{0} + \big(P_l(t) - P_{th} \big) C_l(t)  \ ,  \\
    \vspace{-0.35cm}\\
    V_{a}(t) &= V_{a}^{0} + \big(P_a(t) - \frac{1}{3}P_{th}\big) C_a  \ ,  \\
    \vspace{-0.35cm}\\
    V_{v}(t) &= V_{v}^{0} + P_v(t) C_v \ , \\
    \vspace{-0.35cm}\\
    V_{r}(t) &= V_{r}^{0} + \big(P_r(t) - P_{th} \big) C_r(t) \ ,  \\
    \vspace{-0.35cm}\\
    V_{pa}(t) &= V_{pa}^{0} + \big(P_{pa}(t) - P_{th} \big) C_{pa} \ ,  \\
    \vspace{-0.35cm}\\
    V_{pv}(t) &= V_{pv}^{0} + \big(P_{pv}(t) - P_{th}\big) C_{pv} \ ,
\end{cases}
\end{equation}

and their sum: $\sum V(t) = V_{l}(t) + V_{a}(t) + V_{v}(t)+V_{r}(t)+V_{pa}(t)+V_{pv}(t)$ at any given time $t$ should be consistent with the previously specified constant $V_{tot} = 5000$ (mL).
This constraint provides a practical way to validate the numerical solver for the CVSim-6 system, as demonstrated in Figure~\ref{fig:cvsim6-RK4-explicit-implicit}.

\section{Details of neural network modeling}\label{sec: NN-details}

We build and train inVAErt networks using the \texttt{PyTorch} open-source machine learning platform~\cite{paszke2019pytorch}, together with an Adam optimizer~\cite{kingma2014adam} and a StepLR learning rate scheduler.
The fully-connected MLP networks are the building-blocks of the neural emulator $NN_e$, variational encoder $NN_v$, and the decoder $NN_d$. As for the density estimator $NN_f$, we implement the Real-NVP (real-valued non-volume preserving transformations) normalizing flow architecture~\cite{dinh2016density}.
For the network activations, except for the density estimator, where the recommended functions in~\cite{dinh2016density} are adopted, we use the smoothed ReLU function \emph{Swish}~\cite{ramachandran2017searching} due to its superior performance over other choices.

\subsection{Study of structural non-identifiability}\label{sec: NN-details-cvsim6-str}

In this section, we summarize details of the network model utilized in Section~\ref{sec: cvsim6-str} (see~\Cref{table:cvsim6-no-noise-hyper,table:CVSim6-no-noise training}).

\begin{table}[ht!]
{\small
\begin{center}
\begin{tabular}{@{} l c c c c @{}}
\toprule
Module & Hidden units & Num. of layers & Activation & Total parameters\\
\midrule
$NN_e$ & 60 & 8 & \texttt{SiLU} & 24376\\
$NN_v$ & 32 & 6 & \texttt{SiLU} & 6246\\
$NN_d$ & 64 & 6 & \texttt{SiLU} & 20439\\
\bottomrule
\end{tabular}

\bigskip
\begin{tabular}{@{} l c c c c c c @{}}
\toprule
Module & Hidden units & Layers per block & Blocks & Activation & BatchNorm & Total parameters\\
\midrule
$NN_f$ & 24 & 4 & 16 & \texttt{Tanh, ReLU} & \texttt{False} & 51712\\
\bottomrule
\end{tabular}
\end{center}
\caption{Summary of InVAErt modules for the CVSim-6 system with synthetic data.}
\label{table:cvsim6-no-noise-hyper}}
\end{table}

\begin{table}[ht!]
{\small
\begin{center}
\begin{tabular}{@{} l l l l l l @{}}
\toprule
Module & Minibatch size & Initial lr & Total Epochs & Epochs per decay & Decay rate\\
\midrule
$NN_e$ & 256 & $1\cdot 10^{-3}$ & 20000 & 100 & 0.98\\
$NN_f$ & 512 & $2\cdot 10^{-3}$ & 1500 & 200 & 0.85\\
$NN_v+NN_d$ & 256 & $1\cdot 10^{-3}$ & 20000& 100& 0.985\\
\bottomrule
\end{tabular}
\end{center}
\caption{Hyperparameter choices for the network training: CVSim-6 system with synthetic data. The penalties for the loss functions are $[\beta_d, \beta_v, \beta_r] = [1, 2000, 20]$ and the latent space dimensionality is $\dim(\bw) = 19$. }
\label{table:CVSim6-no-noise training}}
\end{table}


\subsection{Study of the EHR dataset}\label{sec: NN-details-cvsim6-EHR}

In this section, we summarize details of the network model and hyperparameters utilized in the study of the EHR dataset (see Table~\ref{table:cvsim6-ehr-model} and Table~\ref{table:CVSim6-ehr-training}).

\begin{table}[ht!]
{\small
\begin{center}
\begin{tabular}{@{} l c c c c @{}}
\toprule
Module & Hidden units & Num. of layers & Activation & Total parameters\\
\midrule
$NN_e$ & 80 & 8 & \texttt{SiLU} & 42096\\
$NN_v$ & 32 & 6 & \texttt{SiLU} & 6246\\
$NN_d$ & 64 & 6 & \texttt{SiLU} & 20439\\
\bottomrule
\end{tabular}

\bigskip
\begin{tabular}{@{} l c c c c c c @{}}
\toprule
Module & Hidden units & Layers per block & Blocks & Activation & BatchNorm & Total parameters\\
\midrule
$NN_f$ & 18 & 4 & 10 & \texttt{Tanh, ReLU} & \texttt{True} & 20280\\
\bottomrule
\end{tabular}
\end{center}
\caption{Summary of InVAErt modules for the study of the EHR dataset.}
\label{table:cvsim6-ehr-model}}
\end{table}

\begin{table}[ht!]
{\small
\begin{center}
\begin{tabular}{@{} l l l l l l @{}}
\toprule
Module & Minibatch size & Initial lr & Total Epochs & Epochs per decay & Decay rate\\
\midrule
$NN_e$ & 256 & $3\cdot 10^{-3}$ & 25000 & 100 & 0.98\\
$NN_f$ & 512 & $4\cdot 10^{-3}$ & 1000 & 200 & 0.5\\
$NN_v+NN_d$ & 512 & $1\cdot 10^{-3}$ & 2000& 100& 0.97\\
\bottomrule
\end{tabular}
\end{center}
\caption{Hyperparameter choices for the network training: CVSim-6 system with the EHR dataset. In the training of density estimator $NN_f$ and the inverse model $NN_v + NN_d$, $l_2$ weight decay is implemented with a penalty of $2\times 10^{-4}$. The penalties for the loss functions are $[\beta_d, \beta_v, \beta_r] = [1, 1000, 10]$ and the latent space dimensionality is $\dim(\bw) = 19$. }
\label{table:CVSim6-ehr-training}}
\end{table}


\section{Additional results}\label{sec: additional results}
\subsection{Forward emulation and density estimation}\label{sec: cvsim6-emu-nf}

In this section, we evaluate the performance of InVAErt components $NN_e$ and $NN_f$ for the CVSim-6 system studied in Section~\ref{sec: cvsim6-str}.
First, we test the trained neural emulator $NN_e$ via the validation dataset, and summarize the absolute prediction error for each output component in Table~\ref{table:emulator-acc-component}, where a noticeable overall accuracy can be found.

\begin{table}[ht!]
{\small
\begin{center}
\begin{tabular}{@{} l c c c c c @{}}
\toprule
& $Hr$ (bpm) & $P_{a,sys}$ (mmHg) & $P_{a,dia}$ (mmHg) & $P_{r,sys}$ (mmHg) & $P_{r,dia}$ (mmHg) \\
\midrule
Average & 9.81$\cdot 10^{-3}$ & 4.72$\cdot 10^{-2}$ & 4.52$\cdot 10^{-2}$ & 1.20$\cdot 10^{-2}$ & 3.44$\cdot 10^{-3}$ \\
Max     & 7.57$\cdot 10^{-2}$ & 6.68$\cdot 10^{-1}$ & 3.92$\cdot 10^{-1}$ & 2.35$\cdot 10^{-1}$ & 2.33$\cdot 10^{-2}$ \\
Std     & 8.59$\cdot 10^{-3}$ & 4.65$\cdot 10^{-2}$ & 4.02$\cdot 10^{-2}$ & 1.26$\cdot 10^{-2}$ & 2.81$\cdot 10^{-3}$ \\
\bottomrule
\end{tabular}
\bigskip

\begin{tabular}{@{} l c c c c c @{}}
\toprule
 & $P_{pa,sys}$ (mmHg) & $P_{pa,dia}$ (mmHg) & $P_{r,edp}$ (mmHg) & $P_{w}$ (mmHg) & $P_{cvp}$ (mmHg) \\
\midrule
Average & 1.10$\cdot 10^{-2}$ & 9.77$\cdot 10^{-3}$ & 2.74$\cdot 10^{-3}$ & 7.59$\cdot 10^{-3}$ & 3.68$\cdot 10^{-3}$ \\
Max     & 2.87$\cdot 10^{-1}$ & 1.86$\cdot 10^{-1}$ & 7.82$\cdot 10^{-2}$ & 7.38$\cdot 10^{-2}$ & 4.24$\cdot 10^{-2}$ \\
Std     & 1.18$\cdot 10^{-2}$ & 9.94$\cdot 10^{-3}$ & 3.34$\cdot 10^{-3}$ & 7.30$\cdot 10^{-3}$ & 3.49$\cdot 10^{-3}$ \\
\bottomrule
\end{tabular}
\bigskip

\begin{tabular}{@{} l c c c c c c @{}}
\toprule
 & $V_{l,sys}$ (mL) & $V_{l,dia}$ (mL) & LVEF (-) & CO (L/min) & SVR (dyn$\cdot$s/cm$^5$) & PVR (dyn$\cdot$s/cm$^5$) \\
\midrule
Average & 2.33$\cdot 10^{-2}$ & 4.45$\cdot 10^{-2}$ & 1.08$\cdot 10^{-4}$ & 1.91$\cdot 10^{-3}$ & 1.81$\cdot 10^{-1}$ & 1.95$\cdot 10^{-2}$ \\
Max     & 2.88$\cdot 10^{-1}$ & 4.51$\cdot 10^{-1}$ & 2.04$\cdot 10^{-3}$ & 2.00$\cdot 10^{-2}$ & 1.62$\cdot 10^{+0}$ & 2.12$\cdot 10^{-1}$ \\
Std     & 2.38$\cdot 10^{-2}$ & 4.33$\cdot 10^{-2}$ & 1.01$\cdot 10^{-4}$ & 1.82$\cdot 10^{-3}$ & 1.60$\cdot 10^{-1}$ & 1.84$\cdot 10^{-2}$ \\
\bottomrule
\end{tabular}
\bigskip
\caption{Statistics of the emulator prediction error across the validation dataset of size 4,000.}
\label{table:emulator-acc-component}
\end{center}}
\end{table}

\begin{figure}[H]
\begin{subfigure}[b]{0.121\textwidth}
        \centering
         \includegraphics[scale=0.154]{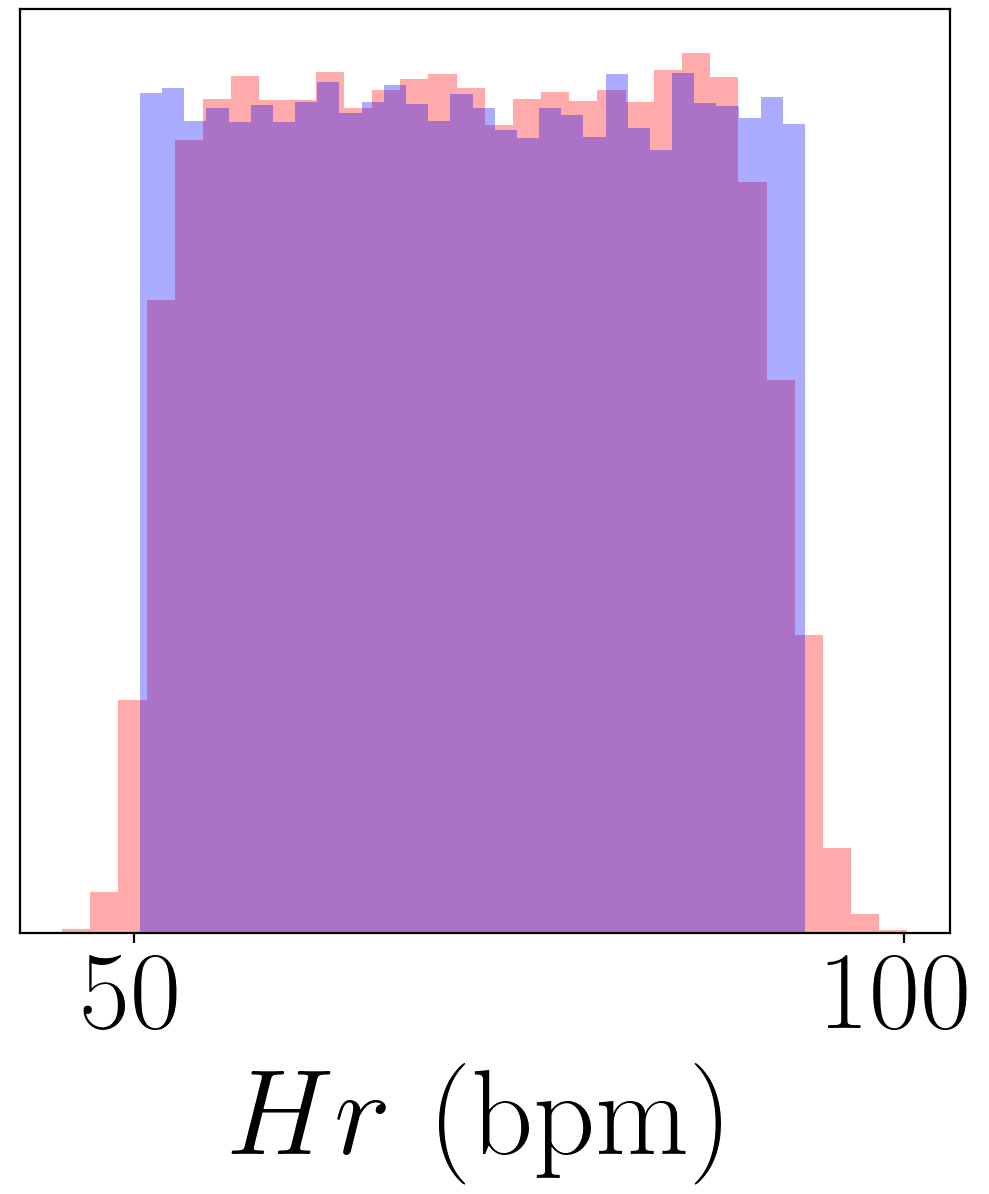}
\end{subfigure}
\begin{subfigure}[b]{0.121\textwidth}
        \centering
         \includegraphics[scale=0.154]{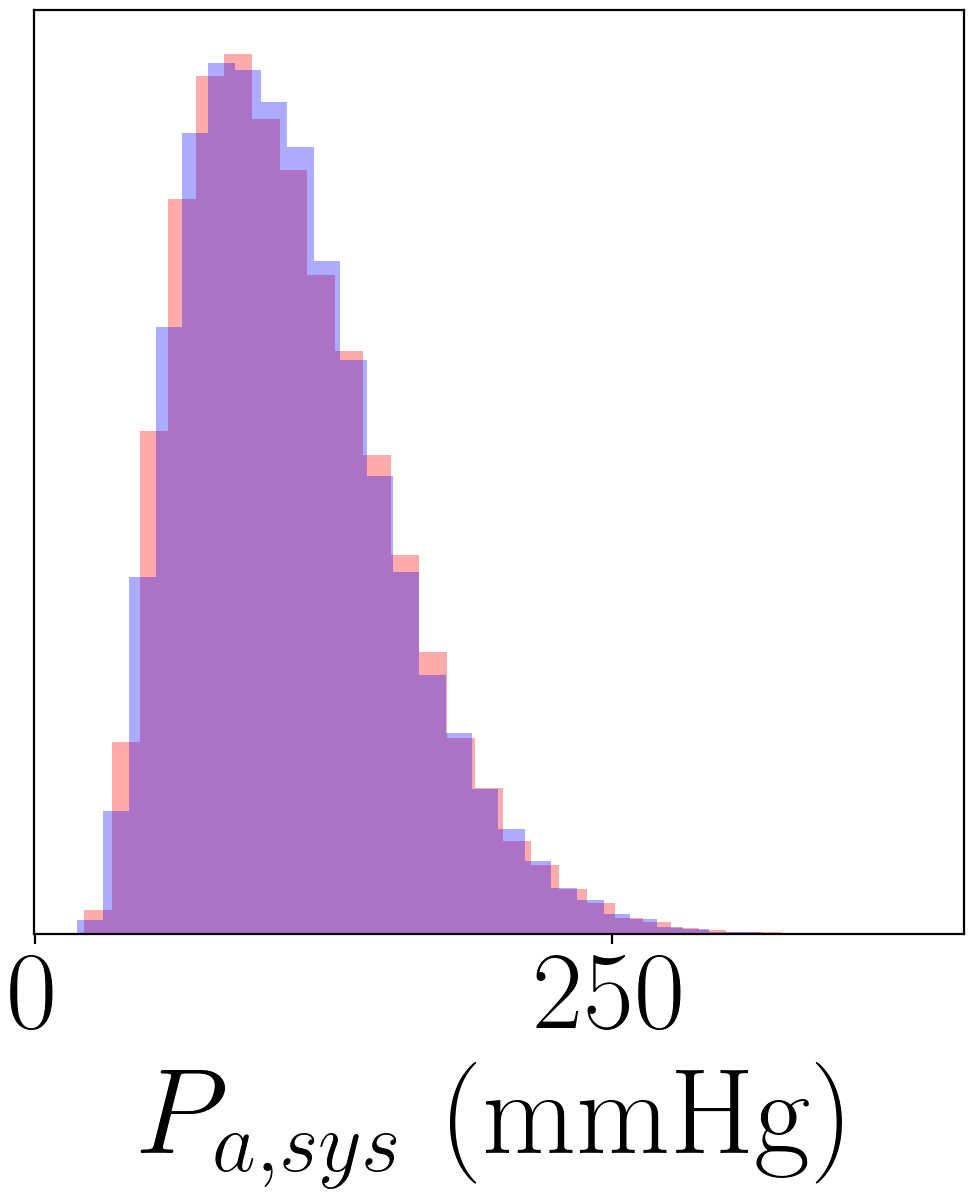}
\end{subfigure}
\begin{subfigure}[b]{0.121\textwidth}
        \centering
         \includegraphics[scale=0.154]{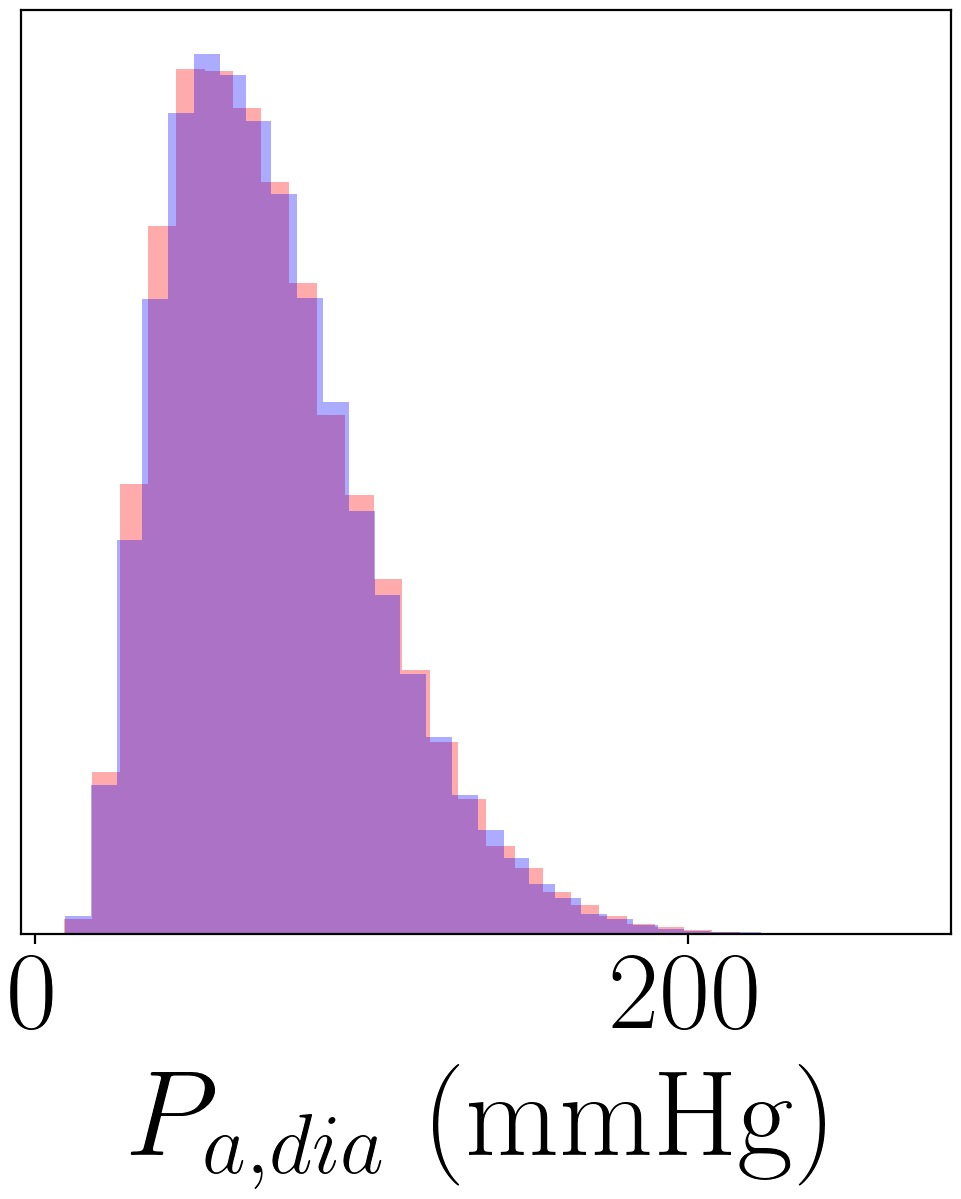}
\end{subfigure}
\begin{subfigure}[b]{0.121\textwidth}
        \centering
         \includegraphics[scale=0.154]{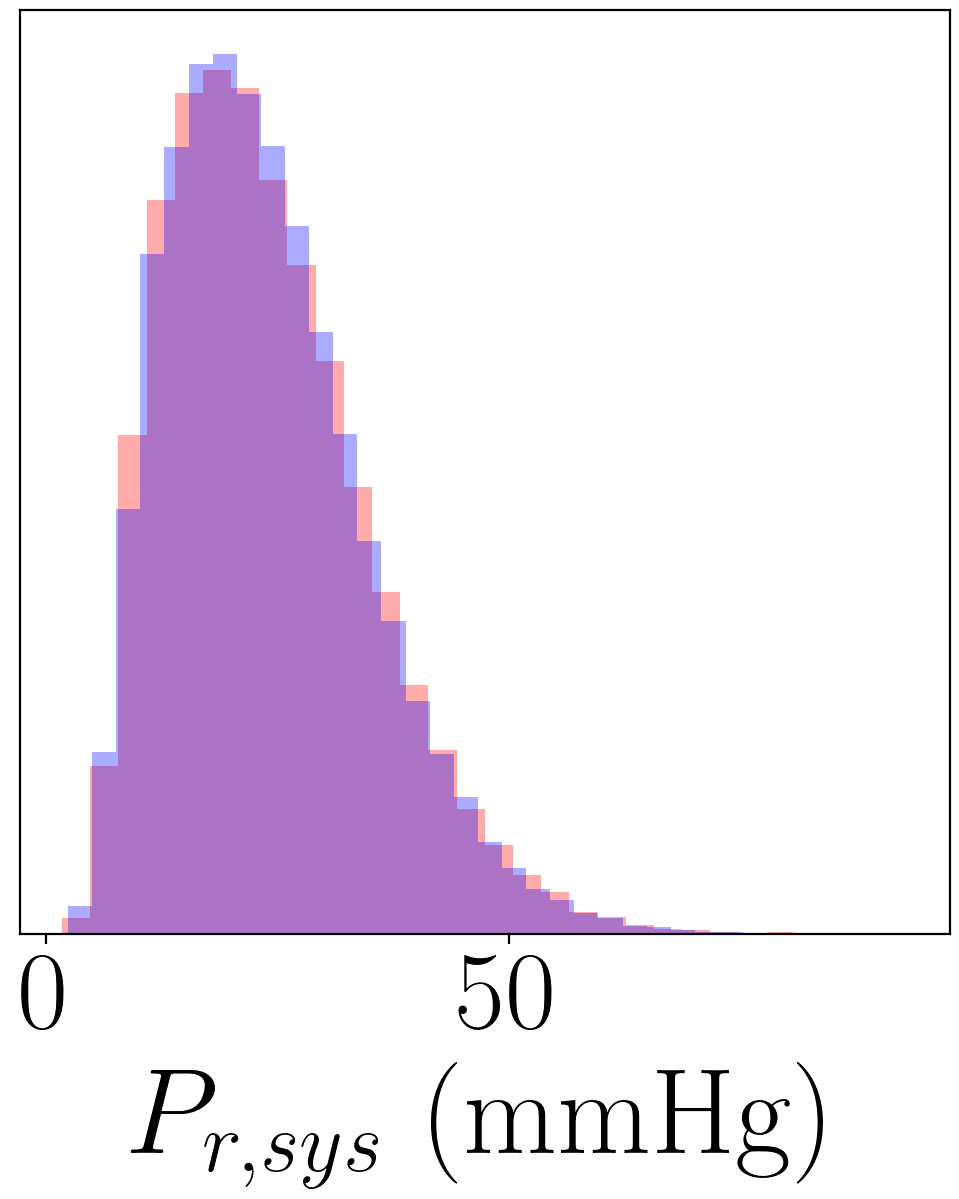}
\end{subfigure}
\begin{subfigure}[b]{0.121\textwidth}
        \centering
         \includegraphics[scale=0.154]{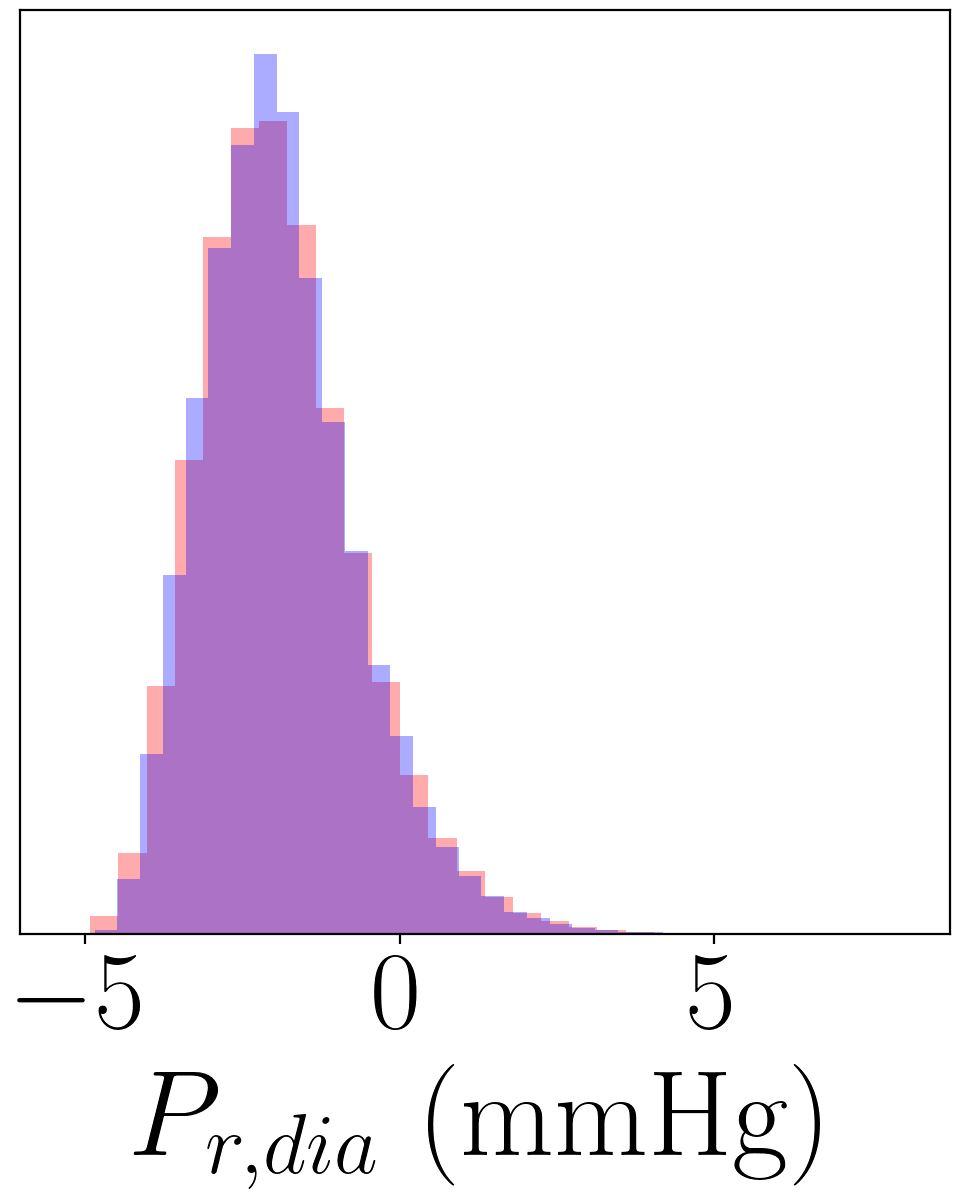}
\end{subfigure}
\begin{subfigure}[b]{0.121\textwidth}
        \centering
         \includegraphics[scale=0.154]{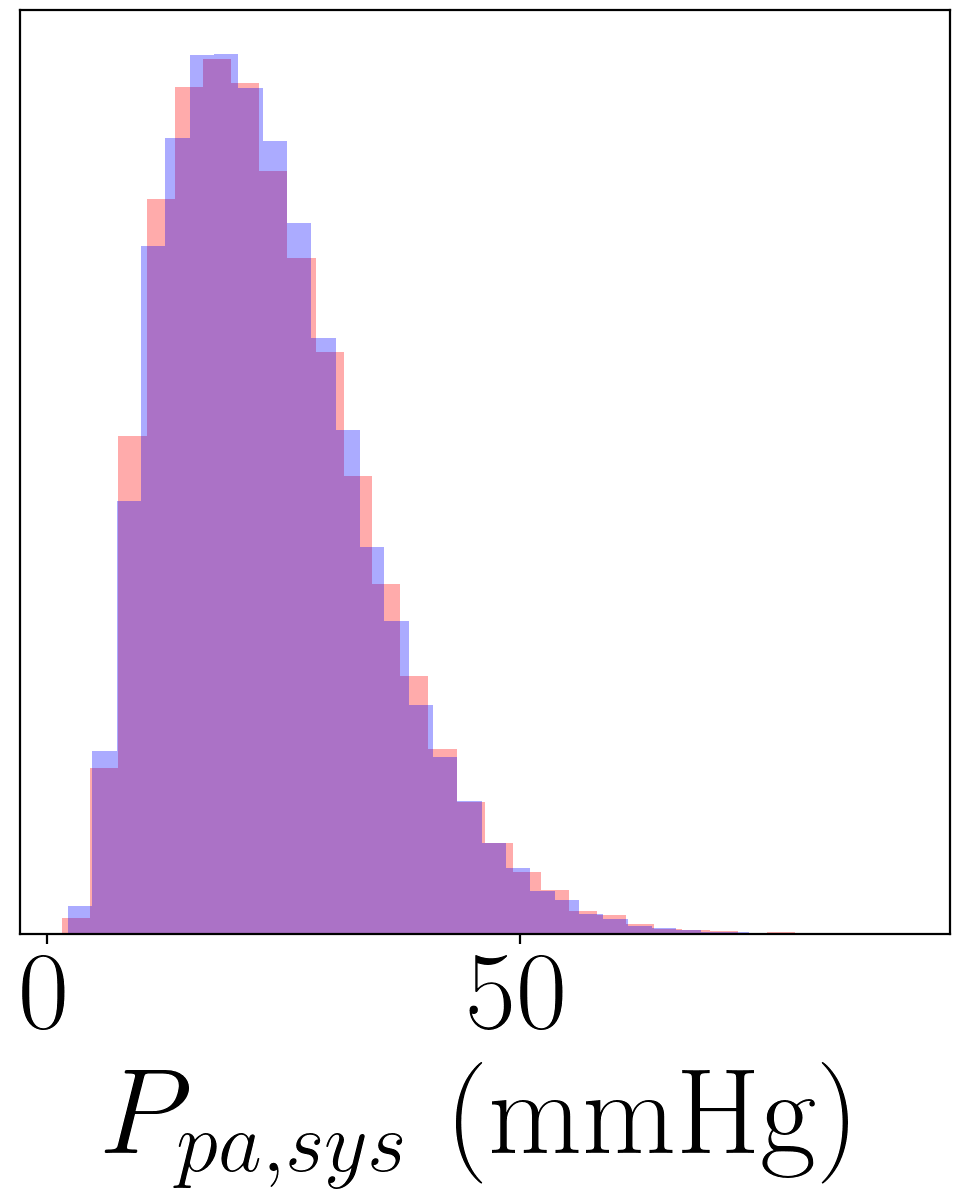}
\end{subfigure}
\begin{subfigure}[b]{0.121\textwidth}
        \centering
         \includegraphics[scale=0.154]{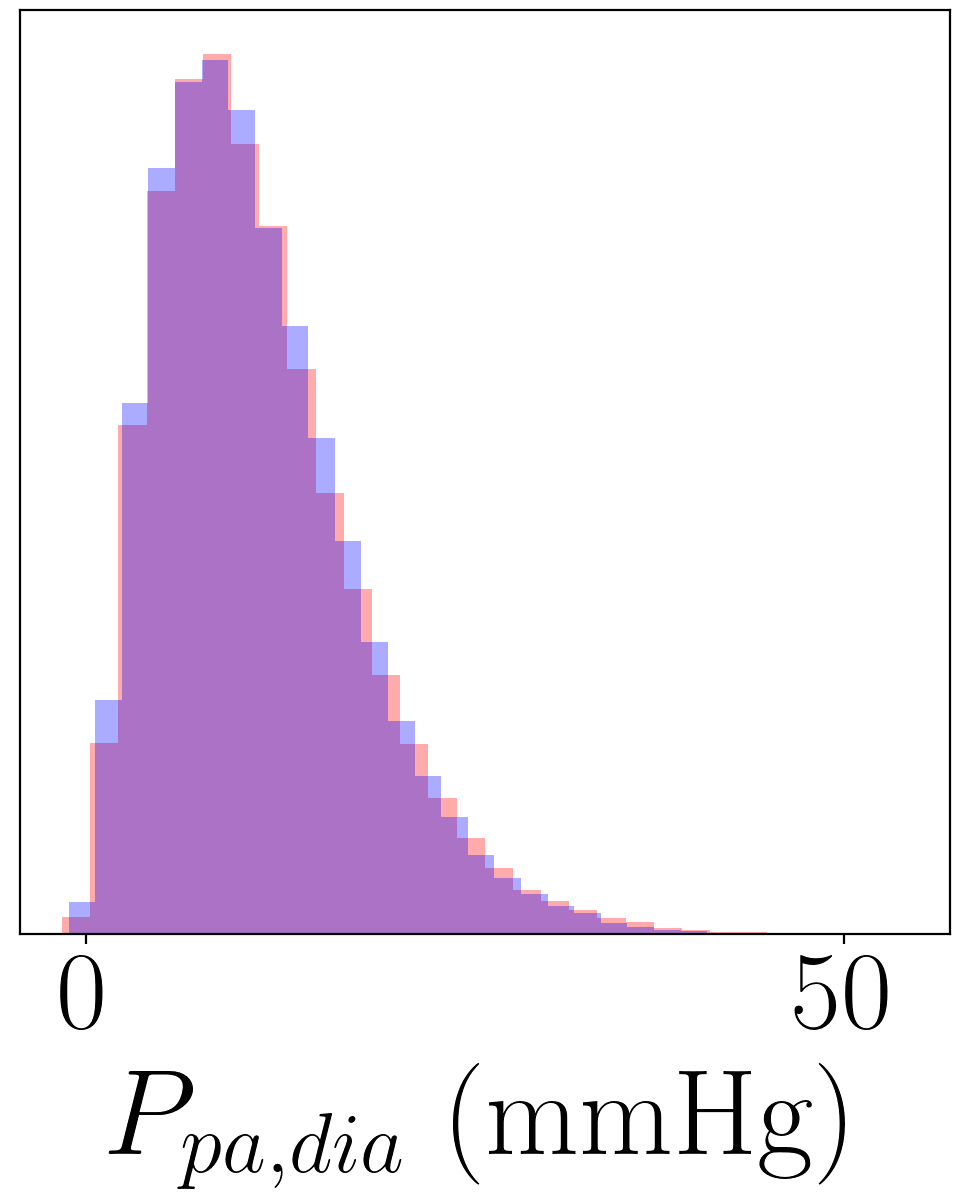}
\end{subfigure}
\begin{subfigure}[b]{0.121\textwidth}
        \centering
         \includegraphics[scale=0.154]{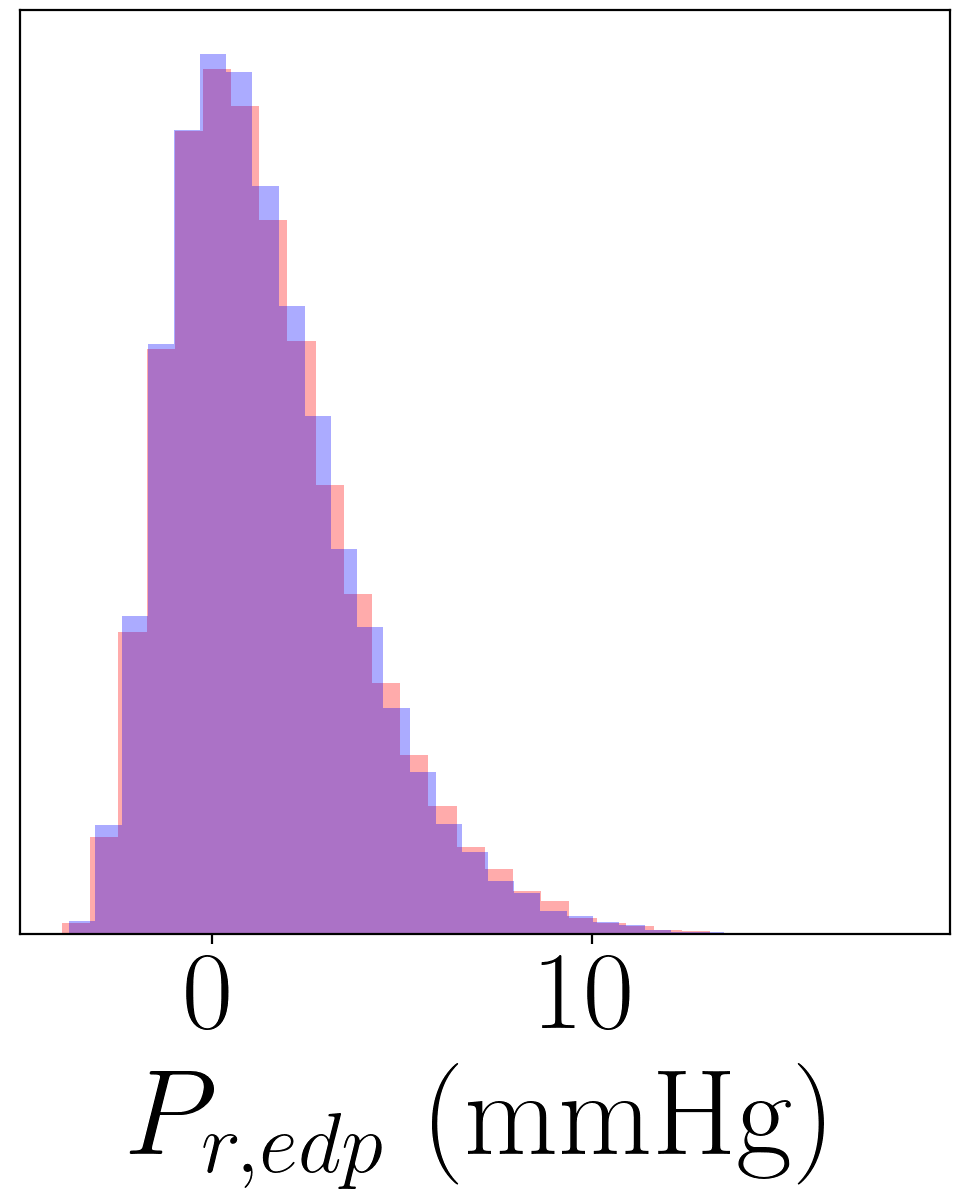}
\end{subfigure}\\

\begin{subfigure}[b]{0.121\textwidth}
        \centering
         \includegraphics[scale=0.154]{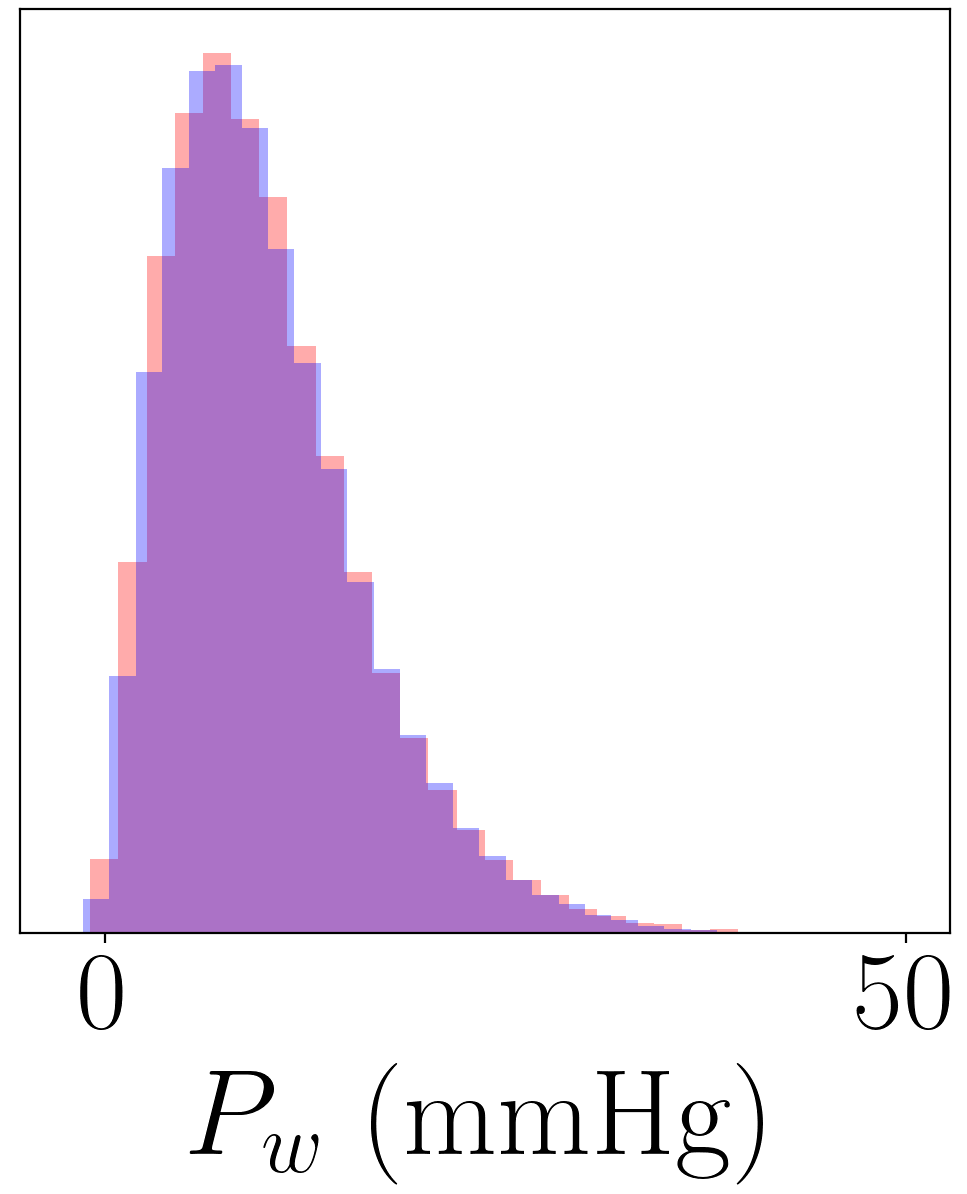}
\end{subfigure}
\begin{subfigure}[b]{0.121\textwidth}
        \centering
         \includegraphics[scale=0.154]{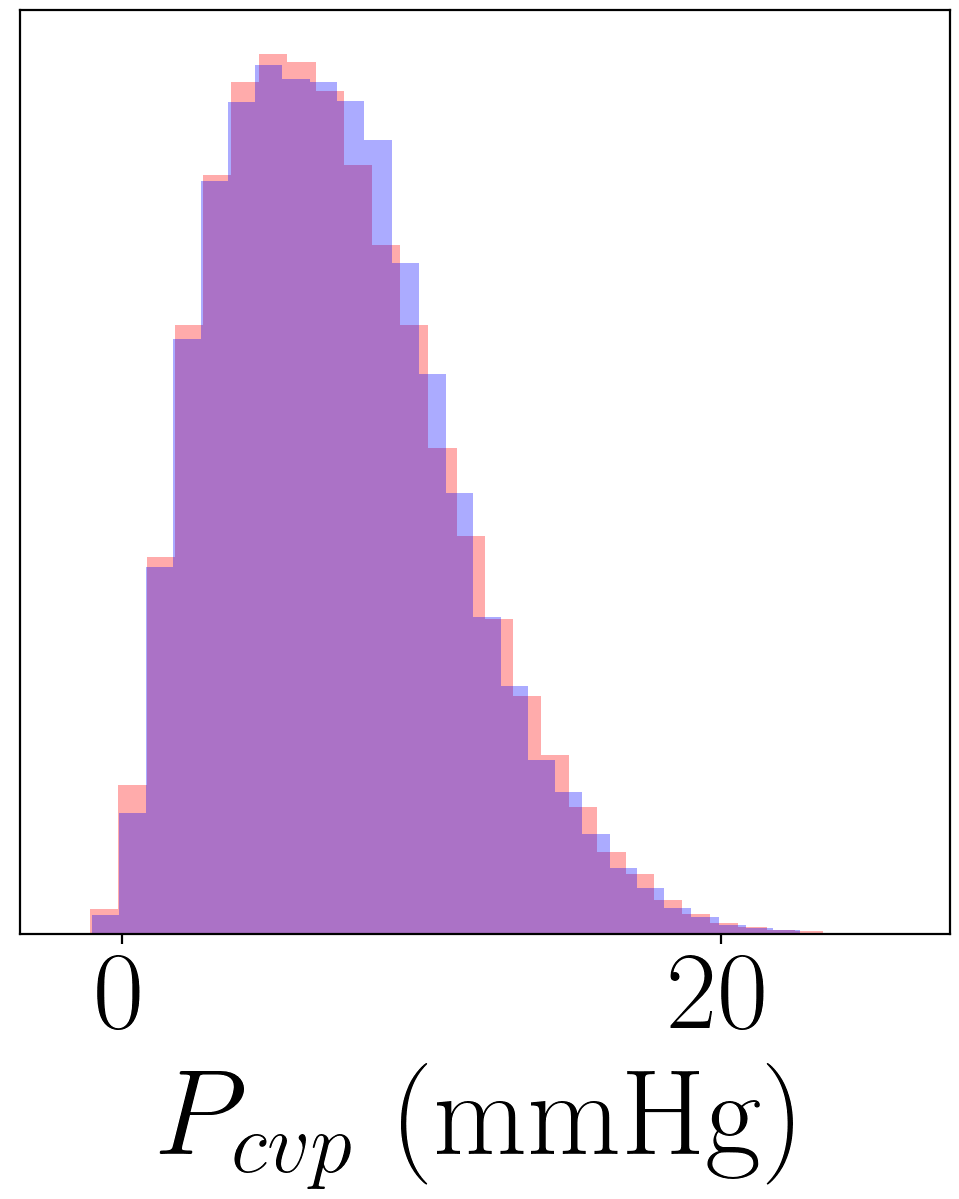}
\end{subfigure}
\begin{subfigure}[b]{0.121\textwidth}
        \centering
         \includegraphics[scale=0.154]{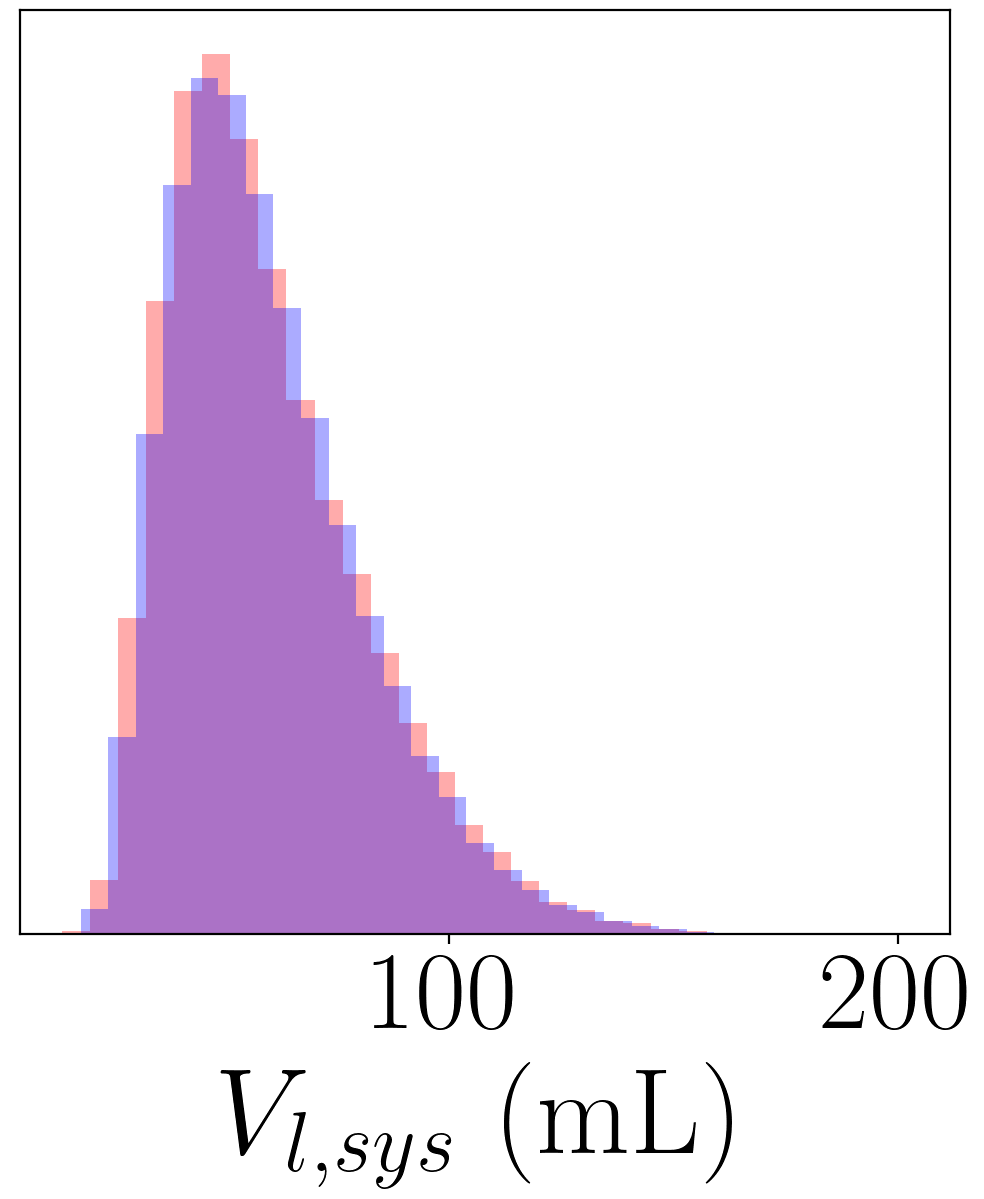}
\end{subfigure}
\begin{subfigure}[b]{0.121\textwidth}
        \centering
         \includegraphics[scale=0.154]{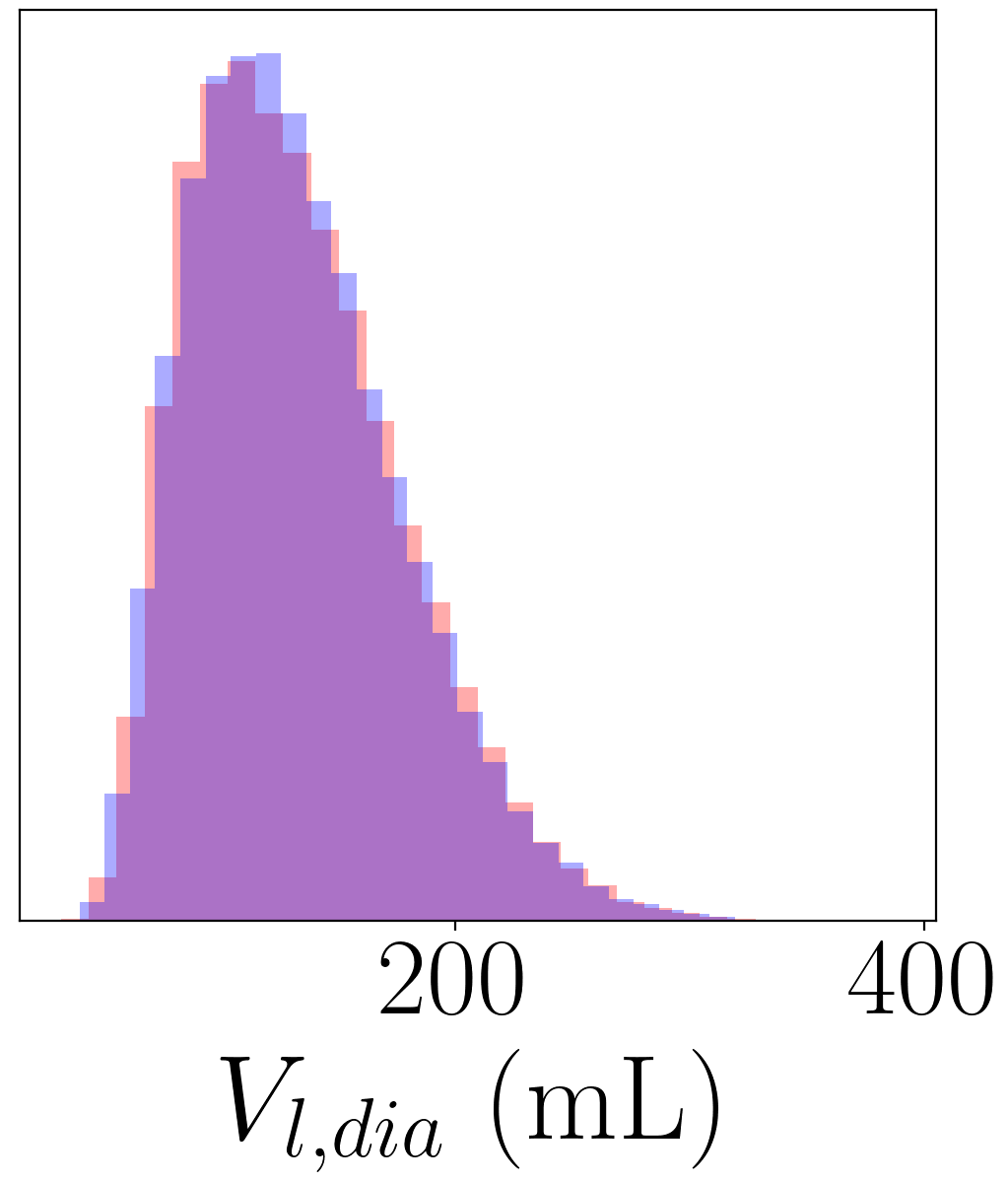}
\end{subfigure}
\begin{subfigure}[b]{0.121\textwidth}
        \centering
         \includegraphics[scale=0.154]{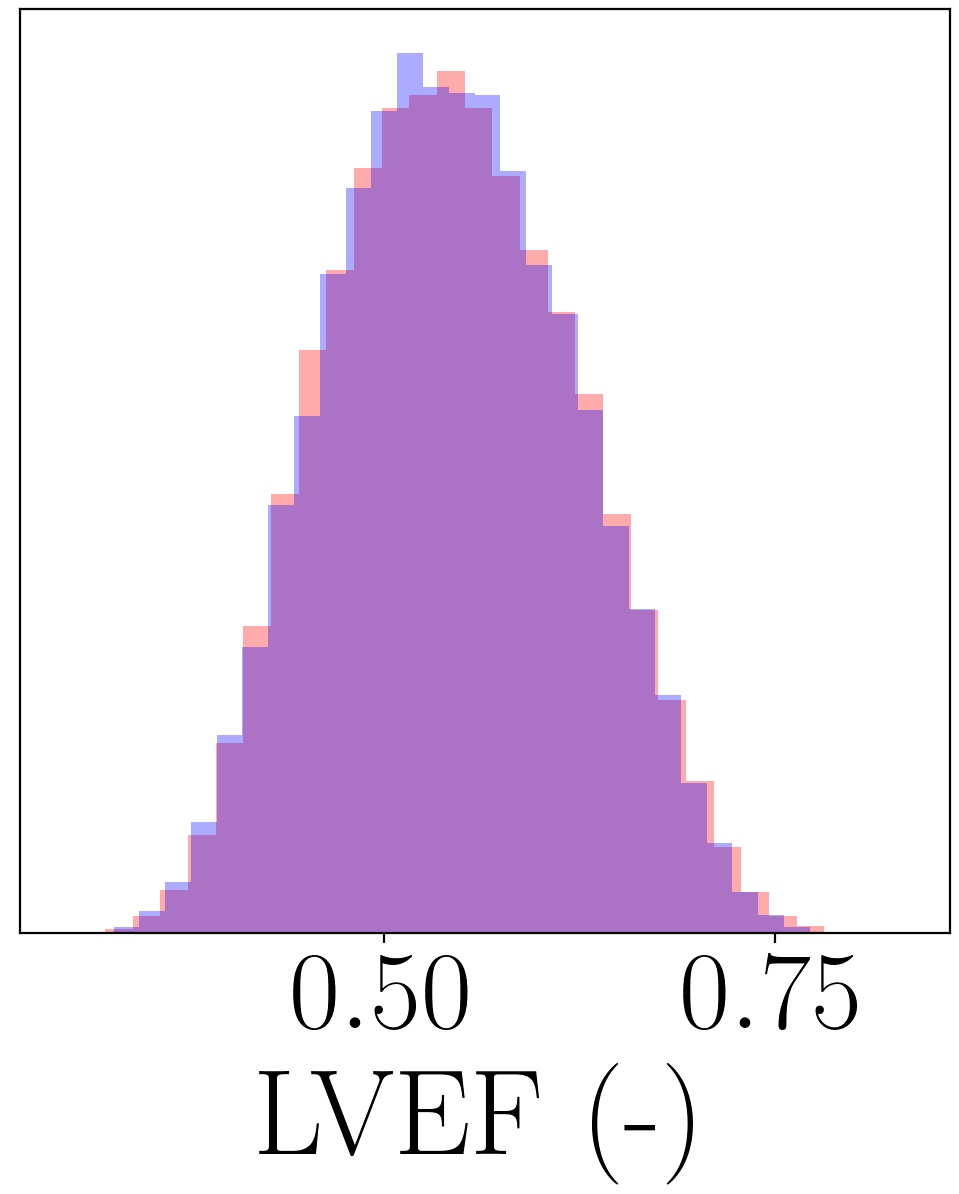}
\end{subfigure}
\begin{subfigure}[b]{0.121\textwidth}
        \centering
         \includegraphics[scale=0.154]{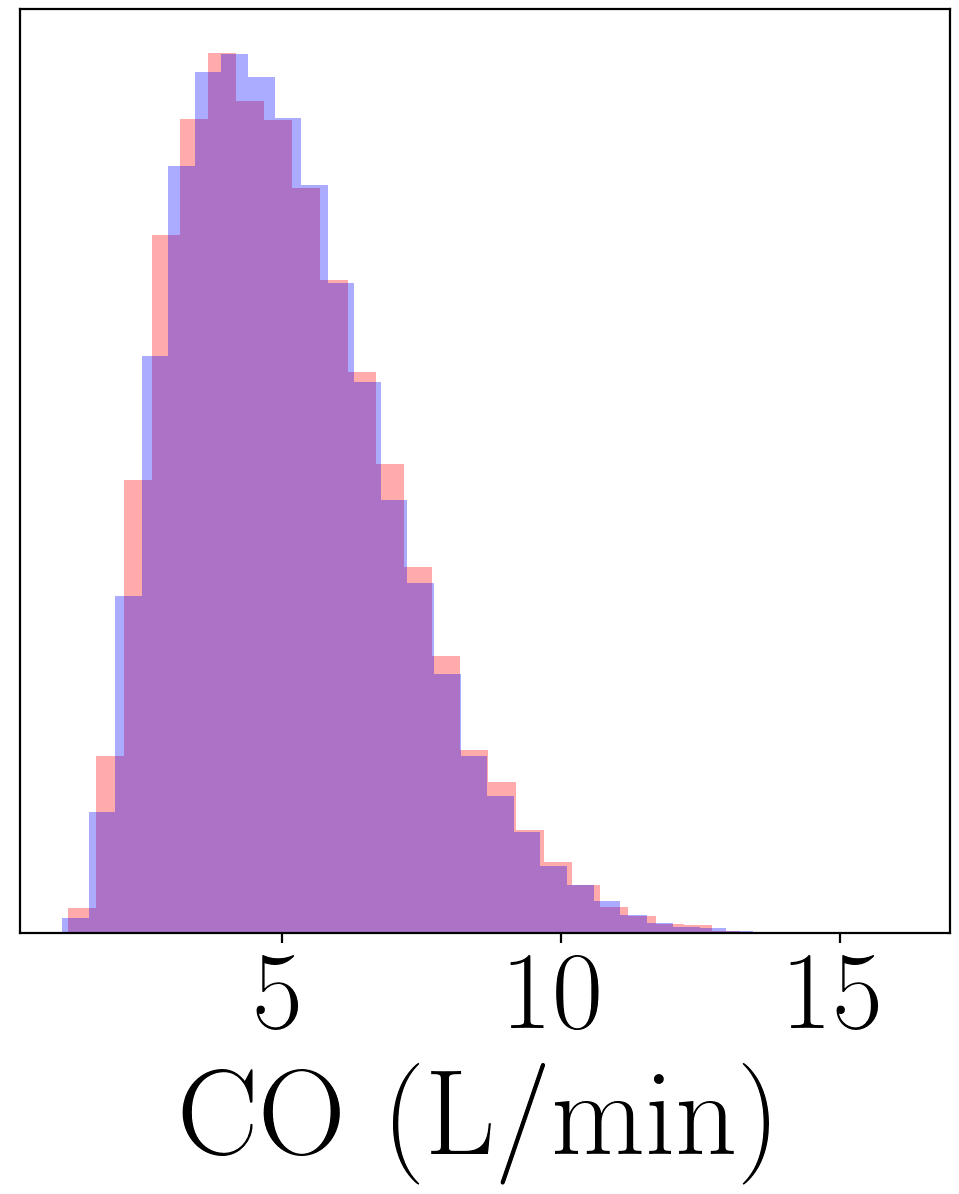}
\end{subfigure}
\begin{subfigure}[b]{0.121\textwidth}
        \centering
         \includegraphics[scale=0.154]{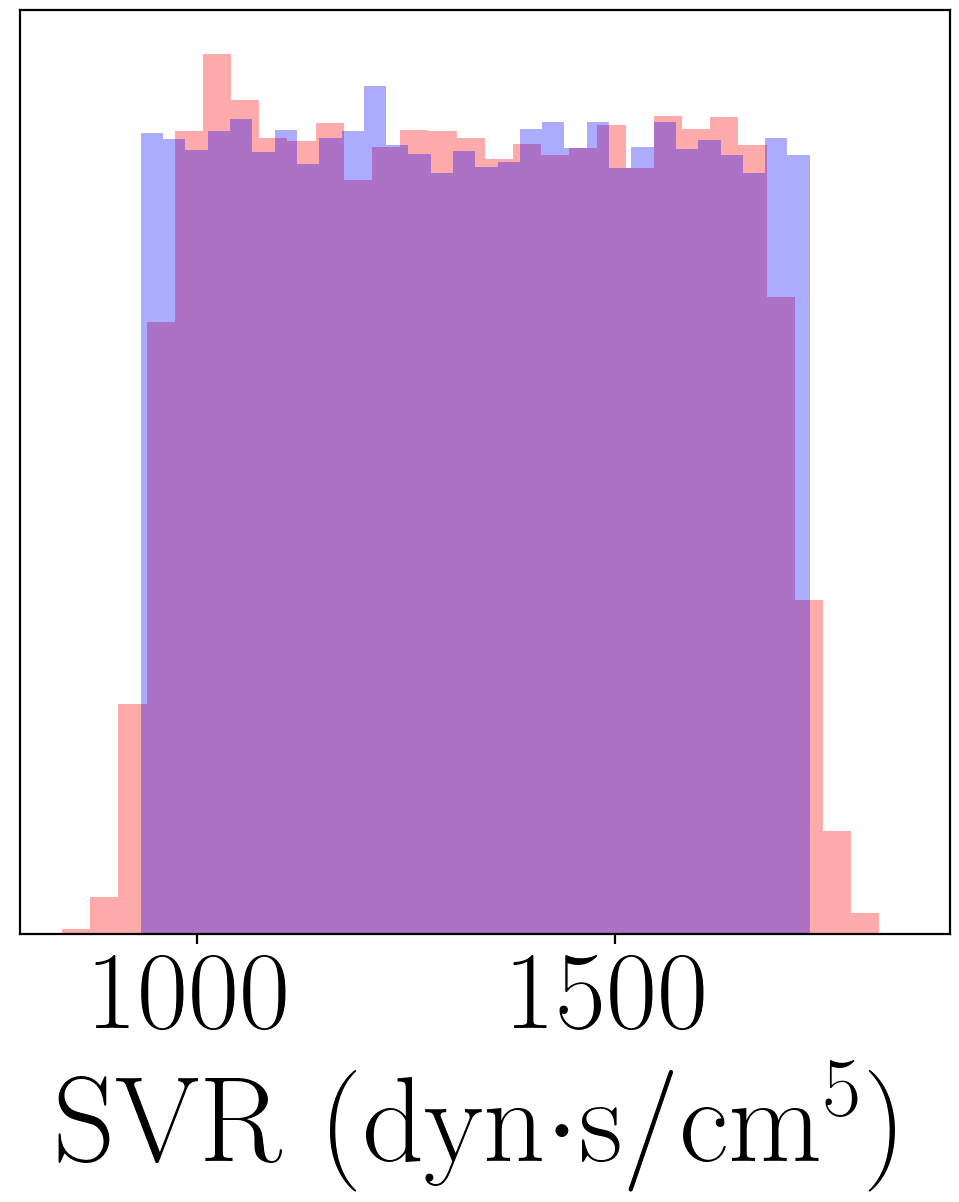}
\end{subfigure}
\begin{subfigure}[b]{0.121\textwidth}
        \centering
         \includegraphics[scale=0.154]{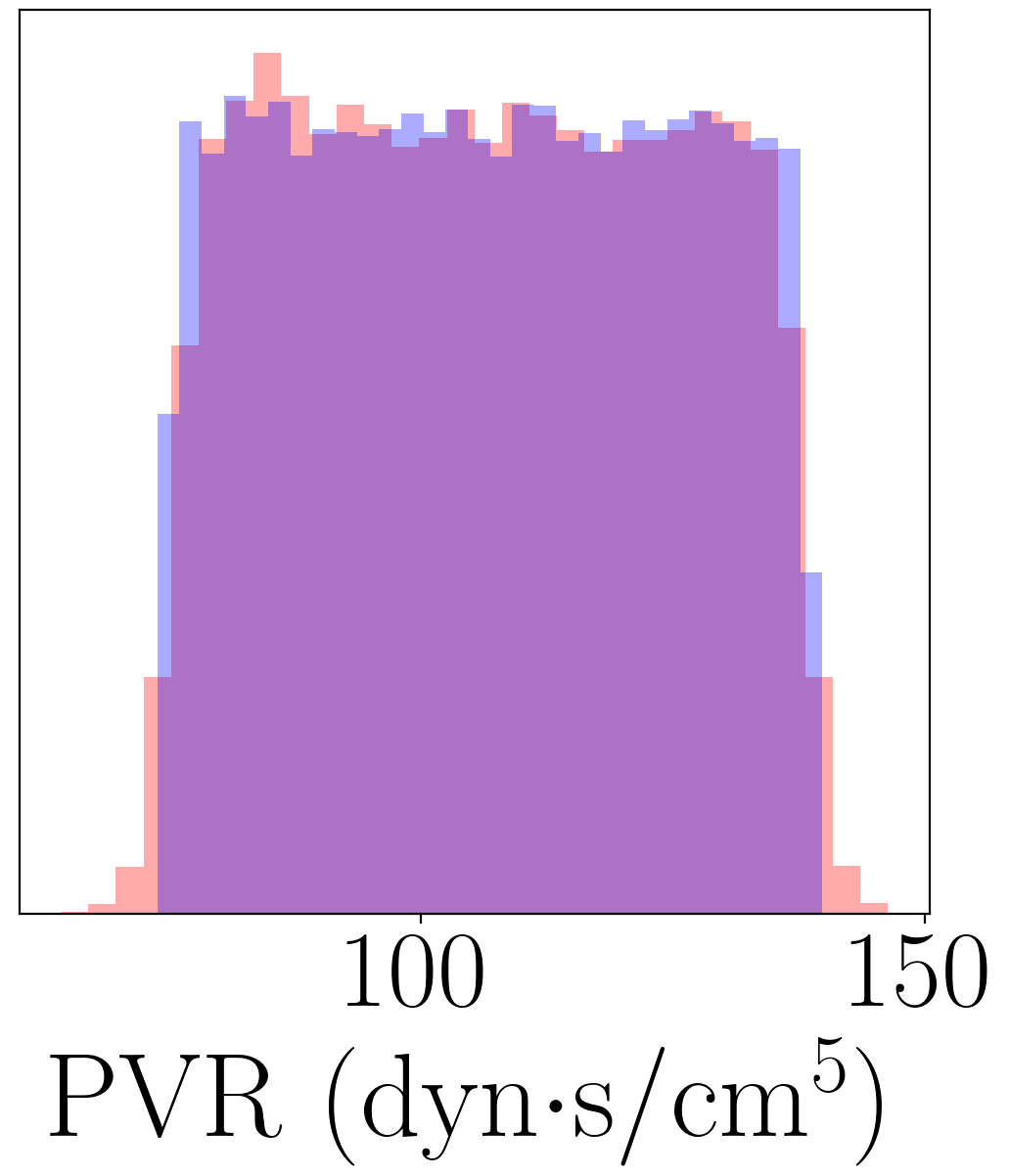}
\end{subfigure}
\caption{Marginal density estimation of the CVSim-6 system outputs. Blue: exact data distribution. Red: predicted distribution via the trained Real-NVP density estimator $NN_f$.}
\label{fig: nf-results}
\end{figure}

Next, Figure~\ref{fig: nf-results} compares the marginal distributions estimated by the trained normalizing flow model $NN_f$ (red), with the exact data distributions (blue), for each output component of the CVSim-6 system.
Overall, the alignment is close, despite some predicted quantity, such as $Hr$, exhibits a slightly heavier tail.

\end{document}